\documentclass[dvips,preprint]{imsart}  

\usepackage{fancyhdr}
\usepackage[T1]{fontenc}
\usepackage[ansinew]{inputenc}
\usepackage{dsfont}
\usepackage{amssymb}
\usepackage{amsfonts}
\usepackage{showlabels}
\usepackage{oldgerm}
\usepackage{multicol}
\usepackage{textcomp}  
\usepackage{amsthm,amsmath}
\usepackage[numbers]{natbib}

\RequirePackage[colorlinks,citecolor=blue,urlcolor=blue]{hyperref}
\usepackage{hypernat}

\startlocaldefs

\newtheorem{defi}{Definition}[section]
\newtheorem{thm}[defi]{Theorem}
\newtheorem{lem}[defi]{Lemma}
\newtheorem{cor}[defi]{Corollary}

\newtheorem{rem}[defi]{Remark}
\newtheorem{rems}[defi]{Remarks}

\newtheorem{prop}[defi]{Proposition}

\newcommand{\pr}{\textsc{Proof:}\\}

\newcommand{\prt}{\textsc{Proof of Theorem }}
\newcommand{\prl}{\textsc{Proof of Lemma }}

\newcommand{\prf}{\textsc{Proof of Formula }}

\newcommand{\Hy}{\mathcal{H}}
\newcommand{\Al}{\mathcal{A}}
\newcommand{\aam}{\alpha^{(m)}_{\Al}}
\newcommand{\ahm}{\alpha^{(m)}_{\Hy}}
\newcommand{\bam}{\beta^{(m)}_{\Al}}
\newcommand{\bhm}{\beta^{(m)}_{\Hy}}

\newcommand{\be}{\vspace{-0.15cm}\begin{equation}}
\newcommand{\ee}{\vspace{-0.15cm}\end{equation}}
\newcommand{\bea}{\vspace{-0.15cm}\begin{eqnarray}}
\newcommand{\eea}{\vspace{-0.15cm}\end{eqnarray}}
\newcommand{\bal}{\beta_{\Al}}
\newcommand{\bhy}{\beta_{\Hy}}
\newcommand{\Pna}{P_{\Al,n}}

\newcommand{\Peinsh}{P_{\Hy,1}}
\newcommand{\Peinsa}{P_{\Al,1}}
\newcommand{\Pnullh}{P_{\Hy,0}}   
\newcommand{\Pnulla}{P_{\Al,0}}   
\newcommand{\Pnh}{P_{\Hy,n}}

\newcommand{\Pa}{P_{\Al}}
\newcommand{\Ph}{P_{\Hy}}
\newcommand{\fa}{f_{\Al}}
\newcommand{\fh}{f_{\Hy}}

\newcommand{\Pmat}{P_{\Al,\left\lfloor \sigma^{2}mt\right\rfloor}^{(m)}}
\newcommand{\Pmht}{P_{\Hy,\left\lfloor \sigma^{2}mt\right\rfloor}^{(m)}}

\newcommand{\oa}{\overline{a}}
\newcommand{\ua}{\underline{a}}
\newcommand{\urho}{\underline{\rho}}
\newcommand{\orho}{\overline{\rho}}

\newcommand{\uu}{\underline{\Upsilon}}
\newcommand{\ou}{\overline{\Upsilon}}

\newcommand{\aal}{\alpha_{\Al}}
\newcommand{\ahy}{\alpha_{\Hy}}
\newcommand{\qua}{\left(\bal,\bhy,\aal,\ahy\right)}     
\newcommand{\qui}{\left(\bal,\bhy,\aal,\ahy,\lambda\right)}     
\newcommand{\quaset}{\mathcal{P}}                               
\newcommand{\quasetNI}{\mathcal{P}_{\textrm{NI}}}                
\newcommand{\quasetSP}{\mathcal{P}_{\textrm{SP}}}               
\newcommand{\quiset}{\mathcal{P} \times ]0,1[}

\newcommand{\quasetSPeins}{\mathcal{P}_{\textrm{SP},1}}               
\newcommand{\quasetSPzwei}{\mathcal{P}_{\textrm{SP},2}}                  
\newcommand{\quasetSPdrei}{\mathcal{P}_{\textrm{SP,3}}}                 
\newcommand{\quasetSPdreiab}{\mathcal{P}_{\textrm{SP,3ab}}}

\newcommand{\quisetSPdreia}{\mathcal{P}_{\textrm{SP,3a}}^{\lambda,\leq 0}}               
\newcommand{\quisetSPdreib}{\mathcal{P}_{\textrm{SP,3b}}^{\lambda,> 0}}               
\newcommand{\quasetSPdreic}{\mathcal{P}_{\textrm{SP,3c}}}                
\newcommand{\quasetSPdreid}{\mathcal{P}_{\textrm{SP,3d}}}                
\newcommand{\quasetSPvier}{\mathcal{P}_{\textrm{SP,4}}}                
\newcommand{\quasetSPcomp}{(\mathcal{P}_{\textrm{SP}}\backslash \mathcal{P}_{\textrm{SP,1}})}     
\newcommand{\quasetSPcompvar}{\mathcal{P}_{\textrm{SP}}\backslash \mathcal{P}_{\textrm{SP,1}}}     

\newcommand{\largeint}{\overline{\mathbb{N}}}

\endlocaldefs

\begin{document}

\begin{frontmatter}

\title{Some distance bounds of branching processes and their diffusion limits}
\runtitle{Distances of branching processes}

\begin{aug}
\author{\fnms{Niels B.} \snm{Kammerer}\thanksref{t1}\ead[label=e1]{niels.kammerer@wiso.uni-erlangen.de}} 
\and  
\author{\fnms{Wolfgang} \snm{Stummer}
\corref{}\ead[label=e2]{stummer@mi.uni-erlangen.de}
}
\affiliation{University of Erlangen--N\"urnberg}

\address{Department of Mathematics\\
University of Erlangen--N\"urnberg\\
Bismarckstrasse $1 \frac{1}{2}$\\
91054 Erlangen, Germany\\
\printead{e1}\\[-0.32cm]
\phantom{E-mail:\ }\printead*{e2}\\
}

\thankstext{t1}{Supported by the ''Studienstiftung des deutschen Volkes''.}
\runauthor{Niels B. Kammerer and Wolfgang Stummer}

\end{aug}


\begin{abstract}
We compute exact values respectively bounds of  
``distances'' -- 
in the sense of (transforms of) power divergences and relative entropy -- between 
two discrete-time 
\emph{Galton-Watson branching processes with immigration} GWI for which
the offspring as well as the immigration is arbitrarily Poisson-distributed
(leading to arbitrary type of criticality).
Implications for asymptotic distinguishability behaviour in terms of contiguity and entire separation of the involved GWI are given, too.
Furthermore, we determine the corresponding limit quantities for the context in which the two GWI converge  
to Feller-type 
branching diffusion processes, as the time-lags between observations tend to zero. 
Some applications to (static random environment like) Bayesian decision making and Neyman-Pearson testing are presented as well.
\end{abstract}



\end{frontmatter}

\section{Introduction}
\label{sec.1}

It is well known that ``distances'' in form of (relative-entropy covering) \emph{power divergences} 
between finite measures are 
important for probability theory and statistics as well as their applications to various different research fields 
such as  physics, information theory, econometrics, biology, 
speech and image recognition, 
transportation of (sorts of) ``mass'', etc. For probability measures $P_{\Hy}$, $P_{\Al}$
on a measurable space $\left(\Omega,\mathcal{A}\right)$  
 and
parameter $\lambda \in \mathbb{R}$ these power divergences -- also known as Cressie-Read measures
respectively generalized cross-entropy family -- are defined as (see e.g.\ Liese and Vajda \cite{LieVa87}, \cite{LieVa06})
\bea  
I_{\lambda}\left(P_{\Al}||P_{\Hy}\right) \ := \ 
\left\{  \begin{array}{cl}
I\left(P_{\Al}||P_{\Hy}\right),  &   \hskip1cm {\rm if} \ \ 
\lambda = 1, \\[0.2cm]
\frac{1}{\lambda (\lambda-1)} \left( H_{\lambda}\left(P_{\Al}||P_{\Hy}\right) -1 \right),  & \hskip1cm {\rm if} \ \ \lambda 
\in \mathbb{R} \backslash {\{}0,1{\}}, \nonumber \\[0.2cm]
I\left(P_{\Hy}||P_{\Al}\right),  &  \hskip1cm {\rm if} \ \ 
\lambda = 0, \\[0.2cm]
\end{array}  \right.\\[-0.82cm]
\label{fo.powerdivdef}
\eea
\noindent
where
 \vspace{-0.2cm}
 \be
I\left(P_{\Al}||P_{\Hy}\right) \ := \ \int_{\{p_{\Hy}>0\}} p_{\Al} \, \log\frac{p_{\Al}}{p_{\Hy}} \ d\mu \ + \ \infty \cdot P_{\Al}(p_{\Hy}=0)
\label{fo.relent}
\ee
is the relative entropy (Kullback-Leibler information divergence) and
 \vspace{-0.1cm}
\be
H_{\lambda}\left(P_{\Al}||P_{\Hy}\right) \ := \ \int_{\Omega} p_{\Al}^{\lambda} \, p_{\Hy}^{1-\lambda} \ d\mu 
\label{fo.hellallg}
\ee
is the Hellinger integral of order $\lambda \in \mathbb{R}\backslash {\{}0,1{\}}$; for this, we assume as usual without loss of generality
that the probability measures $P_{\Hy}$, $P_{\Al}$
are dominated by some $\sigma-$finite measure $\mu$, with densities 
\be
p_{\Al} \:= \ \frac{\mathrm{d}P_{\Al}}{\mathrm{d}\mu} \qquad \mathrm{and} \qquad 
p_{\Hy} \:= \ \frac{\mathrm{d}P_{\Hy}}{\mathrm{d}\mu}
\nonumber
\ee
defined on $\Omega$ (the zeros of $p_{\Hy}$, $p_{\Al}$ are handled in \eqref{fo.relent}, \eqref{fo.hellallg} with the usual conventions). 
Apart from the relative entropy, other prominent examples of power divergences are the squared
Hellinger distance $\frac{1}{2} \, I_{1/2}\left(P_{\Al}||P_{\Hy}\right)$
and  Pearson's $\chi^2-$divergence $2\, I_{2}\left(P_{\Al}||P_{\Hy}\right)$.
Extensive studies about basic and advanced general facts on power divergences, Hellinger integrals and
the related Renyi divergences of order $\lambda \in \mathbb{R}\backslash {\{}0,1{\}}$ 
\be
R_{\lambda}\left(P_{\Al}||P_{\Hy}\right) \ := \ \frac{1}{\lambda (\lambda-1)} \, \log H_{\lambda}\left(P_{\Al}||P_{\Hy}\right) \ ,  
\qquad  \textrm{with } \log\, 0 = - \infty,
\notag
\ee 
can be found e.g.\ in Liese and Vajda \cite{LieVa87}, \cite{LieVa06}, Jacod and Shiryaev \cite{JacSh87}. For instance, the integrals in \eqref{fo.relent} and \eqref{fo.hellallg}
do not depend on the choice of $\mu$. As far as finiteness is concerned, for $\lambda \in ]0,1[$ one gets the rudimentary bounds
\be
0 \ \leq \ I_{\lambda}\left(P_{\Al}||P_{\Hy}\right) \ \leq \frac{1}{\lambda (\lambda-1)} \ ,
\label{fo.powbound}
\ee
where the lower bound is achieved if and only if $P_{\Al}=P_{\Hy}$, and the upper bound
is achieved if and only if $P_{\Al} \bot P_{\Hy}$ (singularity).
For $\lambda \notin ]0,1[$, the power divergences $I_{\lambda}\left(P_{\Al}||P_{\Hy}\right)$ 
and 
Hellinger integrals $H_{\lambda}\left(P_{\Al}||P_{\Hy}\right)$
might be infinite,
depending on the particular setup. For the sake of brevity, we only deal here with the case
$\lambda \in [0,1]$; the case $\lambda \notin [0,1]$ will appear elsewhere. 

\bigskip
\noindent 
Apart from the extensive literature on the relative-entropy cases $\lambda (1-\lambda) = 0$, for $\lambda (1-\lambda) \ne 0$ 
the evaluation of \emph{power divergences $I_{\lambda}$} -- respectively their straightforward transforms 
such as Hellinger integrals $H_{\lambda}$ and Renyi divergences $R_{\lambda}$ -- 
have been investigated for various different contexts of (probability distributions associated with)
\emph{stochastic processes}, such as 
processes with independent increments (see e.g.\ Newman  
\cite{New73}, Liese \cite{Lie82}, Memin and Shiryaev \cite{Mem85}, Jacod and Shiryaev \cite{JacSh87}, Liese and Vajda \cite{LieVa87},
Linkov and Shevlyakov \cite{LinkSh98}), 
 Poisson point processes (see e.g.\  
 Liese \cite{Lie85}, 
 Jacod and Shiryaev \cite{JacSh87}, Liese and Vajda \cite{LieVa87}), 
diffusion prcoesses respectively solutions of stochastic differential equations with continuous paths (see e.g.\ 
Kabanov et al.\ \cite{Kab86}, Liese \cite{Lie86}, 
Jacod and Shiryaev \cite{JacSh87}, Liese and Vajda \cite{LieVa87}, Vajda \cite{Vaj90}, Stummer \cite{Stu99}, \cite{Stu01}, 
Stummer and Vajda \cite{StuVa07}); further related literature can be found e.g.\ in references of the abovementioned papers and books.

\bigskip
\noindent
Another important 
class of time-dynamic
models is given by discrete-time branching processes,
in particluar Galton-Watson processes without immigration GW respectively with immigration GWI, which have numerous applications in 
biotechnology, population genetics, internet traffic research, clinical trials, asset price modelling and derivative pricing.
(Transforms of) Power divergences have been used 
for supercritical Galton-Watson processes without immigration SUPGW 
for instance as follows:
Feigin and Passy \cite{Fei81} study the problem to find
an offspring distribution which is closest (in terms of 
relative entropy type distance) to the original offspring distribution
and under which ultimate extinction is certain. Furthermore, 
Mordecki \cite{Mor94} gives an equivalent characterization for the stable convergence of the corresponding 
log-likelihood process to a mixed Gaussian limit, in terms of conditions on Hellinger integrals 
of the involved offspring laws. 
Moreover, 
Sriram and Vidyashankar \cite{Sri00} study the properties of offspring-distribution-parameters
which minimize the squared Hellinger distance $\frac{1}{2} \, I_{1/2}$ between
the model offspring distribution and the corresponding non-parametric maximum likelihood estimator
of Guttorp \cite{Gut91}.
For
the 
setup of GWI with Poisson offspring and nonstochastic immigration of constant value 1,
Linkov and Lunyova \cite{Link96}
investigate the asymptotics of Hellinger integrals 
in order to deduce large deviation assertions in hypotheses testing problems.

\bigskip
\noindent
In contrast to the abovementioned contexts, this paper pursues the 
following main goals:

\begin{enumerate}

\item[(MG1)] for any time horizon and any criticality scenario, to compute (non-rudimentary) lower and upper bounds 
-- and sometimes even exact values -- of the Hellinger integrals 
$H_{\lambda}\left(P_{\Al}||P_{\Hy}\right)$
and power divergences $I_{\lambda}\left(P_{\Al}||P_{\Hy}\right)$ \, ($\lambda \in [0,1]$) of  two
Galton-Watson branching processes 
$P_{\Al}$, $P_{\Hy}$ with Poisson($\bal$) respectively Poisson($\bhy$) distributed offspring 
as well as Poisson($\aal$) respectively Poisson($\ahy$) distributed immigration.
As a side effect, we also aim for corresponding 
asymptotic distinguishability results in terms of contiguity and entire separation.

\item[(MG2)] to compute the corresponding limit quantities for the context in which 
(a proper rescalation of) the two
Galton-Watson processes with immigration converge  
to \emph{Feller}-type branching diffusion processes, 
as the time-lags between the generation-size observations tend to zero.

\item[(MG3)]  as an exemplary field of application, to 
indicate how to
use the results of (MG1)  
for Bayesian decision making and Neyman-Pearson testing
based on the sample path observations of the GWI-generation sizes, when the hypothesis law
is given by $P_{\Hy}$ and the alternative law by  $P_{\Al}$; in a certain sense,
this can also be 
interpreted in terms of
a rudimentary static random environment. 

\end{enumerate}

\noindent
Because of the involved Poisson distributions, these goals 
(which are potentially reasonable also for other types of offspring resp.\ immigration distributions)
can be tackled with a high degree of 
tractability, 
which is worked out in detail with the following structure:
we first deal with the non-relative-entropy
case $\lambda (1-\lambda) \ne 0$. Section \ref{secSETUP} contains the first basic 
result concerning Goal (MG1),
which is then deepened in Section \ref{secDET} in order to obtain -- parameter constellation dependent -- 
\emph{recursively computable} exact values respectively 
\emph{recursively computable} lower and upper bounds of $H_{\lambda}\left(P_{\Al}||P_{\Hy}\right)$.
Additionally, we construct related \emph{closed-form} bounds in Section \ref{secCFB},
which will also be used to achieve (the Hellinger-integral part of) Goal (MG2) in Section \ref{sec.diflim}.
The power divergences $I_{\lambda}\left(P_{\Al}||P_{\Hy}\right)$ are treated in Section \ref{sec.ent},
complemented with the relative-entropy cases $\lambda (1-\lambda) = 0$ of the Goals (MG1), (MG2).
The subsequent Section \ref{sec.dec} is concerned with Goal (MG3), whereas the Appendix contains
main proofs and auxiliary lemmas.




\section{Process setup and first basic result}
\label{secSETUP}

\noindent
Let $X_{n}$ denote the $n$th generation size of a discrete-time Galton-Watson process
with immigration GWI. We use the recursive description 
 \vspace{-0.1cm}
\be\label{033}
X_{0} := \omega_{0} \in \mathbb{N}; \qquad
X_{n}=\sum_{k=1}^{X_{n-1}}Y_{n-1,k}~+~\widetilde{Y}_{n}, \qquad n \in \mathbb{N}, 
 \vspace{-0.1cm}
\ee
where $Y_{n-1,k}$ is the number of offspring of the $k$th object (e.g.\ organism, person) 
within the $(n-1)$th generation, and  
$\widetilde{Y}_{n}$ denotes the number of immigrating objects in the $n$th generation. Notice that we employ an arbitrary deterministic initial generation size $X_{0}$.
We always assume that under the law $P_{\Hy}$ (e.g.\ a hypothesis), 

\begin{itemize}

\item the collection $Y:= \left\{ Y_{n-1,k}, \, n  \in \mathbb{N}, k  \in \mathbb{N} \right\}$
consists of independent and identically distributed (i.i.d.) random variables which are
Poisson distributed with parameter $\bhy >0$,

\item the collection $\widetilde{Y}:= \left\{ \widetilde{Y}_{n}, \, n  \in \mathbb{N} \right\}$
consists of i.i.d. random variables which are
Poisson distributed with parameter $\ahy \geq 0$ (where $\ahy = 0$ stands for the degenerate
case of having no immigration),

\item $Y$ and $\widetilde{Y}$ are independent.

\end{itemize}

\noindent
In contrast, under the law $P_{\Al}$ (e.g.\ an alternative) the same is supposed to hold with parameters $\bal >0$ 
(instead of $\bhy >0$) and $\aal \geq 0$ (instead of $\ahy \geq 0$). 
Furthermore, let $(\mathcal{F}_{n})_{n \in \mathbb{N}}$ be the 
corresponding canonical filtration generated by $X :=(X_{n})_{n \in \mathbb{N}}$.

\bigskip
\noindent
Basic and advanced facts on GWI (introduced by Heathcote \cite{Hae65}) can be found e.g.\ in the monographs of 
Athreya and Ney \cite{AN72}, Jagers \cite{Jag75},
Asmussen and Hering \cite{As83}, Haccou \cite{Hac05}; see also e.g.\ Heyde and Seneta \cite{HeySen72}, \cite{HeySen72}, Basawa and Rao \cite{Bas80}, 
Basawa and Scott \cite{Bas83}, Sankaranarayanan \cite{San89}, Wei and Winnicki \cite{Wei90}, Winnicki \cite{Win91},
Guttorp \cite{Gut91} as well as 
Yanev \cite{Yan08} 
(and also the references therein all those) for 
adjacent
fundamental statistical issues including 
the involved technical respectively conceptual challenges.

\bigskip
\noindent
For the sake of brevity, wherever we introduce or discuss corresponding quantities \emph{simultaneously} for both the hypothesis $\Hy$ and the alternative $\Al$, 
we will use the subscript $\bullet$ as a synonym for either the symbol 
$\Hy$ or  
$\Al$.
For illustration, recall the well-known fact that the corresponding conditional probabilities $P_{\bullet}(X_{n}=\cdot~|X_{n-1}=k)$ 
are again Poisson-distributed, with parameter $\beta_{\bullet}\cdot k+\alpha_{\bullet}$.
In oder to achieve a 
transparently representable
structure of our results, we subsume 
the involved parameters as follows: let $\quasetSP$ be the set of all constellations $\qua$ of real-valued
parameters $\bal >0$, $\bhy >0$, $\aal > 0$, $\ahy>0$, such that $\bal \ne \bhy$
or $\aal \ne \ahy$ (or both). Furthermore, we write $\quasetNI$ for the
set of all $\qua$ of real-valued parameters $\bal >0$, $\bhy >0$, $\aal = \ahy =0$, such that
$\bal \ne \bhy$; this corresponds to the important special case of having no immigration. 
The resulting disjoint union will be denoted by $\quaset = \quasetSP \cup \quasetNI$.
A typical situation for applications in our mind is that one particular constellation $\qua \in \quaset$ 
(e.g.\ obtained from theoretical or previous statistical investigations)
is 
fixed, whereas -- in contrast -- the parameter $\lambda \in 
]0,1[$ for the 
Hellinger integral or the power divergence might be chosen freely, e.g.\ depending
on which ``probability distance''  
one decides to choose for further
analysis. At this point, let us emphasize that \emph{in general} 
we will not make assumptions of the form $\beta_{\bullet} \gtreqqless 1$, i.e.\ upon the type of criticality.\\[0.1cm]

\noindent
To start with our investigations, we define the extinction time $\tau:= \min\{l \in \mathbb{N}: X_{m}=0$  for all integers  
$m \geq l\}$ if this minimum exists, and $\tau:=\infty$ else. Correspondingly, let $\mathcal{B} := \{\tau < \infty\}$
be the extinction set. It is well known that in the case $\quasetNI$
one gets $P_{\bullet}(\mathcal{B})=1$ \, if $0<\beta_{\bullet} \leq 1$ \, and \, 
$P_{\bullet}(\mathcal{B}) \in \, ]0,1[$ \, if $\beta_{\bullet} > 1$. In contrast, for 
$\quasetSP$ there always holds $P_{\bullet}(\mathcal{B})=0$.
Furthermore, for $\quasetSP$ the two laws $P_{\Hy}$ and $P_{\Al}$ are equivalent, whereas
for $\quasetNI$ the two restrictions  $\left.\Ph\right|_{\mathcal{B}}$ and  $\left.\Pa\right|_{\mathcal{B}}$ 
are equivalent 
(see e.g.\ Lemma 1.1.3 of Guttorp \cite{Gut91}); with a slight abuse of noation we shall henceforth omit
$\left.\, \right|_{\mathcal{B}}$ \, . Consistently, for fixed time $n \in \mathbb{N}_{0}$ 
we introduce $\Pna := \left.\Pa\right|_{\mathcal{F}_{n}}$
and $\Pnh := \left.\Ph\right|_{\mathcal{F}_{n}}$ as well as
the corresponding Radon-Nikodym-derivative 
\be\label{def.Zn}
Z_{n}~:=~\frac{\textrm{d}\Pna}{\textrm{d}\Pnh} \ .
\ee
Clearly, $Z_{0}=1$.
By using the ``rate functions'' $f_{\bullet}(x)=\beta_{\bullet}\, x+\alpha_{\bullet}$ ($x\in[0,\infty[$), 
a version of \eqref{def.Zn} can be easily determined
by calculating for each $\omega=(\omega_{0},...,\omega_{n})\in\Omega_{n}:= \mathbb{N}_{0}^{n}$ 
\bea 
Z_{n}(\omega) = 
 \prod_{k=1}^{n}Z_{n,k}(\omega)  & & \text{with}~Z_{n,k}(\omega):=
 \exp\Big\{-\big(f_{\Al}(\omega_{k-1})-f_{\Hy}(\omega_{k-1})\big)\Big\}\left[\frac{f_{\Al}(\omega_{k-1})}{f_{\Hy}(\omega_{k-1})}\right]^{\omega_{k}},
 \notag
\eea
where for the last term we use the convention $\left(\frac{0}{0}\right)^{x}=1$ for all $x\in \mathbb{N}_{0}$.
\noindent Furthermore, we define for each $\omega\in\Omega_{n}
$ 
 \vspace{-0.2cm}
\be
Z^{(\lambda)}_{n,k}(\omega):=\exp\Big\{-\big(\lambda f_{\Al}(\omega_{k-1})+(1-\lambda)f_{\Hy}(\omega_{k-1})\big)\Big\}
\ \frac{\left[\left(f_{\Al}(\omega_{k-1})\right)^{\lambda}\left(f_{\Hy}(\omega_{k-1})\right)^{1-\lambda}\right]^{\omega_{k}}}{\omega_{k}!} 
\label{fo.Znk}
\ee
with the convention $\frac{\left(0\right)^{0}}{0!}=1$ for the last term.
Accordingly, with the choice $\mu = \Pnh$ one obtains from \eqref{fo.hellallg}
the Hellinger integral $H_{\lambda}\left(\Pnulla||\Pnullh\right) =1$, as well as  
for all $\qui \in \quiset$    
\be
H_{\lambda}\left(\Peinsa||\Peinsh\right) = 
\exp\Big\{\left(f_{\Al}(\omega_{0})\right)^{\lambda}\left(f_{\Hy}(\omega_{0})\right)^{(1-\lambda)}-(\lambda f_{\Al}(\omega_{0})+(1-\lambda)f_{\Hy}(\omega_{0}))\Big\}
\label{fo.heins}
\ee
and for all $n \in \mathbb{N}\backslash\{1\}$
\bea
&& \hspace{-0.5cm} H_{\lambda}\left(\Pna||\Pnh\right)~=~E\Pnh\big[(Z_{n})^{\lambda}\big]
~=~\sum_{\omega_{1}=0}^{\infty}\cdots\sum_{\omega_{n}=0}^{\infty}\prod_{k=1}^{n}Z^{(\lambda)}_{n,k}
(\omega) \notag\\
&&  \hspace{-0.5cm}= \ \sum_{\omega_{1}=0}^{\infty}\cdots\sum_{\omega_{n-1}=0}^{\infty}\prod_{k=1}^{n-1}Z^{(\lambda)}_{n,k}
(\omega)\cdot e^{-(\lambda f_{\Al}(\omega_{n-1})+(1-\lambda)f_{\Hy}(\omega_{n-1}))}\sum_{\omega_{n}=0}^{\infty}\frac{\left[\left(f_{\Al}(\omega_{n-1})\right)^{\lambda}\left(f_{\Hy}(\omega_{n-1})\right)^{1-\lambda}\right]^{\omega_{n}}}{\omega_{n}!}
\notag\\
&&  \ = \ \sum_{\omega_{1}=0}^{\infty}\cdots\sum_{\omega_{n-1}=0}^{\infty}\prod_{k=1}^{n-1}Z^{(\lambda)}_{n,k}
(\omega)\cdot e^{\left(f_{\Al}(\omega_{n-1})\right)^{\lambda}\left(f_{\Hy}(\omega_{n-1})\right)^{1-\lambda}-(\lambda f_{\Al}(\omega_{n-1})+(1-\lambda)f_{\Hy}(\omega_{n-1}))} \ .\label{fo.hellzw}
\eea
From \eqref{fo.hellzw}, one can see that a crucial role for the exact calculation (respectively the derivation of bounds)  
of the Hellinger integral is played by the functions defined for $x \in[0,\infty[$
\be\label{defphi}
\phi_{\lambda}(x)~:=~\varphi_{\lambda}(x)-f_{\lambda}(x)\ , \qquad \textrm{with} 
\vspace{-0.3cm}
\ee
\be\label{defvarphi}
\varphi_{\lambda}(x)~:=~\left(f_{\Al}(x)\right)^{\lambda}\left(f_{\Hy}(x)\right)^{1-\lambda} \qquad \textrm{and} 
\vspace{-0.2cm}
\ee
\be\label{defflambda}
f_{\lambda}(x)~:=~\lambda f_{\Al}(x)+(1-\lambda)f_{\Hy}(x)~=~ \alpha_{\lambda} + \beta_{\lambda} \, x \ ,
\ee
where we have used the $\lambda$-\emph{weighted-averages} $\beta_{\lambda}=\lambda\cdot\beta_{\Al}+(1-\lambda)\cdot\beta_{\Hy}$ and $\alpha_{\lambda}=\lambda\cdot\alpha_{\Al}+(1-\lambda)\cdot\alpha_{\Hy}$.
According to Lemma \ref{lem2} in Appendix \ref{App3}, it follows for $\lambda\in ]0,1[$ 
that $\phi_{\lambda}(x)\leq 0$ for all $x \in[0,\infty[$, and that $\phi(x)=0$ iff $\fa(x)=\fh(x)$. This is consistent
with the corresponding generally valid upper bound 
 \vspace{-0.1cm}
\be
H_{\lambda}\left(\Pna||\Pnh\right) \leq 1 \ .
\label{fo.hellgenup}
\ee

\noindent As a first indication for our proposed method, let us start by illuminating
the simplest case $\lambda\in ]0,1[$ and $\gamma:= 
\alpha_{\Hy}\beta_{\Al}-\alpha_{\Al}\beta_{\Hy}=0$. This means that $\qua \in \quasetNI \cup \quasetSPeins$,
where $\quasetSPeins$ is the set of all (componentwise) strictly positive $\qua$ 
with $\bal \ne \bhy$, $\aal \ne \ahy$ and $\frac{\bal}{\bhy} = \frac{\aal}{\ahy} \ne 1$.
In this situation, \emph{all} the three functions \eqref{defphi} to \eqref{defflambda} are linear. Indeed,
\be
\varphi_{\lambda}(x) \ = \ p_{\lambda}^{E} + q_{\lambda}^{E} \, x
\notag
\ee
with $p_{\lambda}^{E}:= \aal^{\lambda}\, \ahy^{1-\lambda}$ and $q_{\lambda}^{E}:= \bal^{\lambda}\, \bhy^{1-\lambda}$ (where the index E stands for \underline{e}xact linearity).
Clearly, $q_{\lambda}^{E} >0$ on $\quasetNI \cup \quasetSPeins$,
as well as $p_{\lambda}^{E} >0$ on $\quasetSPeins$ respectively $p_{\lambda}^{E} =0$ on $\quasetNI$.
Furthermore,
\be
\phi_{\lambda}(x) \ = \ 
r_{\lambda}^{E} + s_{\lambda}^{E} \, x
\notag
\ee
with $r_{\lambda}^{E}:=p_{\lambda}^{E} - \alpha_{\lambda} = 
\aal^{\lambda}\, \ahy^{1-\lambda} - (\lambda \aal +(1-\lambda) \ahy)$ and  
$s_{\lambda}^{E} := q_{\lambda}^{E}- \beta_{\lambda} = 
\bal^{\lambda}\, \bhy^{1-\lambda} - (\lambda \bal +(1-\lambda) \bhy)$.
Due to Lemma \ref{lem2} one knows $s_{\lambda}^{E} <0$ on $\quasetNI \cup \quasetSPeins$ \, ,
as well as $r_{\lambda}^{E} <0$ on $\quasetSPeins$ respectively $r_{\lambda}^{E} =0$ on $\quasetNI$.\\[0.1cm]

\noindent
As it will be seen later on, such kind of linearity properties are useful for the recursive handling of the Hellinger integrals. 
However, only on the parameter set $\quasetNI \cup \quasetSPeins$ the functions
$\varphi_{\lambda}$ and $\phi_{\lambda}$ are linear. Hence, in the general case $\qui \in \quaset \times ]0,1[$ we aim for linear lower and upper bounds
\be
\varphi_{\lambda}^{L}(x) := p_{\lambda}^{L} + q_{\lambda}^{L} \, x \ \leq \ 
\varphi_{\lambda}(x) 
\ \leq \ \varphi_{\lambda}^{U}(x) := p_{\lambda}^{U} + q_{\lambda}^{U} \, x \ ,
\label{fo.varphibou}
\ee
\hspace{-0.8cm}  $\qquad x\in[0,\infty[$ (ultimately, $x\in \mathbb{N}_{0}$), which lead to
\be
\phi_{\lambda}^{L}(x)   :=   r_{\lambda}^{L} + s_{\lambda}^{L}  x 
 :=  (p_{\lambda}^{L} - \alpha_{\lambda}) + (q_{\lambda}^{L} - \beta_{\lambda})  x
 \leq  
\phi_{\lambda}(x) \leq  
\phi_{\lambda}^{U}(x)  := r_{\lambda}^{U} + s_{\lambda}^{U}  x 
 :=  (p_{\lambda}^{U} - \alpha_{\lambda}) + (q_{\lambda}^{U} - \beta_{\lambda})  x  , 
\label{fo.phibou}
\ee
\hspace{-0.8cm} $\qquad x\in[0,\infty[$ (ultimately, $x\in \mathbb{N}_{0}$). Of course, the involved 
slopes and intercepts should satisfy reasonable restrictions. 
For instance, because of the nonnegativity of $\varphi_{\lambda}$ we require $p_{\lambda}^{U} \geq p_{\lambda}^{L} \geq 0$, 
$q_{\lambda}^{U} \geq q_{\lambda}^{L} \geq 0$
(leading to the nonnegativity of $\varphi_{\lambda}^{L}$, $\varphi_{\lambda}^{U}$). 
Furthermore, \eqref{fo.hellzw} and \eqref{fo.hellgenup} suggest that
$p_{\lambda}^{L} \leq \alpha_{\lambda}$, $q_{\lambda}^{L} \leq \beta_{\lambda}$
which leads to the nonpositivity of $\phi_{\lambda}^{L}$.
Moreover, it is assumed that
\be
\textit{at least one} \textrm{ of the two inequalities }  p_{\lambda}^{U} < \alpha_{\lambda}, \
q_{\lambda}^{U} < \beta_{\lambda}  \textrm{ holds,}
\label{fo.twoin} 
\ee
and hence $\phi_{\lambda}^{U}(x)<0$ for some
(but not necessarily all) $x \in [0,\infty[$.
Notice that in \eqref{fo.twoin} we do not demand the validity of both inequalities, which might lead to
the effect that the constructed Hellinger integral upper bounds have to be cut off at $1$
for some (but not all) observation horizons $n\in \mathbb{N}$; see \eqref{fo.genbounds} below.
For the formulation of our first assertions on Hellinger integrals, we make use of the following notation:


\begin{defi}\label{defseqnull}
For all $\qui \in \quaset \times ]0,1[$ and all $p \in [0,\infty[$, $q \in[0,\infty[$, let us define 
the 
sequences 
$\left(a_{n}^{(q)}\right)_{n\in\mathbb{N}_{0}}$ and $\left(b_{n}^{(p,q)}\right)_{n\in\mathbb{N}_{0}}$ 
recursively by
\bea
a^{(q)}_{0}
:= 0 &;& \qquad a^{(q)}_{n}
\ := \ e^{a^{(q)}_{n-1}}\cdot q - \beta_{\lambda},~~n\in\mathbb{N},\label{defan}\\[0.05cm]
b^{(p,q)}_{0}
:= 0 &;& \qquad b^{(p,q)}_{n}
\ := \ e^{a^{(q)}_{n-1}}\cdot p-\alpha_{\lambda}
,~~n\in\mathbb{N}.
\label{defbn}
\eea
\end{defi}
Notice the interrelation $a^{(q_{\lambda}^{A})}_{1} = s_{\lambda}^{A}$ and $b^{(p_{\lambda}^{A},q_{\lambda}^{A})}_{1}= r_{\lambda}^{A}$ for $A \in\{E, L, U \}$.
Clearly, for $q \in ]0,\infty[$, $p \in [0,\infty[$, one has the linear interrelation
\be 
b^{(p,q)}_{n} = \frac{p}{q} \, a_{n}^{(q)} \, + \, \frac{p}{q} \, \beta_{\lambda} \, - \, \alpha_{\lambda},~~n\in\mathbb{N}.
\label{fo.anbn}
\ee



\noindent
Accordingly, we obtain fundamental 
Hellinger integral evaluations:

\begin{thm}\label{thm2}
(a) For all $\qui \in (\quasetNI \cup \quasetSPeins) \times ]0,1[$,  
all initial 
population sizes $\omega_{0}\in\mathbb{N}$ and  all observation horizons $n \in \mathbb{N}$ one can 
recursively
compute the 
\textbf{exact value} 
\be
H_{\lambda}(\Pna||\Pnh) \ = \ 
\exp\Big\{
a^{(q_{\lambda}^{E})}_{n}\, \omega_{0} \, + \, 
\frac{\aal}{\bal} \, \sum_{k=1}^{n} a^{(q_{\lambda}^{E})}_{k} 
\Big\} \ =: \ V_{\lambda,n} ,
\label{fo.genequality}
\ee
where $\frac{\aal}{\bal}$ can be equivalently replaced by $\frac{\ahy}{\bhy}$.
Recall that $q_{\lambda}^{E}:= \bal^{\lambda}\, \bhy^{1-\lambda}$. \\[0.05cm]
\noindent 
(b) For \textit{all} $\qui \in \quasetSPcomp \times ]0,1[$, 
\textit{all} coefficients $p_{\lambda}^{U} \in [0,\infty[$, $q_{\lambda}^{U}\in [0,\infty[$,
$p_{\lambda}^{L} \in [0,\min\{p_{\lambda}^{U},\alpha_{\lambda}\}]$, 
$q_{\lambda}^{L} \in [0,\min\{q_{\lambda}^{U},\beta_{\lambda}\}]$, 
such that \eqref{fo.varphibou} holds for all $x \in \mathbb{N}_{0}$ as well as \eqref{fo.twoin},
\textit{all} initial 
population sizes $\omega_{0}\in\mathbb{N}$ and 
\textit{all} observation horizons $n \in \mathbb{N}$ one gets the \textbf{recursive} (i.e.\ recursively computable) \textbf{bounds} $~ B_{\lambda,n}^{L}<H_{\lambda}(\Pna||\Pnh)<B_{\lambda,n}^{U}$~, where
\be
B_{\lambda,n}^{L} \ := \ 
\exp\Big\{
a^{(q_{\lambda}^{L})}_{n}\, \omega_{0} \, + \, 
\sum_{k=1}^{n} b^{(p_{\lambda}^{L},q_{\lambda}^{L})}_{k}
\Big\}
\quad
\textrm{and} \quad B_{\lambda,n}^{U} \ := \ \min \left\{
\exp\Big\{
a^{(q_{\lambda}^{U})}_{n}\, \omega_{0} \, + \, 
\sum_{k=1}^{n} b^{(p_{\lambda}^{U},q_{\lambda}^{U})}_{k}
\Big\}\, , \, 1 \right\}  \ .
\label{fo.genbounds}
\ee
\end{thm} 
%

\begin{rem}
From the proof below one can see that both parts of Theorem \ref{thm2} remain true
for the cases $\lambda \notin [0,1]$. For the (to our context)
incompatible setup of GWI with Poisson offspring but nonstochastic immigration of constant value 1,
the exact values of the corresponding Hellinger integrals (i.e.\ an ``analogue'' of part (a))
was established in Linkov and Lunyova \cite{Link96}. 

\end{rem}

%
%
\noindent
\pr We first prove the upper bound $B_{\lambda,n}^{U}$.
Let us fix $\qui$, $p_{\lambda}^{U}$, $q_{\lambda}^{U}$, $\omega_{0}\in\mathbb{N}$ as described in part (b). 
From \eqref{fo.heins}, \eqref{defphi}, \eqref{defvarphi}, \eqref{defflambda} and \eqref{fo.varphibou}
one gets immediately $B_{\lambda,1}^{U}$, and with the help of \eqref{fo.hellzw}
 for all observation horizons $n \in \mathbb{N}\backslash\{1\}$ (with the obvious shortcut
 for $n=2$)
\bea
&& H_{\lambda}\left(\Pna||\Pnh\right)
\ = \ \sum_{\omega_{1}=0}^{\infty}\cdots\sum_{\omega_{n-1}=0}^{\infty} \, \prod_{k=1}^{n-1}Z^{(\lambda)}_{n,k}
(\omega)
\cdot\exp\Big\{
\varphi_{\lambda}(\omega_{n-1}) -  f_{\lambda}(\omega_{n-1})
\Big\} \notag\\
&&
< \ \sum_{\omega_{1}=0}^{\infty}\cdots\sum_{\omega_{n-1}=0}^{\infty} \, \prod_{k=1}^{n-1}Z^{(\lambda)}_{n,k}
(\omega)
\cdot\exp\Big\{
(p_{\lambda}^{U} - \alpha_{\lambda}) + 
(q_{\lambda}^{U} - \beta_{\lambda})  \, \omega_{n-1}
\Big\} \notag \\
&&
= \ \sum_{\omega_{1}=0}^{\infty}\cdots\sum_{\omega_{n-1}=0}^{\infty} \, \prod_{k=1}^{n-1}Z^{(\lambda)}_{n,k}
(\omega)
\cdot\exp\Big\{
b^{(p_{\lambda}^{U},q_{\lambda}^{U})}_{1} + 
a^{(q_{\lambda}^{U})}_{1}  \, \omega_{n-1}
\Big\} \notag 
\\ 
&&
= \ \exp\Big\{ b^{(p_{\lambda}^{U},q_{\lambda}^{U})}_{1} \Big\} \ 
\sum_{\omega_{1}=0}^{\infty}\cdots\sum_{\omega_{n-2}=0}^{\infty} \, \prod_{k=1}^{n-2}Z^{(\lambda)}_{n,k}
(\omega)
\cdot\exp\Big\{
\exp\left\{ a^{(q_{\lambda}^{U})}_{1} \right\}
\ \varphi_{\lambda}(\omega_{n-2}) -  f_{\lambda}(\omega_{n-2})
\Big\} \notag\\
&&
< \ \exp\Big\{ b^{(p_{\lambda}^{U},q_{\lambda}^{U})}_{1} \Big\} \ 
\sum_{\omega_{1}=0}^{\infty}\cdots\sum_{\omega_{n-2}=0}^{\infty} \, \prod_{k=1}^{n-2}Z^{(\lambda)}_{n,k}
(\omega) \notag \\
&&
\cdot\exp\Big\{
\left( \exp\left\{ a^{(q_{\lambda}^{U})}_{1} \right\} \, p_{\lambda}^{U} - \alpha_{\lambda} \right)
+ \left( \exp\left\{ a^{(q_{\lambda}^{U})}_{1} \right\} \, q_{\lambda}^{U} - \beta_{\lambda} \right)
\cdot \omega_{n-2}
\Big\} \notag \\
&&
< \ \exp\Big\{ b^{(p_{\lambda}^{U},q_{\lambda}^{U})}_{1} \Big\} \ 
\sum_{\omega_{1}=0}^{\infty}\cdots\sum_{\omega_{n-2}=0}^{\infty} \, \prod_{k=1}^{n-2}Z^{(\lambda)}_{n,k}
(\omega) 
\cdot\exp\Big\{
b^{(p_{\lambda}^{U},q_{\lambda}^{U})}_{2} + a^{(q_{\lambda}^{U})}_{2} \, \omega_{n-2}
\Big\} \notag \\
&&
< \ \cdots \ < \
\exp\Big\{
a^{(q_{\lambda}^{U})}_{n}\, \omega_{0} \, + \, 
\sum_{k=1}^{n} b^{(p_{\lambda}^{U},q_{\lambda}^{U})}_{k}
\Big\} \ .
\label{fo.proofup}
\eea
Notice that for the strictness of the above inequalities 
we have used the fact that $\phi_{\lambda}(x) < \phi_{\lambda}^{U}(x)$  for some (in fact, all but at most two) 
$x \in \mathbb{N}_{0}$ (cf.\ (p-xiv) below). 
Since for some admissible choices of $p_{\lambda}^{U},q_{\lambda}^{U}$
and some $n \in \mathbb{N}$
the last term in \eqref{fo.proofup} can become larger than 1, 
one needs to take into
account the cutoff-point $1$ arising from \eqref{fo.hellgenup}. Notice that without
assumption \eqref{fo.twoin}, the last term in \eqref{fo.proofup} would always be larger than 1
(and thus useless). 
The lower bound $B_{\lambda,n}^{L}$ of part (b), as well as the exact value of part (a) follow
from \eqref{fo.hellzw} in an analoguous manner by employing $p_{\lambda}^{L},q_{\lambda}^{L}$ and 
$p_{\lambda}^{E},q_{\lambda}^{E}$ respectively. Furthermore, we use the fact that 
for $\qui \in (\quasetNI \cup \quasetSPeins) \times ]0,1[$ one gets from \eqref{fo.anbn}
the relation
$b^{(p_{\lambda}^{E},q_{\lambda}^{E})}_{n} = \frac{\aal}{\bal} \,  a^{(q_{\lambda}^{E})}_{n}$. 
For the sake of brevity, the corresponding straightforward details are omitted here.
Although we take the minimum of the upper bound derived in \eqref{fo.proofup} and 1, the inequality 
$B_{\lambda,n}^{L}<B_{\lambda,n}^{U}$ is nevertheless valid: the reason is that for constituting a lower bound, 
the parameters $p_{\lambda}^{L},q_{\lambda}^{L}$ must fulfil either the conditions
[$p_{\lambda}^{L}<0
$ and $
q_{\lambda}^{L}\leq0$] or  [$p_{\lambda}^{L}\leq0
$ and $
q_{\lambda}^{L}<0$] (or both).
\qed\\





\section{Detailed analyses}
\label{secDET}
For part (b) in Theorem \ref{thm2}, we have assumed the existence of reasonable linear lower and upper bounds of 
$\varphi_{\lambda}$ and $\phi_{\lambda}$. In the following, we shall carry out a more detailed analysis
addressing questions upon the non-uniqueness (and thus, flexibility) of the coefficients
$p_{\lambda}^{L}$, $q_{\lambda}^{L}$, $p_{\lambda}^{U}$, $q_{\lambda}^{U}$ in \eqref{fo.varphibou}, their ``optimal respectively reasonable choices'',
as well as the corresponding behaviour of the Hellinger integrals $H_{\lambda}(\Pna||\Pnh)$ as the observation 
horizon $n$ increases and finally converges to $\infty$. Of course, the answers to these
questions will depend on the (e.g.\ fixed) value of $\qua$ and the (e.g.\ selectable) value of $\lambda$.\\[0.1cm]

\noindent
Before starting a closer inspection, notice by induction the general fact
that for $\qui \in \quaset \times ]0,1[$ and $q \in ]0,\infty[$ the principal behaviour of the sequence $\left(a^{(q)}_{n}\right)_{n\in\mathbb{N}}$
is strongly governed by its 
first element:

\begin{enumerate} 

\item[(p-i)] $a^{(q)}_{n} \equiv 0$, 
\hspace{7.75cm} if \ $a^{(q)}_{1} = q - \beta_{\lambda} = 0$ \ (i.e.\ $q = \beta_{\lambda}$),

\item[(p-ii)] $\left(a^{(q)}_{n}\right)_{n\in\mathbb{N}}$ is strictly negative and strictly decreasing,
\hspace{0.9cm} if \ $a^{(q)}_{1} < 0$, 

\item[(p-iii)]  $\left(a^{(q)}_{n}\right)_{n\in\mathbb{N}}$ is strictly positive and strictly increasing,
\hspace{1.1cm} if \ $a^{(q)}_{1} > 0$. \\[0.1cm]

\end{enumerate}

\noindent 
Due to the linear interrelation \eqref{fo.anbn}, the monotonicity carries over to
the sequence $\left(b_{n}^{(p,q)}\right)_{n\in\mathbb{N}_{0}}$ 
($p \in [0,\infty[$, $q \in ]0,\infty[$) 
in the following way:
\begin{enumerate}

\item[(p-iv)] $b_{n}^{(0,q)} \equiv \,  - \alpha_{\lambda} \ < 0$,

\item[(p-v)] $b_{n}^{(p,q)} \equiv  \, p \, - \, \alpha_{\lambda}$, 
\hspace{3.65cm} if \ $q = \beta_{\lambda}$,

\item[(p-vi)] $\left(b_{n}^{(p,q)}\right)_{n\in\mathbb{N}}$ is strictly decreasing,
\hspace{0.9cm} if \ $q < \beta_{\lambda}$, 

\item[(p-vii)]  $\left(b_{n}^{(p,q)}\right)_{n\in\mathbb{N}}$ is strictly increasing,
\hspace{1.0cm} if \ $q > \beta_{\lambda}$. \\[0.1cm]

\end{enumerate}

\noindent
Notice that the sign of $b_{n}^{(p,q)}$ might not be same as the sign of $a^{(q)}_{n}$ (see e.g.\ (p-i), (p-iv)). 
Finally, for the remaining case one trivially gets 
\begin{enumerate}

\item[(p-viii)] $a^{(0)}_{n} \equiv - \beta_{\lambda}$, \qquad  $b_{n}^{(p,0)} \equiv e^{- \beta_{\lambda}}\cdot p-\alpha_{\lambda}$ \ \  ($p \geq 0$).

\end{enumerate}

\noindent
Moreover, for $\qui \in \quaset \times ]0,1[$ and $q \in ]0,\infty[$ we shall sometimes use the function 
\be\label{defxi}
\xi^{(q)}_{\lambda}(x)~:=~q\cdot e^{x}-\beta_{\lambda}, \quad x \in \mathbb{R},
\ee
which has the following obvious properties: 
\begin{enumerate}

\item[(p-ix)] $\xi^{(q)}_{\lambda}$ is strictly increasing, strictly conxex and smooth,

\item[(p-x)] $\lim_{x \rightarrow - \infty} \,  \xi^{(q)}_{\lambda}(x) \ = \ -\beta_{\lambda} \ < 0$, \
$\lim_{x \rightarrow \infty} \, \xi^{(q)}_{\lambda}(x) \ = \ \infty$.

\end{enumerate}

\noindent
With these auxilliary basic facts in hand, let us now start our detailed investigations of the time-behaviour $n \mapsto H_{\lambda}(\Pna||\Pnh)$ for the exactly
treatable case (a) in Theorem \ref{thm2}.




\subsection{Detailed analysis of the exact values}
\label{secDETEX}
~\\ \textbf{(aNI) \, The non-immigration case $\qui  \in \quasetNI \times ]0,1[$:}\\[0.05cm] 

\noindent
Recall that for this set-up we derived
$q_{\lambda}^{E}:= \bal^{\lambda}\, \bhy^{1-\lambda}> 0$ 
and 
$p_{\lambda}^{E}:= \aal^{\lambda}\, \ahy^{1-\lambda} = 0$.
According to Lemma \ref{lem2}, one has $q_{\lambda}^{E} < \beta_{\lambda}$
and thus, \,  $\left(a^{(q_{\lambda}^{E})}_{n}\right)_{n\in\mathbb{N}}$ is strictly negative as well as strictly decreasing.
Furthermore, because of (p-ix), (p-x) and  $a^{(q_{\lambda}^{E})}_{1} < 0$, the function 
$\xi^{(q_{\lambda}^{E})}_{\lambda}$ hits on $]-\infty,0]$ the straight line $id(x):=x$ once and only once.
Consequently, $\left(a^{(q_{\lambda}^{E})}_{n}\right)_{n\in\mathbb{N}}$ converges to 
the unique solution $x_0^{(q_{\lambda}^{E})} \in ]-\beta_{\lambda},a^{(q_{\lambda}^{E})}_{1}[$ of the equation 
\be
\xi^{(q_{\lambda}^{E})}_{\lambda}(x)~=~ q_{\lambda}^{E}\cdot e^x - \beta_{\lambda} = x , 
\quad x < 0.
\label{fo.zero}
\ee  
Summing up, we have shown the following detailed behaviour of Hellinger integrals:

\begin{prop}\label{propNI}
For all $\qui \in \quasetNI  \times ]0,1[$  
and all initial 
population sizes $\omega_{0}\in\mathbb{N}$ there holds
\bea
& (a) & \quad H_{\lambda}(\Peinsa||\Peinsh) \ = \  
\exp\Big\{
\Big(\bal^{\lambda}\, \bhy^{1-\lambda} - \lambda \bal - (1-\lambda) \bhy   \Big)\, x_{0} 
\Big\}
\ < \ 1, \hspace{3cm} ~ \notag\\
& (b) & \quad \textrm{the sequence } \left(H_{\lambda}(\Pna||\Pnh)\right)_{n\in\mathbb{N}} \ \textrm{ given by}\notag\\
&&  \hspace{2cm} H_{\lambda}(\Pna||\Pnh) \ = \ 
\exp\Big\{
a^{(q_{\lambda}^{E})}_{n}\, \omega_{0} 
\Big\} 
\ =: \ 
V_{\lambda,n} \notag\\
&& \quad \textrm{is strictly decreasing,}\notag\\
& (c) & \quad \lim_{n \rightarrow \infty} \, H_{\lambda}(\Pna||\Pnh) \ = \  
\exp\Big\{
x_{0}^{(q_{\lambda}^{E})} \, \omega_{0} 
\Big\} \in \, ]0,1[ \ , \notag\\
& (d) & \quad \lim_{n \rightarrow \infty} \, \frac{1}{n} \,  \log H_{\lambda}(\Pna||\Pnh) \ = \  
0 \ .  \notag
\eea
\end{prop}

\newpage

\noindent
\textbf{(aEF) \, The ``equal-fraction-case'' $\qui \in \quasetSPeins \times ]0,1[$:}\\[-0.05cm] 

\noindent
Again, one has
$q_{\lambda}^{E}:= \bal^{\lambda}\, \bhy^{1-\lambda}> 0$. 
Furthermore, 
$p_{\lambda}^{E}:= \aal^{\lambda}\, \ahy^{1-\lambda} > 0$, which leads to the abovementioned relation
$b^{(p_{\lambda}^{E},q_{\lambda}^{E})}_{n} = \frac{\aal}{\bal} \,  a^{(q_{\lambda}^{E})}_{n}$.
Hence, the results about the sequence $\left(a^{(q_{\lambda}^{E})}_{n}\right)_{n\in\mathbb{N}}$
coincide with those of the non-immigration case. This implies also that
the sequence  $\left(\sum_{k=1}^{n} \, a^{(q_{\lambda}^{E})}_{k}\right)_{n\in\mathbb{N}}$
is strictly negative, strictly decreasing and  converges to $-\infty$.
Hence, we get 

\begin{prop}\label{propPSP1}
For all $\qui \in \quasetSPeins \times ]0,1[$  
and all initial 
population sizes $\omega_{0}\in\mathbb{N}$ there holds 
\bea
& (a) & \quad H_{\lambda}(\Peinsa||\Peinsh) \ = \  
\exp\left\{
\Big(\bal^{\lambda}\, \bhy^{1-\lambda} - \lambda \bal - (1-\lambda) \bhy   \Big)\, \left(\omega_{0} + \frac{\aal}{\bal} \right) 
\right\}
\ < \ 1, \hspace{1cm} ~\notag\\
& (b) & \quad \textrm{the sequence} \ \left(H_{\lambda}(\Pna||\Pnh)\right)_{n\in\mathbb{N}} \ \textrm{given by}\notag\\
&& \hspace{2cm} H_{\lambda}(\Pna||\Pnh) \ = \ 
\exp\left\{
a^{(q_{\lambda}^{E})}_{n}\, \omega_{0} \, + \, 
\frac{\aal}{\bal} \, \sum_{k=1}^{n} a^{(q_{\lambda}^{E})}_{k} 
\right\}  \ =: \ 
V_{\lambda,n} \notag\\
&& \quad \textrm{is strictly decreasing,} \notag\\
& (c) & \quad \lim_{n \rightarrow \infty} \, H_{\lambda}(\Pna||\Pnh) \ = \  
 0 \ ,  \notag\\
 & (d) & \quad \lim_{n \rightarrow \infty} \, \frac{1}{n} \,  \log H_{\lambda}(\Pna||\Pnh) \ = \  
\frac{\aal}{\bal} \ x_{0}^{(q_{\lambda}^{E})} \ . \notag 
\eea
 
\end{prop}

\begin{rem}
For the (to our context)
incompatible setup of GWI with Poisson offspring but nonstochastic immigration of constant value 1,
an ``analogue'' of part (d) of Proposition \ref{propPSP1}
was established in Linkov and Lunyova \cite{Link96}.
\end{rem}




\subsection{Detailed analysis of the lower bounds}
\label{secDETLOW}
In this section we  
assume $\qui \in \quasetSPcomp \times ]0,1[$
and thus $\aal >0$,  $\ahy >0$, $\frac{\aal}{\ahy} \ne \frac{\bal}{\bhy}$, $\gamma \ne 0$,
$f_{\Al}(x)>0$, $f_{\Hy}(x)>0$ ($x \in [0,\infty[$).
Concerning \eqref{fo.phibou}, let us derive a lower linear bound $\phi_{\lambda}^{L}(\cdot)$ of $\phi_{\lambda}(\cdot)$
which is optimal. In order to achieve this, one can use the following straightforward properties
of $\phi_{\lambda}(x)$, $x\in[0,\infty[$ \, (cf.\ \eqref{defphi}):

\begin{enumerate}

\item[(p-xi)]  $\phi_{\lambda}(0) = \aal^{\lambda} \ahy^{1-\lambda} - \alpha_{\lambda} \leq 0$ (cf.\ Lemma \ref{lem2}),
with equality iff $\aal=\ahy$ (together with $\bal \ne \bhy$).

\item[(p-xii)]  $\phi'_{\lambda}(x)~=~\lambda\beta_{\Al}\left(f_{\Al}(x)\right)^{\lambda-1}\left(f_{\Hy}(x)\right)^{1-\lambda}+(1-\lambda)\beta_{\Hy}\left(f_{\Al}(x)\right)^{\lambda}\left(f_{\Hy}(x)\right)^{-\lambda} -  \beta_{\lambda} \  > \ - \beta_{\lambda}$.

\item[(p-xiii)] $\lim_{x \rightarrow \infty} \phi'_{\lambda}(x) = \bal^{\lambda} \bhy^{1-\lambda} - \beta_{\lambda} \leq 0$ (cf.\ Lemma \ref{lem2}),
with equality iff $\bal=\bhy$ (together with $\aal \ne \ahy$).

\item[(p-xiv)] $\phi''_{\lambda}(x)~=~-\lambda(1-\lambda)\left(f_{\Al}(x)\right)^{\lambda-2}\left(f_{\Hy}(x)\right)^{-\lambda-1}
\gamma^2 < 0$, i.e.\ the function $\phi_{\lambda}(\cdot)$ is strictly concave;  \ notice that 
$\phi'_{\lambda}(0)
=~\lambda\beta_{\Al}\left(\aal/\ahy\right)^{\lambda-1}+(1-\lambda)\beta_{\Hy}\left(\aal/\ahy\right)^{\lambda} -  \beta_{\lambda}$
 can be \textit{either}
negative (e.g.\ for $\qui=(4,2,3,1,0.5)$), \textit{or} zero 
(e.g.\ for $\qui=(4,2,4,1,0.5)$), \textit{or} positive (e.g.\ for $\qui=(4,2,5,1,0.5)$).  
Accordingly, the strict decreasingness and continuity of 
$\phi'_{\lambda}(\cdot)$ as well as (p-xiii) imply that $\phi_{\lambda}(\cdot)$ can be 
\textit{either} strictly decreasing, \textit{or} 
can obtain its global maximum on $]0,\infty[$, \textit{or} -- only in the case $\bal=\bhy$ -- can be strictly increasing.

\item[(p-xv)] $\lim_{x\rightarrow\infty}\Big(\phi_{\lambda}(x)-\left(\widetilde{r_{\lambda}} + \widetilde{s_{\lambda}}\, x\right)\Big)~=~0$ \ 
for 
$\widetilde{r_{\lambda}} :=\lambda\alpha_{\Al}\left[\left(\frac{\beta_{\Al}}{\beta_{\Hy}}\right)^{\lambda-1}-1\right]+(1-\lambda)
\alpha_{\Hy}\left[\left(\frac{\beta_{\Al}}{\beta_{\Hy}}\right)^{\lambda} -1 \right]$ 
and $\widetilde{s_{\lambda}} := \beta_{\Al}^{\lambda}\beta_{\Hy}^{1-\lambda} - \beta_{\lambda} \leq 0$;
notice that $\widetilde{s_{\lambda}}=0$ iff $\bal=\bhy$ (together with $\aal \ne \ahy$). Furthermore, $\phi_{\lambda}(0) < \widetilde{r_{\lambda}}$
(cf.\ Lemma \ref{lem2}). If $\aal=\ahy$ (and thus $\bal\ne\bhy$) then the intercept $\widetilde{r_{\lambda}}$ is strictly positive,
whereas for the case $\aal\ne\ahy$ the intercept $\widetilde{r_{\lambda}}$ can take any sign
(take e.g.\ $\qui =(3.7,0.9,2.0,1.0,0.5)$ for $\widetilde{r_{\lambda}}>0$,
$\qui =(3.6,0.9,2.0,1.0,0.5)$ for $\widetilde{r_{\lambda}}=0$,
$\qui =(3.5,0.9,2.0,1.0,0.5)$ for $\widetilde{r_{\lambda}}<0$). 
\end{enumerate}

\noindent
From (p-xi) to (p-xv) it is easy to see that for all current parameter constellations  
the particular choices 
$$p_{\lambda}^{L} := \aal^{\lambda} \ahy^{1-\lambda} > 0, \qquad
q_{\lambda}^{L} := \bal^{\lambda} \bhy^{1-\lambda} > 0 \ 
$$ 
-- 
which correspond to the choices
$$r_{\lambda}^{L} := \aal^{\lambda} \ahy^{1-\lambda} - \alpha_{\lambda} \leq 0, \qquad
s_{\lambda}^{L} := \bal^{\lambda} \bhy^{1-\lambda} - \beta_{\lambda} \leq 0
$$  
in \eqref{fo.phibou} (and at least one of the two last inequalities is strict) -- 
lead to the tightest lower bound $B_{\lambda,n}^{L}$ for $H_{\lambda}(\Pna||\Pnh)$ in
\eqref{fo.genbounds}. 
This situation coincides partially with those in Section \ref{secDETEX}. Formally,
$p_{\lambda}^{L}= p_{\lambda}^{E}$ and $q_{\lambda}^{L}= q_{\lambda}^{E}$,
but because of $\gamma \ne 0$
the relation
$b^{(p_{\lambda}^{L},q_{\lambda}^{L})}_{n} = \frac{\aal}{\bal} \,  a^{(q_{\lambda}^{L})}_{n}$
is in general not valid anymore and has to be replaced by the relation (cf.\ \eqref{fo.anbn})
\be
b^{(p_{\lambda}^{L},q_{\lambda}^{L})}_{n} \ = \ 
\left(\frac{\aal}{\bal}\right)^{\lambda} 
\left(\frac{\ahy}{\bhy}\right)^{1-\lambda} \,  
a_{n}^{(q_{\lambda}^{L})} \, + \, \left(\frac{\aal}{\bal}\right)^{\lambda} 
\left(\frac{\ahy}{\bhy}\right)^{1-\lambda} \, \beta_{\lambda} \, - \, \alpha_{\lambda},~~n\in\mathbb{N}.
\label{fo.anbnLOW}
\ee
\noindent
Hence, for a better distinguishability and easier reference
we stick to the $L-$notation here. Nevertheless, the behaviour of the sequence
$\left(a^{(q_{\lambda}^{L})}_{n}\right)_{n\in\mathbb{N}}$ coincides
exactly with that of the sequence $\left(a^{(q_{\lambda}^{E})}_{n}\right)_{n\in\mathbb{N}}$
in the Subsections \ref{secDETEX}(aNI), (aEF). In particular $\left(a^{(q_{\lambda}^{L})}_{n}\right)_{n\in\mathbb{N}}$ is strictly negative,
strictly decreasing and converges to the unique solution $x_0^{(q_{\lambda}^{L})} \in ]-\infty,a^{(q_{\lambda}^{L})}_{1}[$ of the equation 
\be
\xi^{(q_{\lambda}^{L})}_{\lambda}(x)~=~ q_{\lambda}^{L}\cdot e^x - \beta_{\lambda} = x , 
\quad x < 0 \ .
\label{fo.zeroLOW}
\ee  
Consequently, because of \eqref{fo.anbnLOW} and
$b^{(p_{\lambda}^{L},q_{\lambda}^{L})}_{1} = 
\aal^{\lambda} \ahy^{1-\lambda} - \alpha_{\lambda} \leq 0$ (cf.\ \eqref{defbn}),
the sequence\\  
$\left(b^{(p_{\lambda}^{L},q_{\lambda}^{L})}_{n}\right)_{n\in\mathbb{N}\backslash\{1\}}$
is strictly negative and strictly decreasing. As in Subsection \ref{secDETEX}(aEF), we
obtain

\begin{prop}\label{propLOW}
For all $\qui \in \quasetSPcomp \times ]0,1[$ 
and all initial 
population sizes $\omega_{0}\in\mathbb{N}$ there holds
\bea
& (a) &  B_{\lambda,1}^{L} \ := \  
\exp\Big\{
\Big(\bal^{\lambda}\, \bhy^{1-\lambda} - \lambda \bal - (1-\lambda) \bhy   \Big)\, \omega_{0} + 
\Big(\aal^{\lambda}\, \ahy^{1-\lambda} - \lambda \aal - (1-\lambda) \ahy   \Big) 
\Big\}
 <  1,  \notag\\
& (b) &  \textrm{the sequence} \ \left(B_{\lambda,n}^{L}\right)_{n\in\mathbb{N}} \ \textrm{of lower bounds for} 
\ \left(H_{\lambda}(\Pna||\Pnh)\right)_{n\in\mathbb{N}} \ \textrm{given}\notag\\
&& \ \textrm{by} \  B_{\lambda,n}^{L} \ := \  
\exp\Big\{
a^{(p_{\lambda}^{L})}_{n}\, \omega_{0} \, + \, 
\sum_{k=1}^{n} b^{(p_{\lambda}^{L},q_{\lambda}^{L})}_{k} 
\Big\} \ \textrm{\, is strictly decreasing,}\notag\\
& (c) &  \lim_{n \rightarrow \infty} \, B_{\lambda,n}^{L} \ = \  
 0 \ ,  \notag\\
 & (d) &  \lim_{n \rightarrow \infty} \, \frac{1}{n} \, \log B_{\lambda,n}^{L} \ = \ 
\frac{p_{\lambda}^{L}}{q_{\lambda}^{L}} \, \left(x_0^{(q_{\lambda}^{L})} + \beta_{\lambda}\right)
- \alpha_{\lambda} 
\ .  
\notag
\eea
\end{prop}




\subsection{Detailed analysis of the upper bounds}
\label{secDETUP}
As above, we again assume $\qui \in \quasetSPcomp \times ]0,1[$ throughout this section.
In contrast to the treatment of the lower bounds in Section \ref{secDETLOW},
the finetuning of the upper bounds is more involved. 
Because of the strict concavity of the function $\phi_{\lambda}(\cdot)$ (cf.\ (p-xiv)),
there is in general no overall best linear upper bound of $\phi_{\lambda}(\cdot)$
within the framework \eqref{fo.phibou}. Different reasonable goals might lead to 
different reasonable choices of $p_{\lambda}^{U}$, $q_{\lambda}^{U}$ (and thus
of $r_{\lambda}^{U}$, $s_{\lambda}^{U}$) which might imply different
behaviour of the corresponding sequence $\left(B_{\lambda,n}^{U}\right)_{n\in\mathbb{N}}$
of upper bounds in \eqref{fo.genbounds}. This can be conjectured from the following
immediate monotonicity properties: 

\begin{enumerate} 

\item[(p-xvi)] $0 \leq q_{1} < q_{2} \ \Longrightarrow \ a^{(q_{1})}_{n} < \, a^{(q_{2})}_{n}$ for all $n \in \mathbb{N}$. 

\item[(p-xvii)] 
Trivially, $b_{n}^{(0,q_{1})} = \, b_{n}^{(0,q_{2})} \equiv - \alpha_{\lambda}$. 
In contrast, let $p \in ]0,\infty[$ be fixed; 
then,  $0 \leq q_{1} < q_{2} \ \Longrightarrow \ b_{n}^{(p,q_{1})} < \, b_{n}^{(p,q_{2})}$ for all $n \in \mathbb{N}$.

\item[(p-xviii)] Let $q \in [0,\infty[$ be fixed. 
Then, $0 \leq p_{1} < p_{2} \ \Longrightarrow \ b_{n}^{(p_{1},q)} < \, b_{n}^{(p_{2},q)}$ for all $n \in \mathbb{N}$.

\item[(p-xix)]  
$0 \leq p_{1} < p_{2}, \, 0 \leq q_{1} < q_{2} \ \Longrightarrow \ b_{n}^{(p_{1},q_{1})} < \, b_{n}^{(p_{2},q_{2})}$ 
for all $n \in \mathbb{N}$.

\item[(p-xx)] For the case $0 \leq p_{1} < p_{2}, \, 0 \leq q_{2} < q_{1}$  
there is in general no dominance assertion for $b_{n}^{(p_{1},q_{1})}$, \, $b_{n}^{(p_{2},q_{2})}$ which holds for all $n \in \mathbb{N}$;
take e.g.\ $\qui=(1,0.6,3,3,0.5)$, $p_{1}=3.4641$, $q_{1}=0.7785$ (for which $\phi_{\lambda}^{U}(\cdot)$ corresponds
to the secant line through the points $\phi_{\lambda}(0)$ and $\phi_{\lambda}(1)$), as well as  
$p_{2}=3.4857$, $q_{2}=0.7746$ 
(for which $\phi_{\lambda}^{U}(\cdot)$ corresponds
to the asymptote of $\phi_{\lambda}$), and inspect the first six values of of the corresponding $b_n-$sequence. 
\end{enumerate}

\noindent
The properties (p-xvi) to (p-xx) have corresponding effects on the behaviour\\ 
$(p_{\lambda}^{U},q_{\lambda}^{U}) \mapsto B_{\lambda,n}^{U} \ = \  \min \left\{
\exp\Big\{
a^{(q_{\lambda}^{U})}_{n}\, \omega_{0} \, + \, 
\sum_{k=1}^{n} b^{(p_{\lambda}^{U},q_{\lambda}^{U})}_{k}
\Big\}\, , \, 1 \right\}$ (cf.\ \eqref{fo.genbounds})
of the upper bounds. For instance, for any fixed admissible intercept $p_{\lambda}^{U}$
one would always choose the smallest admissible $q_{\lambda}^{U}$ in order to achieve the smallest
possible upper bound;  
due to (p-xiv) this implies that on the ultimately relevant subdomain $\mathbb{N}_{0}$ 
the linear function $\phi_{\lambda}^{U}(\cdot)$ should hit 
$\phi_{\lambda}(\cdot)$ in at least one but at most two points (tangent or secant line).
Furthermore, we require for the rest of the section that $p_{\lambda}^{U} >0$ and $q_{\lambda}^{U}>0$, 
because otherwise $r_{\lambda}^{U} < \phi_{\lambda}(0)$ and $s_{\lambda}^{U}<\widetilde{s_{\lambda}}$ (cf.\ (p-xv))
which contradicts to the nature of linear upper bounds of $\phi_{\lambda}$.
\\[0.1cm] 

\noindent
The (only partially restricted) choice of parameters $p_{\lambda}^{U}$, $q_{\lambda}^{U}$ for the upper bounds $B_{\lambda,n}^{U}$
can be made according to different, partially incompatible (``optimality-'' respectively ``goodness-'') criteria,
such as:\\
(Ga) \, very good tightness for $n \geq N$ for some fixed large $N \in \mathbb{N}$, or\\ 
(Gb) \, for a fixed initial population size $\omega_{0} \in \mathbb{N}$ there holds $B_{\lambda,n}^{U} < 1$ for all $n \in \mathbb{N}$, or \\
(Gc) \, there holds $B_{\lambda,n}^{U} < 1$ for all $n \in \mathbb{N}$ 
and all $\omega_{0} \in \mathbb{N}$ (strict improvement of the general upper \\
\indent \hspace{0.4cm} bound \eqref{fo.hellgenup}).\\
For the sake of brevity, we investigate only goal (Gc) (with the exception of Subsection \ref{secDETUP}(a7) 
 and Theorem \ref{thm.entexUP})
which can be achieved if (and ``nearly but not fully'' iff) 
\eqref{fo.twoin} holds; 
this can be seen from
\be
B_{\lambda,1}^{U} \ = \  
\min \left\{
\exp\Big\{
a^{(q_{\lambda}^{U})}_{1}\, \omega_{0} \, + \, 
 b^{(p_{\lambda}^{U},q_{\lambda}^{U})}_{1}
 \Big\}\, , \, 1 \right\}  \ = \ 
\min \left\{
\exp\Big\{
(q - \beta_{\lambda})\, \omega_{0} + 
(p-\alpha_{\lambda}) 
\Big\}\, , \, 1 \right\} 
\notag
\ee
and the properties (p-i) to (p-vii). 
Furthermore, (p-xiv) and (p-xv) imply that the slope  $s_{\lambda} := q_{\lambda}^{U} - \beta_{\lambda}$ 
in \eqref{fo.phibou} should be greater
or equal to the limit slope $\widetilde{s_{\lambda}}$ which leads to the restriction
$q_{\lambda}^{U}\geq \bal^{\lambda} \bhy^{1-\lambda}$.
Moreover, since $s_{\lambda}\leq 0$, the intercept 
$r_{\lambda} := p_{\lambda}^{U} - \alpha_{\lambda}$ in \eqref{fo.phibou} should be greater or equal to $\phi_{\lambda}(0)$ and
thus, $p_{\lambda}^{U}\geq \aal^{\lambda} \ahy^{1-\lambda}$.
By comparing the above established lower and upper parameter-bounds, from Lemma \ref{lem2}
it follows that
the case $q_{\lambda}^{U} < \beta_{\lambda}$ automatically implies $\bal \ne \bhy$ whereas  
the case $p_{\lambda}^{U} < \alpha_{\lambda}$ leads to $\aal \ne \ahy$.  In consistence with (p-xiv),
various different parameter constellations can lead to different Hellinger-integral-upper-bound details, which we investigate 
in the following.\\[0.0cm]

\noindent
\textbf{
(a1) \, The case 
$\quasetSPzwei$ of all (componentwise) strictly positive $\qua$ with $\bhy \ne \bal$, $\aal = \ahy$
} \\[0.05cm] 

\noindent
We have $\phi_{\lambda}(0)=0$ (cf.\ (p-xi)), $\phi'_{\lambda}(0)=0$ (cf.\ (p-xii)). Thus, the only admissible intercept choice
is $r_{\lambda}^{U} = 0 = p_{\lambda}^{U} - \alpha_{\lambda} = b^{(p_{\lambda}^{U},q_{\lambda}^{U})}_{1}$ 
(i.e.\ $p_{\lambda}^{U} = \alpha_{\lambda} = \alpha_{\bullet}>0$), and the minimal admissible slope 
which implies \eqref{fo.phibou} for $x \in \mathbb{N}$ is given by
$s_{\lambda}^{U} = \frac{\phi_{\lambda}(1) - \phi_{\lambda}(0)}{1-0} = q_{\lambda}^{U} - \beta_{\lambda}
= a^{(q_{\lambda}^{U})}_{1} <0$ (i.e.\ $q_{\lambda}^{U} = 
(\alpha_{\bullet}+\bal)^{\lambda}  (\alpha_{\bullet}+\bhy)^{1-\lambda} - \alpha_{\bullet}>0$).
Analogously to Subsection \ref{secDETEX}(aNI), one can derive
that $\left(a^{(q_{\lambda}^{U})}_{n}\right)_{n\in\mathbb{N}}$ is strictly negative, strictly decreasing,
and converges to 
the unique solution $x_0^{(q_{\lambda}^{U})} \in ]-\infty,a^{(q_{\lambda}^{U})}_{1}[$ of the equation 
\be\label{fo.zeroUPa1}
\xi_{\lambda}^{(q_{\lambda}^{U})}(x)~=~q_{\lambda}^{U}\cdot e^x - \beta_{\lambda}~=~x,~~x<0 \ . 
\ee
Moreover, in the same manner as 
in Section \ref{secDETLOW}, the sequence  
$\left(b^{(p_{\lambda}^{U},q_{\lambda}^{U})}_{n}\right)_{n\in\mathbb{N}\backslash\{1\}}$
is strictly negative and strictly decreasing. This leads to  

\begin{prop}\label{propUPa1}
For all $\qui \in \quasetSPzwei \times ]0,1[$ 
and all initial  
population sizes $\omega_{0}\in\mathbb{N}$ there holds
\bea
& (a) & B_{\lambda,1}^{U} \ := \  
\exp\left\{
\Big( q_{\lambda}^{U} - \beta_{\lambda} \Big)\, \omega_{0} + 
\Big( p_{\lambda}^{U} - \alpha_{\lambda}  \Big) 
\right\}
\ < \ 1,  \notag\\
& (b) & \textrm{the sequence} \ \left(B_{\lambda,n}^{U}\right)_{n\in\mathbb{N}} \ \textrm{of upper bounds for}
\ \left(H_{\lambda}(\Pna||\Pnh)\right)_{n\in\mathbb{N}} \ \textrm{given} \hspace{5.0cm} \ \notag\\
& & \textrm{by} \ B_{\lambda,n}^{U} \ := \ 
\exp\Big\{
a^{(q_{\lambda}^{U})}_{n}\, \omega_{0} \, + \, 
\sum_{k=1}^{n} b^{(p_{\lambda}^{U},q_{\lambda}^{U})}_{k} 
\Big\} \ \textrm{\ is strictly decreasing,}\notag\\
& (c) & \lim_{n \rightarrow \infty} \, B_{\lambda,n}^{U} \ = \  
 0 \ = \ \lim_{n \rightarrow \infty} \,H_{\lambda}(\Pna||\Pnh),\notag\\
 & (d) & \lim_{n \rightarrow \infty} \, \frac{1}{n} \, \log B_{\lambda,n}^{U} \ = \ 
\frac{p_{\lambda}^{U}}{q_{\lambda}^{U}} \, \left(x_0^{(q_{\lambda}^{U})} + \beta_{\lambda}\right)
- \alpha_{\lambda} 
\ .  \notag
\eea
\end{prop}

\noindent In contrast to $\quasetSPzwei$, the constellation $\quasetSPdrei$ of all (componentwise) strictly positive $\qua$ with $\aal \ne \ahy$,
$\bal \ne \bhy$ and $\frac{\aal}{\ahy} \ne \frac{\bal}{\bhy}$ is divided into three main parts as follows:
because of Lemma \ref{lem2} one gets 
on the domain $]0,\infty[$  
the relation
$\phi_{\lambda}(x) =0$  
iff $f_{\Al}(x)= f_{\Hy}(x)$ 
iff 
$x = x^{*} := \frac{\ahy -\aal}{\bal - \bhy} > 0$. Accordingly, (for reasons which will be explained below) 
we denote by $\quasetSPdreiab$ resp.\ $\quasetSPdreic$
resp.\ $\quasetSPdreid$ the subset of $\quasetSPdrei$ for which $\frac{\ahy -\aal}{\bal - \bhy} < 0$ resp.\
$\frac{\ahy -\aal}{\bal - \bhy} \in ]0,\infty[\backslash \mathbb{N}$ resp.\ 
$\frac{\ahy -\aal}{\bal - \bhy} \in \mathbb{N}$; notice that the case $\frac{\ahy -\aal}{\bal - \bhy}=0$
can not appear within $\quasetSPdrei$. For further investigations let us first divide
the set $\quasetSPdreiab \times ]0,1[$ of quintuples $\qui$ into two parts $\quisetSPdreia$ and $\quisetSPdreib$:\\[0.1cm]

\noindent
\textbf{
(a2) \, The case 
$\quisetSPdreia$ of all $\qui \in \quasetSPdreiab \times ]0,1[$ for which
$
\lambda\beta_{\Al}\left(\aal/\ahy\right)^{\lambda-1}+(1-\lambda)\beta_{\Hy}\left(\aal/\ahy\right)^{\lambda} -  \beta_{\lambda}
\leq 0$ \ holds 
} \\[0.05cm]

\noindent  
From (p-xi) and (p-xii), one gets $\phi_{\lambda}(0) < 0$ and $\phi'_{\lambda}(0) \leq 0$. For the latter, both the strict
negativity as well as the vanishing can appear in the current parameter setup,
take e.g.\ $\qui =(1.8,0.9,2.8,0.7,0.5)$ for $\phi'_{\lambda}(0) = 0$ and 
$\qui =(1.8,0.9,2.7,0.7,0.5)$  for $\phi'_{\lambda}(0) < 0$.
In the current setup,  $\phi_{\lambda}$ is 
a strictly negative, strictly decreasing, and -- due to (p-xiv) -- strictly concave function (and thus,
the assumption $\frac{\ahy -\aal}{\bal - \bhy} < 0$ is superfluous here). In contrast to Subsection (a1),
one has the flexibility to choose the intercept $p_{\lambda}^{U}$ from the nonempty interval $[\aal^{\lambda} \ahy^{1-\lambda},\alpha_{\lambda}]$
and the slope $q_{\lambda}^{U}$ from the nonempty interval $[\bal^{\lambda} \bhy^{1-\lambda},\beta_{\lambda}]$,
subject to the constraints that $(p_{\lambda}^{U},q_{\lambda}^{U})\ne(\alpha_{\lambda},\beta_{\lambda})$ and
$\phi_{\lambda}(x) 
\ \leq \ (p_{\lambda}^{U} - \alpha_{\lambda}) + (q_{\lambda}^{U} - \beta_{\lambda}) \, x  
$ (cf.\ \eqref{fo.phibou}). Of course, one way to obtain a reasonable choice of intercept and slope is the search for the
optimum
\be
\left(\overline{p_{\lambda}^{U}},\overline{q_{\lambda}^{U}}\right) :=
argmin_{(p,q)} \left\{ \exp\Big\{
a^{(q)}_{n}\, \omega_{0} \, + \, 
\sum_{k=1}^{n} b^{(p,q)}_{k}
\Big\}  \right\}
\label{fo.optG1}
\ee
subject to the abovementioned constraints. However, the corresponding result generally depends on the choice of the initial
population size $\omega_{0}$ and the observation horizon $n$. Hence, there is in general no overall
optimal choice of $p_{\lambda}^{U}$, $q_{\lambda}^{U}$ (without the incorporation of further goal-dependent constraints
such as $\lim_{n \rightarrow \infty} \, B_{\lambda,n}^{U} \ = \  0 $ in case of $\lim_{n \rightarrow \infty} \,H_{\lambda}(\Pna||\Pnh)
=0$). By the way, due to the recursive nature of the sequences in $\eqref{fo.optG1}$ and the nontriviality of the constraints,
this optimization problem seems to be not straightforward to solve, in general.
\\[0.1cm]

\noindent
Inspired from Subsection (a1), a more pragmatic but yet reasonable choice is the following:  take any intercept $p_{\lambda}^{U} \in [\aal^{\lambda} \ahy^{1-\lambda},\alpha_{\lambda}]$
such that $(p_{\lambda}^{U}-\alpha_{\lambda}) + 2 (\phi_{\lambda}(1) - (p_{\lambda}^{U}-\alpha_{\lambda}))  \geq \phi_{\lambda}(2)$
\ (i.e.\ $2 \left(\aal + \bal\right)^{\lambda}\, \left(\ahy + \bhy\right)^{1-\lambda} - p_{\lambda}^{U} + \alpha_{\lambda}
\geq  \left(\aal + 2\bal\right)^{\lambda}\, \left(\ahy + 2\bhy\right)^{1-\lambda}$) \ and  \ $q_{\lambda}^{U}:= \phi_{\lambda}(1) - (p_{\lambda}^{U}-\alpha_{\lambda}) + \beta_{\lambda}
= \left(\aal + \bal\right)^{\lambda}\, \left(\ahy + \bhy\right)^{1-\lambda} - p_{\lambda}^{U}$, which corresponds to 
a linear function $\phi_{\lambda}^{U}$ which is\\
\noindent (a) \, nonpositive  on $\mathbb{N}_{0}$ and strictly negative on $\mathbb{N}$,\\
\noindent (b) \, larger than or equal to $\phi_{\lambda}$ on $\mathbb{N}_{0}$, strictly larger than $\phi_{\lambda}$ on $\mathbb{N}\backslash\{1,2\}$, 
and equal to $\phi_{\lambda}$ at\\ 
\indent \hspace{0.28cm} the point
$x=1$ (``discrete tangent or secant line through $x=1$'').\\[0.15cm] 
\noindent 
One can easily see that (due to the restriction \eqref{fo.varphibou}) not all $p_{\lambda}^{U}\in[\aal^{\lambda}\ahy^{1-\lambda},\alpha_{\lambda}]$ 
might qualify for the current purpose.
For the particular choice
$p_{\lambda}^{U} = \aal^{\lambda} \ahy^{1-\lambda}$ and $q_{\lambda}^{U} = \left(\aal + \bal\right)^{\lambda}\, \left(\ahy + 
\bhy\right)^{1-\lambda} - \aal^{\lambda} \ahy^{1-\lambda}$
one obtains 
$r_{\lambda}^{U} = p_{\lambda}^{U} - \alpha_{\lambda} = b^{(p_{\lambda}^{U},q_{\lambda}^{U})}_{1} < 0$ 
(cf.\ Lemma \ref{lem2}) and 
$s_{\lambda}^{U} = q_{\lambda}^{U} - \beta_{\lambda} = \phi_{\lambda}(1)- \phi_{\lambda}(0)
= a^{(q_{\lambda}^{U})}_{1} <0$ (secant line through $\phi_{\lambda}(0)$ and $\phi_{\lambda}(1)$).
Hence, analogously to Subsection (a1) one can derive
that $\left(a^{(q_{\lambda}^{U})}_{n}\right)_{n\in\mathbb{N}}$ is strictly negative, strictly decreasing,
and converges to 
the unique solution $x_0^{(q_{\lambda}^{U})} \in ]-\infty,a^{(q_{\lambda}^{U})}_{1}[$ of equation \eqref{fo.zeroUPa1}.
Moreover, the sequence  
$\left(b^{(p_{\lambda}^{U},q_{\lambda}^{U})}_{n}\right)_{n\in\mathbb{N}\backslash\{1\}}$
is strictly negative and strictly decreasing. Thus,  
all the assertions (a), (b), (c), (d) of Proposition \ref{propUPa1} hold
for all $\qui \in \quisetSPdreia$ 
and all initial 
population sizes $\omega_{0}\in\mathbb{N}$.\\

\noindent
\textbf{
(a3) \, The case 
$\quisetSPdreib$ of all $\qui \in \quasetSPdreiab \times ]0,1[$ for which
$
\lambda\beta_{\Al}\left(\aal/\ahy\right)^{\lambda-1}+(1-\lambda)\beta_{\Hy}\left(\aal/\ahy\right)^{\lambda} -  \beta_{\lambda}
> 0$ \ holds
} \\[0.05cm] 

\noindent
In this situation (which appears e.g.\ for $\qui =(1.8,0.9,2.9,0.7,0.5)$)
one gets from (p-xi) and (p-xii) the two inequalities $\phi_{\lambda}(0) < 0$ and $\phi'_{\lambda}(0) > 0$. 
Furthermore, 
in accordance with the arguments in the forefront of Subsection (a2), 
$\phi_{\lambda}$ is 
a strictly negative, strictly concave, hump-shaped (cf.\ (p-xiii)) function. One can proceed similarly to (a2).
Indeed, let $x_{\max} := \textrm{argmax}_{x\in ]0,\infty[} \phi_{\lambda}(x)$ which is the unique solution
of 
\be
\lambda \bal      \left[ \left(\frac{f_{\Al}(x)}{f_{\Hy}(x)}\right)^{\lambda-1} - 1 \right] \, + \, 
(1-\lambda) \bhy  \left[ \left(\frac{f_{\Al}(x)}{f_{\Hy}(x)}\right)^{\lambda} - 1 \right]
\ = \ 0 \ , \qquad x \in ]0,\infty[ \ ,
\label{fo.xmax}
\ee
(cf.\ (p-xii), (p-xiv)); notice that $x^{*}$ formally satisfies the equation \eqref{fo.xmax} but 
does not qualify because of the current restriction $x^{*}<0$.\\[0.05cm]

\noindent
 Let us first inspect the case 
$\phi_{\lambda}(\lfloor x_{\max} \rfloor) > \phi_{\lambda}(\lfloor x_{\max} \rfloor+1)$, 
where $\lfloor x \rfloor$ denotes the integer part of $x$. Consider the subcase 
$\phi_{\lambda}(\lfloor x_{\max} \rfloor)+\lfloor x_{\max} \rfloor \left(\phi_{\lambda}(\lfloor x_{\max} \rfloor)
-\phi_{\lambda}(\lfloor x_{\max} \rfloor+1)\right)\leq0$, which means that the secant line through 
$\phi_{\lambda}(\lfloor x_{\max} \rfloor)$ and $\phi_{\lambda}(\lfloor x_{\max} \rfloor+1)$ possesses a non-positive intercept. 
In this situation it is reasonable to choose as \textit{intercept} any $p_{\lambda}^{U} - \alpha_{\lambda} = b^{(p_{\lambda}^{U},q_{\lambda}^{U})}_{1} = r_{\lambda}^{U} \in [\phi_{\lambda}(\lfloor x_{\max} \rfloor),\phi_{\lambda}(\lfloor x_{\max} \rfloor)+\lfloor x_{\max} \rfloor
\left(\phi_{\lambda}(\lfloor x_{\max} \rfloor)-\phi_{\lambda}(\lfloor x_{\max} \rfloor+1)\right)]$,
and as corresponding \textit{slope} $q_{\lambda}^{U} - \alpha_{\lambda} = a^{(q_{\lambda}^{U})}_{1} = s_{\lambda}^{U} = \\
\frac{\phi_{\lambda}(\lfloor x_{\max} \rfloor) - r_{\lambda}^{U}}{(\lfloor x_{\max} \rfloor) - 0} \, \leq 0$. 
A larger intercept
would lead to a linear function $\phi_{\lambda}^{U}$ for which \eqref{fo.phibou} is not valid at $\lfloor x_{\max} \rfloor+1$.\\[0.05cm]

\noindent
In the other subcase $\phi_{\lambda}(\lfloor x_{\max} \rfloor)+x_{\max}\left(\phi_{\lambda}(\lfloor x_{\max} \rfloor)-\phi_{\lambda}(\lfloor x_{\max} \rfloor+1)\right)>0$, one can choose any intercept $p_{\lambda}^{U} - \alpha_{\lambda} = b^{(p_{\lambda}^{U},q_{\lambda}^{U})}_{1} = r_{\lambda}^{U} \in [\phi_{\lambda}(\lfloor x_{\max} \rfloor),0]$
and as corresponding slope  $q_{\lambda}^{U} - \alpha_{\lambda} = a^{(q_{\lambda}^{U})}_{1} = s_{\lambda}^{U} = 
\frac{\phi_{\lambda}(\lfloor x_{\max} \rfloor) - r_{\lambda}^{U}}{(\lfloor x_{\max} \rfloor) - 0} \, \leq 0$
(notice that the corresponding line $\phi_{\lambda}^{U}$ is on $]\lfloor x_{\max} \rfloor ,\infty[$
strictly larger than the secant line through 
$\phi_{\lambda}(\lfloor x_{\max} \rfloor)$ and $\phi_{\lambda}(\lfloor x_{\max} \rfloor+1)$).\\

\noindent
If $\phi_{\lambda}(\lfloor x_{\max} \rfloor) \leq \phi_{\lambda}(\lfloor x_{\max} \rfloor+1)$, one can
proceed as above by substituting 
the crucial pair of points $(\lfloor x_{\max} \rfloor, \lfloor x_{\max} \rfloor +1)$
with $(\lfloor x_{\max} \rfloor +1, \lfloor x_{\max} \rfloor +2)$ and examining the
analogous two subcases. \\[0.05cm]

\noindent
With the accordingly derived $p_{\lambda}^{U}$, $q_{\lambda}^{U}$ one gets in all four (sub)cases 
exactly the same
kind of behaviour of the sequences $\left(a^{(q_{\lambda}^{U})}_{n}\right)_{n\in\mathbb{N}}$,   
$\left(b^{(p_{\lambda}^{U},q_{\lambda}^{U})}_{n}\right)_{n\in\mathbb{N}\backslash\{1\}}$
as in Subsection (a2). Hence, 
all the assertions (a), (b), (c), (d) of Proposition \ref{propUPa1} hold 
for all $\qui \in \quisetSPdreib$ 
and all initial  
population sizes $\omega_{0}\in\mathbb{N}$.\\ 

\noindent
\textbf{
(a4) \, The case 
$\quasetSPdreic$ of all (componentwise) strictly positive $\qua$ with $\aal \ne \ahy$,
$\bal \ne \bhy$, $\frac{\aal}{\ahy} \ne \frac{\bal}{\bhy}$ and $\frac{\ahy -\aal}{\bal - \bhy} \in ]0,\infty[ \backslash\mathbb{N}$
} \\[0.05cm]

\noindent
The only difference to Subsection (a3) is that the maximum value of $\phi_{\lambda}(\cdot)$ now achieves $0$, 
at the positive \textit{non-integer} point 
$x_{\max} = x^{*} = \frac{\ahy -\aal}{\bal - \bhy} \in ]0,\infty[ \backslash\mathbb{N}$ (take e.g.\ $\qui =(1.8,0.9,1.1,3.0,0.5)$ as an
example). Due to (p-xi), (p-xii) and (p-xiv) one gets automatically 
$
\lambda\beta_{\Al}\left(\aal/\ahy\right)^{\lambda-1}+(1-\lambda)\beta_{\Hy}\left(\aal/\ahy\right)^{\lambda} -  \beta_{\lambda}
> 0$ for all $\lambda \in ]0,1[$.
This situation can be treated exactly as in (a3). Consequently,
all the assertions (a), (b), (c), (d) of Proposition \ref{propUPa1} hold  
for all $\qui \in \quasetSPdreic \times]0,1[$ 
and all initial  
population sizes $\omega_{0}\in\mathbb{N}$. \\[0.1cm]

\noindent
\textbf{
(a5) \, The case 
$\quasetSPdreid$ of all (componentwise) strictly positive $\qua$ with $\aal \ne \ahy$,
$\bal \ne \bhy$, $\frac{\aal}{\ahy} \ne \frac{\bal}{\bhy}$ and 
$\frac{\ahy -\aal}{\bal - \bhy} \in \mathbb{N}$  
} \\[0.05cm]

 \noindent
The only difference to Subsection (a4) is that the maximum value of $\phi_{\lambda}(\cdot)$ now achieves $0$ at the \textit{integer} point 
$x_{\max} = x^{*} = \frac{\ahy -\aal}{\bal - \bhy} \in \mathbb{N}$ 
(take e.g.\ $\qui =(1.8,0.9,1.2,3.0,0.5)$ as an example).
Under the restriction that $\exp\Big\{
a^{(q_{\lambda}^{U})}_{n}\, \omega_{0} \, + \, 
\sum_{k=1}^{n} b^{(p_{\lambda}^{U},q_{\lambda}^{U})}_{k}
\Big\} \leq 1$ for all $n \in \mathbb{N}$ and all $\omega_{0} \in \mathbb{N}$, our method
leads to the choices $r_{\lambda}^{U}=0$ as well as $s_{\lambda}^{U}=0$.
Consequently, $B_{\lambda,n}^{U} \equiv 1$, which coincides with the general upper bound \eqref{fo.hellgenup}, but violates 
the abovementioned desired goal (Gc).\\ 
However, by using a conceptually different method we can nevertheless prove the convergence
\be \label{fo.hellnull3d}
\lim_{n \rightarrow \infty} \, H_{\lambda}(\Pna||\Pnh) \ = \  0
\ee 
(which will be used for
the study of entire separation below). This will be done in Appendix \ref{App3}.\\[0.2cm] 

\noindent As a next step, let us investigate the last possible parameter constellation:\\ 

\noindent
\textbf{
(a6) \, The case 
 $\quasetSPvier$ of all (componentwise) strictly positive $\qua$ with  $\aal \ne \ahy$,
$\bal = \bhy$
} \\[-0.1cm] 

\noindent
This is the only case where $\phi_{\lambda}(\cdot)$ is strictly negative and strictly increasing,
with $\lim_{x \rightarrow \infty} \phi_{\lambda}(x) = \lim_{x \rightarrow \infty} \phi'_{\lambda}(x)
= 0$, leading to the choices $r_{\lambda}^{U}=0$ as well as $s_{\lambda}^{U}=0$
under the restriction that $\exp\Big\{
a^{(q_{\lambda}^{U})}_{n}\, \omega_{0} \, + \, 
\sum_{k=1}^{n} b^{(p_{\lambda}^{U},q_{\lambda}^{U})}_{k}
\Big\} \leq 1$ for all $n \in \mathbb{N}$ and all $\omega_{0} \in \mathbb{N}$.
Consequently, $B_{\lambda,n}^{U} \equiv 1$, which is consistent with the general upper bound \eqref{fo.hellgenup}, but violates the
abovementioned desired Goal (Gc). Unfortunately, the proof method of \eqref{fo.hellnull3d} can't be carried over 
to the current setup (see Appendix \ref{App3}).\\[0.25cm]
\noindent
\textbf{(a7) \, Alternative bounds for $\quasetSPzwei \cup \quasetSPdreiab 
\cup \quasetSPdreic \cup \quasetSPdreid$
} \\[-0.1cm] 

\noindent
Within this last subsection, let us exceptionally ignore the Goal (Gc). Correspondingly, 
for the derivation of an upper bound $\widetilde{B_{\lambda,n}^{U}}$ one can use the asymptote of $\varphi_{\lambda}$
given in (p-xv) to end up with $\widetilde{p}_{\lambda}^{\, U}:= \widetilde{r_{\lambda}} +  \alpha_{\lambda} = \lambda \, \aal\left(\frac{\bal}{\bhy}\right)^{\lambda-1}+(1-\lambda)\, \ahy\left(\frac{\bal}{\bhy}\right)^{\lambda}$ as well as 
$\widetilde{q}_{\lambda}^{\, U}=
\widetilde{s_{\lambda}} + \beta_{\lambda} = \bal^{\lambda}\bhy^{1-\lambda}$. Clearly, $\widetilde{p}_{\lambda}^{\, U} > 
p_{\lambda}^{L} = \aal^{\lambda}\ahy^{1-\lambda}$
by Lemma \ref{lem2} and $\widetilde{q}_{\lambda}^{\, U}= q_{\lambda}^{L}$. Furthermore, $\widetilde{q}_{\lambda}^{\, U} < \beta_{\lambda}$
and thus \eqref{fo.twoin} holds, since we have excluded $\quasetSPvier$. However -- depending on the
choice of $\qua$ -- the intercept $\widetilde{r_{\lambda}} =  \widetilde{p}_{\lambda}^{\, U} - \alpha_{\lambda}$
may become strictly positive, and hence
\begin{equation*}
\widetilde{B_{\lambda,1}^{U}} \  :=  \ 
\exp\Big\{
a^{(\widetilde{q}_{\lambda}^{\, U})}_{1}\, \omega_{0} \, + \, 
b^{(\widetilde{p}_{\lambda}^{\, U},\widetilde{q}_{\lambda}^{\, U})}_{1}
\Big\} 
\ = \ 
\exp\Big\{
\left( \widetilde{q}_{\lambda}^{\, U} - \beta_{\lambda}  \right) \cdot  \omega_{0} \, + \, 
\widetilde{p}_{\lambda}^{\, U} - \alpha_{\lambda}
\Big\} 
\end{equation*} 
may become larger than 1. However, according to properties (p-ii) and (p-vi) the sequence
\be
n \ \mapsto \ \widetilde{B_{\lambda,n}^{U}} \  :=  \ 
\exp\Big\{
a^{(\widetilde{q}_{\lambda}^{\, U})}_{n}\, \omega_{0} \, + \, 
\sum_{k=1}^{n} b^{(\widetilde{p}_{\lambda}^{\, U},\widetilde{q}_{\lambda}^{\, U})}_{k}
\Big\}
 =  
\exp\Big\{
a^{(\widetilde{q}_{\lambda}^{\, U})}_{n} \omega_{0} \, + \, \frac{\widetilde{p}_{\lambda}^{\, U}}{\widetilde{q}_{\lambda}^{\, U}} \, 
\sum_{k=1}^{n} a^{(\widetilde{q}_{\lambda}^{\, U})}_{k}
+ \left( \frac{\widetilde{p}_{\lambda}^{\, U}}{\widetilde{q}_{\lambda}^{\, U}} \, \beta_{\lambda} - \alpha_{\lambda} \right) \cdot n
\Big\}
\notag
\ee
may become smaller than 1. 
Let us therefore define for all $n \in \mathbb{N}$ and all $\lambda \in ]0,1[$
\begin{equation*}
\widetilde{\widetilde{B_{\lambda,n}^{U}}} \  :=  \ 
\min\left\{  \, \widetilde{B_{\lambda,n}^{U}} , 1   \right\} \ 
\end{equation*} 
which can be used as an upper bound for the case  $\quasetSPdreid \times ]0,1[$. 

\vspace{0.2cm} 
\noindent
For the other cases
$(\quasetSPzwei \times ]0,1[) \cup \quisetSPdreia \cup \quisetSPdreib \cup (\quasetSPdreic \times ]0,1[)$
all the assertions (a),(b),(c) of Proposition \ref{propUPa1} remain valid for
replacing $B_{\lambda,n}^{U}$  by the improved upper bound
\begin{equation*}
B_{\lambda,n}^{U,\text{impr}} \  :=  \ 
\min\left\{ B_{\lambda,n}^{U} \, , \widetilde{B_{\lambda,n}^{U}}   \right\} \ < \ 1 \ .
\end{equation*} 

\noindent
In fact, for all these parameter classes  
there are 
concrete examples such that the upper bound $B_{\lambda,n}^{U,\text{impr}}$ really improves the upper bound $B_{\lambda,n}^{U}$ for all $n\in\mathbb{N}$ (i.e. $\widetilde{B_{\lambda,n}^{U}}<B_{\lambda,n}^{U}$). For $\quasetSPzwei \times ]0,1[$ take e.g. $\qui=(0.8,0.6,2,2,0.5)$ and $\omega_{0}=10$, with $\widetilde{p}_{\lambda}^{\, U}=2.021$, $\widetilde{q}_{\lambda}^{\, U}=0.693$, instead of the proposed choice $p_{\lambda}^{U}=2$ and $q_{\lambda}^{U}=0.698$. For $\quisetSPdreia$ take e.g. $\qui=(0.8,0.6,2,1.9,0.5)$ and $\omega_{0}=10$, with $\widetilde{p}_{\lambda}^{\, U}=1.963$, $\widetilde{q}_{\lambda}^{\, U}=0.693$, instead of the proposed choice $p_{\lambda}^{U}=1.949$ and $q_{\lambda}^{U}=0.696$. For $\quisetSPdreib$ take e.g. $\qui=(0.8,0.6,2,1.1,0.5)$ and $\omega_{0}=10$, with $\widetilde{p}_{\lambda}^{\, U}=1.501$, $\widetilde{q}_{\lambda}^{\, U}=0.693$, instead of the (amongst others proposed) choice $p_{\lambda}^{U}=1.483$ and $q_{\lambda}^{U}=0.699$. For $\quasetSPdreic \times ]0,1[$ take e.g. $\qui=(1,1.5,2,1.8,0.5)$ and $\omega_{0}=10$, with $\widetilde{p}_{\lambda}^{\, U}=1.960$, $\widetilde{q}_{\lambda}^{\, U}=1.225$, instead of the (amongst others proposed) choice $p_{\lambda}^{U}=1.897$ and $q_{\lambda}^{U}=1.249$.




\subsection{
Asymptotic distinguishability}
\label{secCONTIG}

For each $n \in \mathbb{N}_{0}$, let $(\Omega_n, \mathcal{F}_n)$ be a measurable space equipped
with two probability measures $\widehat{P_n}$, $\overline{P_n}$. 
The following two general types of asymptotic distinguishability are well known (see e.g.\ LeCam \cite{LeC86}, 
Liese and Vajda \cite{LieVa87}, Jacod and Shiryaev \cite{JacSh87}, Linkov \cite{Link05}, and the references therein): 
\begin{itemize}
\item[(CEa)] the sequence $(\widehat{P_n})_{n\in \mathbb{N}_{0}}$ 
is contiguous to the sequence $(\overline{P_n})_{n\in \mathbb{N}_{0}}$
-- in symbols, $(\widehat{P_n}) \triangleleft (\overline{P_n})$) -- 
if for all sequences $A_n \in \mathcal{F}_n$ with $\lim_{n \rightarrow \infty} \overline{P_n}(A_n)=0$
there holds $\lim_{n \rightarrow \infty} \widehat{P_n}(A_n)=0$.
\item[(CEb)] the sequences $(\widehat{P_n})_{n\in \mathbb{N}_{0}}$ 
and $(\overline{P_n})_{n\in \mathbb{N}_{0}}$ are entirely separated (completely asymptotically separable) 
-- in symbols, $(\widehat{P_n}) \bigtriangleup (\overline{P_n})$ -- 
if there exist a sequence $n_m \uparrow \infty$ as $m \uparrow \infty$ and for each $m\in \mathbb{N}_{0}$ an
$A_{n_{m}} \in \mathcal{F}_{n_{m}}$ such that $\lim_{m \rightarrow \infty} \widehat{P_{n_{m}}}(A_{n_{m}})=1$ and 
$\lim_{m \rightarrow \infty} \overline{P_{n_{m}}}(A_{n_{m}})=0$.
\end{itemize}
The corresponding negations will be denoted by $\overline{\triangleleft}$ and $\overline{\bigtriangleup}$.
As demonstrated in the abovementioned references for a general context,\\
 (CEb) holds iff
$\liminf_{n \rightarrow \infty} H_{\lambda}\left(\widehat{P_n}||\overline{P_n}\right) =0$
for some (or equivalently, all) $\lambda \in ]0,1[$; \ furthermore, \\ 
(CEa)
holds iff $\liminf_{\lambda  \uparrow 1} \left\{ 
\liminf_{n \rightarrow \infty} H_{\lambda}\left(\widehat{P_n}||\overline{P_n}\right)
\right\} =1$.\\[0.1cm]
 Combining these results with the respective part (c) of Propositions \ref{propNI}, \ref{propPSP1} and \ref{propUPa1} 
as well as the connected investigations of Subsections \ref{secDETUP}(a2) to (a5), we obtain the following



\begin{cor}\label{corCONENT}
(a) For all 
$\qua \in (\quasetSPeins \cup \quasetSPzwei \cup \quasetSPdreiab 
\cup \quasetSPdreic
\cup \quasetSPdreid)$ 
and  all initial population sizes $\omega_{0}\in \mathbb{N}$, the corresponding
sequences $(\Pna)_{n \in \mathbb{N}_{0}}$ and $(\Pnh)_{n \in \mathbb{N}_{0}}$ are entirely separated.\\[0.1cm]
(b) For all $\qua \in \quasetNI$ with $\bal \leq 1$
and all initial population sizes $\omega_{0}\in\mathbb{N}$, 
the sequence $(\Pna)_{n \in \mathbb{N}_{0}}$ is contiguous to $(\Pnh)_{n \in \mathbb{N}_{0}}$.\\[0.1cm]
(c) For all $\qua \in \quasetNI$ with $\bal > 1$
and all initial population sizes $\omega_{0}\in\mathbb{N}$, 
the sequence $(\Pna)_{n \in \mathbb{N}_{0}}$ is neither contiguous to nor entirely separated to $(\Pnh)_{n \in \mathbb{N}_{0}}$.\\[-0.2cm]
\end{cor}

\begin{rems}
(i) Assertion (c) of Corollary \ref{corCONENT} contrasts the case of Gaussian processes with independent
increments where one gets either entire separation or mutual contiguity
(see e.g.\ Liese and Vajda \cite{LieVa87}).\\
(ii) By putting Corollary \ref{corCONENT}(b) and (c) together, we obtain for different ``criticality
pairs'' in the non-immigration case $\quasetNI$ the following asymptotic distinguishability types:\\
$(\Pna) \triangleleft  \triangleright (\Pnh)$ \, if $\bal \leq 1$, $\bhy \leq 1$; \qquad
$(\Pna) \triangleleft  \overline{\triangleright} \,  (\Pnh)$ \, if $\bal \leq 1$, $\bhy > 1$; \\
$(\Pna) \, \overline{\triangleleft}  \triangleright (\Pnh)$ \, if $\bal > 1$, $\bhy \leq 1$; \quad
$(\Pna) \, \overline{\triangleleft} \,   \overline{\triangleright} \,  (\Pnh)$ and $(\Pna) \overline{\bigtriangleup} (\Pnh)$ \, if $\bal > 1$, $\bhy > 1$; \\
in particular, for $\quasetNI$ the sequences $(\Pna)_{n \in \mathbb{N}_{0}}$ and $(\Pnh)_{n \in \mathbb{N}_{0}}$
are not completely asymptotically inseparable (indistinguishable).\\
(iii) In the light of the abovementioned (CEa) resp.\ (CEb) characteriztions by means of Hellinger integral
limits, the finite-time-horizon results on Hellinger integrals given in Theorem \ref{thm2},
Section \ref{secDET} and also in the following Section \ref{secCFB} can loosely be interpreted
as ``finite-sample (rather than asympotic) distinguishability'' assertions.
\end{rems}




\section{Closed-form bounds}\label{secCFB}

Depending on the parameter constellation, we have given bounds 
respectively exact values for the Hellinger integrals, which can be obtained 
with the help of recursions \eqref{defan} (together with \eqref{fo.anbn} respectively (p-viii))
which are ``stepwise fully evaluable'' but generally seem not to admit a closed-form
representation in the observation horizons $n$; consequently, the exact time-behaviour of (the
bounds of) the Hellinger integrals can generally not be seen explicitly. 
To avoid this intransparency (at the expense of losing some precision)
one can approximate \eqref{defan} by a recursion that allows for a closed-form representation.
Accordingly, we shall employ (context-adapted) linear inhomogeneous difference equations  
\bea
\widetilde{a}_{0}
:= \ 0 &;& \qquad \widetilde{a}_{n}
\ := \
\widetilde{\xi}\left(\widetilde{a}_{n-1}\right) \ + \ \rho_{n-1}
,~~n\in\mathbb{N}, \qquad \textrm{with}
\label{fo.defanCLO}
\eea 
 \vspace{-0.5cm}
\bea
\widetilde{\xi}(x) &:=&  c \ + \ d \cdot x \ , 
\hspace{2.8cm} x \in ]-\infty,0]\ ,
\label{fo.xiCLO}
\\
\rho_{n-1} &:= & K_{1}\cdot\varkappa^{n-1} \ + \ K_{2}\cdot\nu^{n-1} \, , \qquad n \in \mathbb{N},
\label{defrho}
\eea
for some constants $c \in ]-\infty,0[$, $d \in ]0,1[$, $K_{1}, K_{2},\varkappa,\nu \in \mathbb{R}$  
with $0 \leq\nu<\varkappa<d$. 
As usual,
one gets the closed-form representation
\be
\widetilde{a}_{n} \ =  \ \widetilde{a}_{n}^{hom} +  \widetilde{c}_{n}
\quad \textrm{with } \ \widetilde{a}_{n}^{hom} = c\cdot\frac{1-d^{n}}{1-d} 
\ \ \ \textrm{and  } \  \widetilde{c}_{n} = K_{1}\cdot\frac{d^{n}-\varkappa^{n}}{d-\varkappa} \ + \ K_{2}\cdot\frac{d^{n}-\nu^{n}}{d-\nu}
\label{fo.repanCLO} 
\ee
which immediately leads for all $n\in\mathbb{N}$ to 
 \vspace{0.1cm}
\be
\sum_{k=1}^{n} \widetilde{a}_{k} \ = \ \left(\frac{K_{1}}{d-\varkappa}+\frac{K_{2}}{d-\nu}-\frac{c}{1-d}\right)\cdot\frac{d\cdot\left(1-d^{n}\right)}{1-d}-\frac{K_{1}\cdot\varkappa\cdot\left(1-\varkappa^{n}\right)}{(d-\varkappa)(1-\varkappa)}-\frac{K_{2}\cdot\nu\cdot\left(1-\nu^{n}\right)}{(d-\nu)(1-\nu)}+\frac{c \cdot n}{1-d} \ .
\label{fo.repsumanCLO}
\ee
Notice that for the special case $K_{2}=-K_{1}>0$ one has from \eqref{defrho} for all integers $n\geq2$  the relation $\rho_{n-1}<0$ and thus 
$\widetilde{a}_{n} - \widetilde{a}_{n}^{hom} <0$, leading to
\be
\widetilde{c}_{n} <0  \quad \textrm{and}  \quad \sum_{k=1}^{n} \widetilde{c}_{n} < 0 \ .
\label{focntil}
\ee
In the following, we appropriately apply \eqref{fo.defanCLO}-\eqref{fo.repsumanCLO} to the different parameter contexts
of Section \ref{secDET}.




\subsection{Closed-form lower bounds}
\label{secCLOSLOW}

Let $\qui \in \quaset \times ]0,1[$. We have seen in the
Sections \ref{secDETEX} and \ref{secDETLOW} that
the determination of the exact values  and the lower bounds
had (more or less) identical structure:
choose $q_{\lambda}^{\bigstar}:= q_{\lambda}^{L}= q_{\lambda}^{E}= \bal^{\lambda} \bhy^{1-\lambda}  >0$,
compute the sequence $\left(a_{n}^{(q_{\lambda}^{\bigstar})}\right)_{n\in\mathbb{N}_{0}}$  
by the nonlinear recursion (cf.\ \eqref{defan}, \eqref{defxi})
\bea
a^{(q_{\lambda}^{\bigstar})}_{0}
:= 0 &;& \qquad a^{(q_{\lambda}^{\bigstar})}_{n}
\ := \
\xi^{(q_{\lambda}^{\bigstar})}_{\lambda}\hspace{-0.12cm}\left(a^{(q_{\lambda}^{\bigstar})}_{n-1}\right) 
,~~n\in\mathbb{N}, \label{defanvar}
\eea 
choose $p_{\lambda}^{\bigstar}:= p_{\lambda}^{L} = p_{\lambda}^{E} = \aal^{\lambda} \ahy^{1-\lambda} \geq 0$,
compute (cf.\ \eqref{fo.anbn})
\be
b^{(p_{\lambda}^{\bigstar},q_{\lambda}^{\bigstar})}_{n} \ = \ 
\left(\frac{\aal}{\bal}\right)^{\lambda} 
\left(\frac{\ahy}{\bhy}\right)^{1-\lambda} \,  
a_{n}^{(q_{\lambda}^{\bigstar})} \, + \, \left(\frac{\aal}{\bal}\right)^{\lambda} 
\left(\frac{\ahy}{\bhy}\right)^{1-\lambda} \, \beta_{\lambda} \, - \, \alpha_{\lambda},~~n\in\mathbb{N}, 
\notag
\ee
and finally end up with (cf.\ \eqref{fo.genbounds}, \eqref{fo.genequality}) $\exp\Big\{
a^{(q_{\lambda}^{\bigstar})}_{n}\, \omega_{0} \, + \, 
\sum_{k=1}^{n} b^{(p_{\lambda}^{\bigstar},q_{\lambda}^{\bigstar})}_{k}
\Big\} 
$ which is either interpreted as bound $B_{\lambda,n}^{L}$ in the parameter case $\qui \in \quasetSPcomp$ $ \times ]0,1[$
or as exact value $V_{\lambda,n}$ in the parameter case $\qui \in (\quasetNI \cup \quasetSPeins) \times ]0,1[$
(where we achieved some further simplifications above). Since 
$\left(a^{(q_{\lambda}^{\bigstar})}_{n}\right)_{n\in\mathbb{N}}$ is strictly negative, strictly decreasing and converges to the unique solution 
$x_0^{(q_{\lambda}^{\bigstar})} \in ]-\infty,a^{(q_{\lambda}^{\bigstar})}_{1}[ \, = \, ]-\infty,q_{\lambda}^{\bigstar}-\beta_{\lambda}[$ of the equation (cf.\ \eqref{fo.zeroLOW}, \eqref{fo.zero})
 \vspace{-0.2cm}
\be
\xi^{(q_{\lambda}^{\bigstar})}_{\lambda}(x)~=~ q_{\lambda}^{\bigstar}\cdot e^x - \beta_{\lambda} = x , 
\quad x < 0 \ ,
\label{fo.zeroSTAR}
\ee  
we 
use the following approximative linear recursion in order to obtain
a closed-form lower bound for both (here identically treatable) cases $L$, $E$:
\bea
\ua_{0}^{(q_{\lambda}^{\bigstar})}
:= 0 &;& \qquad \ua_{n}^{(q_{\lambda}^{\bigstar})}
\ := \
\xi_{\lambda}^{(q_{\lambda}^{\bigstar}),T}\hspace{-0.12cm}\left(\ua^{(q_{\lambda}^{\bigstar})}_{n-1}\right) \ + \ \urho_{n-1}^{(q_{\lambda}^{\bigstar})}
,~~n\in\mathbb{N}, \label{defua} 
\eea 
i.e.\ we replace  
the nonlinear function $\xi^{(q_{\lambda}^{\bigstar})}_{\lambda}(x)~=~q_{\lambda}^{\bigstar}\cdot e^{x}-\beta_{\lambda}$
by the tangent line of $\xi^{(q_{\lambda}^{\bigstar})}_{\lambda}$ at $x=x_0^{(q_{\lambda}^{\bigstar})}$ defined by
\be\label{fo.xiT} 
\xi_{\lambda}^{(q_{\lambda}^{\bigstar}),T}(x)\ := 
\  x_{0}^{(q_{\lambda}^{\bigstar})}\left(1-q_{\lambda}^{\bigstar}\cdot e^{x_{0}^{(q_{\lambda}^{\bigstar})}}\right) \, + \, q_{\lambda}^{\bigstar}\cdot e^{x_{0}^{(q_{\lambda}^{\bigstar})}} \cdot x 
\, , \qquad x\in[x_{0}^{(q_{\lambda}^{\bigstar})},0] \ ,
\ee
and reduce the error we face by adding the 
``correction-term'' 
\be\label{defurho}
\urho_{n-1}^{(q_{\lambda}^{\bigstar})} \ := \ \frac{1}{2}\cdot \left(x_{0}^{(q_{\lambda}^{\bigstar})}\right)^{2}\cdot\left(q_{\lambda}^{\bigstar}\cdot e^{x_{0}^{(q_{\lambda}^{\bigstar})}}\right)^{2n-1} \ > \ 0.
\ee
In other words, by means of the two functions on the domain $[0,\infty[$
\be
q \ \mapsto \  d^{(q),T} := q \cdot e^{x_{0}^{(q)}} \qquad 
q \ \mapsto \  \Gamma^{(q)} :=
\frac{q}{2} \cdot e^{x_{0}^{(q)}}\cdot\left(x_{0}^{(q)}\right)^{2} = 
\frac{d^{(q),T}}{2} \cdot \left(x_{0}^{(q)}\right)^{2}
\label{fo.dTundGamma}
\ee
we use \eqref{fo.defanCLO}, \eqref{fo.xiCLO}, \eqref{defrho}
with constants  
$d:= d^{(q_{\lambda}^{\bigstar}),T}  \in ]0,1[$,  $c:= 
x_{0}^{(q_{\lambda}^{\bigstar})} \cdot \left(1-d^{(q_{\lambda}^{\bigstar}),T}\right)  \in ]-\infty,0[$, 
\ $K_{1}:= \Gamma^{(q_{\lambda}^{\bigstar})} >0$,  
\ $\varkappa:= \left(d^{(q_{\lambda}^{\bigstar}),T}\right)^{2} \in\mathbb{R}\backslash\{d,1\}$,
\ $K_{2}:= 0 $,
\ $\nu:=0$. Let us first present some fundamental properties which will be proved in Appendix \ref{App4}:

\begin{lem}\label{lemua}
For all $\qui \in \quaset\times]0,1[$ there holds:
\bea
& (a) & \ua_{n}^{(q_{\lambda}^{\bigstar})} \ < \ a_{n}^{(q_{\lambda}^{\bigstar})}, \qquad \textrm{for all} \ n\in\mathbb{N}.\notag\\
& (b) & \textrm{The sequence} \ \left(\ua_{n}^{(q_{\lambda}^{\bigstar})}\right)_{n\in\mathbb{N}} \ \textrm{is strictly decreasing.}\hspace{6cm}\notag\\
& (c) & \lim_{n\rightarrow\infty}\ua_{n}^{(q_{\lambda}^{\bigstar})} \ = \ \lim_{n\rightarrow\infty}a_{n}^{(q_{\lambda}^{\bigstar})} \ = \ x_{0}^{(q_{\lambda}^{\bigstar})}.\notag
\eea
\end{lem}

\noindent
Applying Theorem \ref{thm2}, Lemma \ref{lemua} as well as the formulae \eqref{fo.anbn}, \eqref{fo.repanCLO} and \eqref{fo.repsumanCLO}, one gets



\begin{thm}\label{thm3varLOW}
For all $\qui \in \quaset \times ]0,1[$ and  all initial population sizes $\omega_{0}\in\mathbb{N}$ 
the following assertions hold:\\
\noindent
(a) for all observation horizons $n \in \mathbb{N}$ the Hellinger integral can be bounded from below by the closed-form bounds $H_{\lambda}(\Pna||\Pnh) \ > \ C_{\lambda,n}^{L}$ given by
\bea\label{fo.ansumbnLCLO}
&& \hspace{-0.8cm}C_{\lambda,n}^{L} \ := \ \exp\Bigg\{x_{0}^{(q_{\lambda}^{\bigstar})}\left[\omega_{0}-\frac{p_{\lambda}^{\bigstar}}{q_{\lambda}^{\bigstar}}\frac{d^{(q_{\lambda}^{\bigstar}),T}}{1-d^{(q_{\lambda}^{\bigstar}),T}}\right]\left(1-\left(d^{(q_{\lambda}^{\bigstar}),T}\right)^{n}\right) \ + \ \left[\frac{p_{\lambda}^{\bigstar}}{q_{\lambda}^{\bigstar}}\left(\beta_{\lambda}+x_{0}^{(q_{\lambda}^{\bigstar})}\right)-\alpha_{\lambda}\right]\cdot ~n\notag\\
&& \hspace{1.7cm} + \  \underline{\zeta}^{(q_{\lambda}^{\bigstar})}_{n} \cdot\omega_{0} \ + \ \underline{\vartheta}^{(q_{\lambda}^{\bigstar})}_{n}\Bigg\}, \quad \textrm{where for all} \ n\in\mathbb{N}
\eea
\vspace{-0.2cm}
\be\label{def.uzeta}
\underline{\zeta}^{(q_{\lambda}^{\bigstar})}_{n} \ := \ \Gamma^{(q_{\lambda}^{\bigstar})}\cdot\frac{\left(d^{(q_{\lambda}^{\bigstar}),T}\right)^{n-1}\left(1-\left(d^{(q_{\lambda}^{\bigstar}),T}\right)^{n}\right)}{1-d^{(q_{\lambda}^{\bigstar}),T}} \ > \ 0 \qquad \textrm{and}
\ee
\vspace{0.1cm}
\be\label{def.uvartheta}
\underline{\vartheta}^{(q_{\lambda}^{\bigstar})}_{n} \ := \ \frac{p_{\lambda}^{\bigstar}}{q_{\lambda}^{\bigstar}}\cdot\Gamma^{(q_{\lambda}^{\bigstar})}\cdot\frac{\left(1-\left(d^{(q_{\lambda}^{\bigstar}),T}\right)^{n}\right)}{\left(1-d^{(q_{\lambda}^{\bigstar}),T}\right)^{2}}\cdot\left(1-d^{(q_{\lambda}^{\bigstar}),T}\frac{\left(1+\left(d^{(q_{\lambda}^{\bigstar}),T}\right)^{n}\right)}{1+d^{(q_{\lambda}^{\bigstar}),T}}\right) \ > \ 0.
\ee

\noindent
(b)
the sequence $\left(C_{\lambda,n}^{L}\right)_{n\in\mathbb{N}}$ is strictly decreasing.\\
\noindent
(c) for all observation horizons $n \in \mathbb{N}$
\be
C_{\lambda,n}^{L} \ < \ \left\{
\begin{array}{ll}
B_{\lambda,n}^{L}, & \textrm{if } \qui \in \quasetSPcomp \times ]0,1[,\\ 
V_{\lambda,n}, & \textrm{if } \qui \in (\quasetNI \cup \quasetSPeins) \times ]0,1[.
\end{array}
\right.
\notag
\ee 
(d) 
\be
\lim_{n\rightarrow\infty} C_{\lambda,n}^{L} \ = \ \left\{
\begin{array}{ll}
0, & \textrm{if } \qui \in (\quasetNI)^{c} \times ]0,1[,\\
\exp\Big\{
\omega_{0} \, x_{0}^{(q_{\lambda}^{E})} 
\Big\} > 0, & \textrm{if } \qui \in \quasetNI  \times ]0,1[,
\end{array}
\right.
\notag
\ee
which coincides with $\lim_{n \rightarrow \infty} \,H_{\lambda}(\Pna||\Pnh)$ for all 
$\qui \in \left(\quaset\backslash \quasetSPvier \right) \times ]0,1[$.
\be
\hspace{-3.6cm} 
\textrm{(e)} \quad
\lim_{n \rightarrow \infty} \, \frac{1}{n} \, \log C_{\lambda,n}^{L} \ = \ 
\frac{p_{\lambda}^{\bigstar}}{q_{\lambda}^{\bigstar}} \, \left(x_0^{(q_{\lambda}^{\bigstar})} + \beta_{\lambda}\right)
- \alpha_{\lambda} \ , \quad \textrm{which coincides with}
\notag
\ee 
 $\lim_{n \rightarrow \infty} \, \frac{1}{n} \, \log H_{\lambda}(\Pna||\Pnh)$ for all 
$\qui \in \left(\quasetNI \cup \quasetSPeins \right) \times ]0,1[$
respectively with $\lim_{n \rightarrow \infty} \, \frac{1}{n} \, \log B_{\lambda,n}^{L}$ for all 
$\qui \in \left(\quasetSPzwei \cup \quasetSPdreiab 
\cup \quasetSPdreic \right)\times ]0,1[$.
\end{thm}

\begin{rem}
\label{remCLLOW}
Notice that the formula \eqref{fo.ansumbnLCLO} simplifies in the parameter case 
$\qui $ $\in \quasetSPeins \times ]0,1[$, since then it holds $p_{\lambda}^{\bigstar}/q_{\lambda}^{\bigstar}=\frac{\aal}{\bal}=\frac{\ahy}{\bhy}$ and therewith $(p_{\lambda}^{\bigstar}/q_{\lambda}^{\bigstar})\cdot\beta_{\lambda}-\alpha_{\lambda}=0$;
for the case $\qui \in \quasetNI \times ]0,1[$, one can even use the stronger relation $p_{\lambda}^{\bigstar} = 0 = \alpha_{\lambda}$.
\end{rem}

\noindent
In order to get an 
``explicit'' lower bound which does not rely on the implicitly given fixed point $x_{0}^{(q_{\lambda}^{\bigstar})}$, one can replace the latter by a 
close explicit lower approximate $\underline{x}_{0}^{(q_{\lambda}^{\bigstar})} < x_{0}^{(q_{\lambda}^{\bigstar})}$ and proceed completely analogously,
leading to a smaller lower bound (say) $\underline{C}_{\lambda,n}^{L} < C_{\lambda,n}^{L}$ in assertion (a) of Theorem \ref{thm3varLOW};
in the corresponding assertions (b), (c) and (d) one then has to replace $C_{\lambda,n}^{L}$ by $\underline{C}_{\lambda,n}^{L}$
and $x_{0}^{(q_{\lambda}^{E})}$ by $\underline{x}_{0}^{(q_{\lambda}^{E})}$.
For instance, one could choose 
\be\label{defux0}
\underline{x}_{0}^{(q_{\lambda}^{\bigstar})} \ := \ \frac{e^{-
h(q_{\lambda}^{\bigstar})  
}}{q_{\lambda}^{\bigstar}}\cdot\left[\left(1-q_{\lambda}^{\bigstar}\right)-\sqrt{\left(1-q_{\lambda}^{\bigstar}\right)^{2}-2\cdot q_{\lambda}^{\bigstar}\cdot e^{h(q_{\lambda}^{\bigstar})}\cdot\left(q_{\lambda}^{\bigstar}-\beta_{\lambda}\right)}\right],
\ee
where
\be 
h(q_{\lambda}^{\bigstar}) \ := \ \left\{
\begin{array}{ll}
\max\left\{-\beta_{\lambda} \ ; \ \frac{q_{\lambda}^{\bigstar}-\beta_{\lambda}}{1-q_{\lambda}^{\bigstar}}\right\}, & \textrm{if} \ q_{\lambda}^{\bigstar} \ < \ 1,\\
-\beta_{\lambda}, & \textrm{if} \ q_{\lambda}^{\bigstar} \ \geq \ 1;
\end{array}
\right. 
\notag
\ee
this
will be used 
as an auxiliary tool for the diffusion-limit-concerning proof
of Lemma \ref{lem4}(c) in the appendix.
If $q_{\lambda}^{\bigstar}<1$, the term $\frac{q_{\lambda}^{\bigstar}-\beta_{\lambda}}{1-q_{\lambda}^{\bigstar}}$ represents the existing negative intersection of the tangent of $\xi_{\lambda}^{(q_{\lambda}^{\bigstar})}$ at $x=0$ and the bisectrix. Clearly, $h(q_{\lambda}^{\bigstar}) < x_{0}^{(q_{\lambda}^{\bigstar})}$.
By \eqref{defux0}, $\underline{x}_{0}^{(q_{\lambda}^{\bigstar})}$ is the unique negative solution of $\underline{Q}_{\lambda}^{(q_{\lambda}^{\bigstar})}(x)=x$ with the quadratic function
\be 
\underline{Q}_{\lambda}^{(q_{\lambda}^{\bigstar})}(x) \ := \ \frac{q_{\lambda}^{\bigstar}}{2} \, e^{h(q_{\lambda}^{\bigstar})}\cdot x^{2}+q_{\lambda}^{\bigstar}\cdot x+q_{\lambda}^{\bigstar}-\beta_{\lambda}.
\notag
\ee
Notice that $\underline{Q}_{\lambda}^{(q_{\lambda}^{\bigstar})}(0)=\xi_{\lambda}^{(q_{\lambda}^{\bigstar})}(0)$,
$\frac{\textrm{d}\underline{Q}_{\lambda}^{(q_{\lambda}^{\bigstar})}}{\textrm{d}x}\, (0) = \frac{\textrm{d}\xi_{\lambda}^{(q_{\lambda}^{\bigstar})}}{\textrm{d}x}\, (0)$,
$\frac{\textrm{d}^2\underline{Q}_{\lambda}^{(q_{\lambda}^{\bigstar})}}{\textrm{d}x^2}\, (x) < \frac{\textrm{d}^2\xi_{\lambda}^{(q_{\lambda}^{\bigstar})}}{\textrm{d}x^2}\, (x)$ for all 
$x\in[x_{0}^{(q_{\lambda}^{\bigstar})},0]$,
and thus
$\underline{Q}_{\lambda}^{(q_{\lambda}^{\bigstar})}(x) < \xi_{\lambda}^{(q_{\lambda}^{\bigstar})}(x)$ for all $x\in[x_{0}^{(q_{\lambda}^{\bigstar})},0[$,
which leads to the desired $\underline{x}_{0}^{(q_{\lambda}^{\bigstar})} < x_{0}^{(q_{\lambda}^{\bigstar})}$.




\subsection{Closed-form upper bounds}
\label{secCLOSUP}

In order to achieve closed-form upper bounds, we principially proceed as in the previous
Section \ref{secCLOSLOW}. 
However, the situation is now more diverse since we have to start from Section \ref{secDETUP}
which carries much more ``nonuniqueness'' respectively variety than the corresponding
Sections \ref{secDETEX} and \ref{secDETLOW} which we used as a starting point
for the investigations in Section \ref{secCLOSLOW}.\\[0.1cm]

\noindent Notice first that for the subcases $\quasetSPdreiab \times ]0,1[$  
and $\quasetSPdreic \times ]0,1[$
(cf.\ Subsections \ref{secDETUP}(a2),(a3),(a4)) one can achieve a closed-form upper
bound without further investigations: if one chooses $q_{\lambda}^{U}=\beta_{\lambda}$
(and thus, the slope $s_{\lambda}^{U}=0$),
then by properties (p-i), (p-v) one has $a^{(q)}_{n} \equiv 0$ (i.e.\ recursion \eqref{defan} is trivial),
$b_{n}^{(p,q)} \equiv  \, p \, - \, \alpha_{\lambda} <0$ and hence
\be
\label{fo.HI.horizontal}
B_{\lambda,n}^{U} \ = \ 
\exp\Big\{
a^{(q_{\lambda}^{U})}_{n}\, \omega_{0} \, + \, 
\sum_{k=1}^{n} b^{(p_{\lambda}^{U},q_{\lambda}^{U})}_{k}
\Big\} 
\ = \ 
\exp\left\{n \cdot (p_{\lambda}^{U} \, - \, \alpha_{\lambda})   \right\}
\ \stackrel{n\rightarrow\infty}{\longrightarrow} \ 0 \ .
\ee
However, there might exist (and for $\quisetSPdreia$ definitely exists) choices
$(p_{\lambda}^{U},q_{\lambda}^{U})$ which lead
to (fully or eventually partially) tighter upper bounds $B_{\lambda,n}^{U}$ but for which the non-linear recursion \eqref{defan}
is nontrivial. Such potential cases, for which in particular $0<q_{\lambda}^{U}<\beta_{\lambda}$
and $0<p_{\lambda}^{U}\leq\alpha_{\lambda}$ holds, will be treated in the following;
since the parameter constellation $\qui \in (\quasetSPdreid \cup \quasetSPvier) \times ]0,1[$
does not meet this requirement,
let us fix 
$\qui \in (\quasetNI \cup \quasetSPeins \cup \quasetSPzwei \cup \quasetSPdreiab 
\cup \quasetSPdreic) \times ]0,1[$
where we also include the two setups $\quasetNI \cup \quasetSPeins$ for which we want to replace the recursive, non-closed-form exact values
by closed-form upper bounds.
For 
this situation, we determined recursive upper bounds respectively exact values in a (more or less) identical structure
which is also very close to the one given by \eqref{defanvar} to \eqref{fo.zeroSTAR}:
choose $q_{\lambda}^{G}$ for $G=U$ respectively $G=E$ subject to the corresponding parameter case which 
leads to $a^{(q_{\lambda}^{G})}_{1} = s_{\lambda}^{G} = 
q_{\lambda}^{G} - \beta_{\lambda} <0$,
compute the (rest of the) sequence $\left(a_{n}^{(q_{\lambda}^{G})}\right)_{n\in\mathbb{N}_{0}}$  
by the nonlinear recursion (cf.\ \eqref{defan}, \eqref{defxi})
\bea
a^{(q_{\lambda}^{G})}_{0}
:= 0 &;& \qquad a^{(q_{\lambda}^{G})}_{n}
\ := \
\xi^{(q_{\lambda}^{G})}_{\lambda}\hspace{-0.12cm}\left(a^{(q_{\lambda}^{G})}_{n-1}\right) 
,~~n\in\mathbb{N}, \label{defanUP}
\eea 
choose $p_{\lambda}^{G}$ 
subject to the corresponding parameter case and evaluate 
\be
b^{(p_{\lambda}^{G},q_{\lambda}^{G})}_{n} \ = \ 
\frac{p_{\lambda}^{G}}{q_{\lambda}^{G}} \, a_{n}^{(q_{\lambda}^{G})} \, + \, \frac{p_{\lambda}^{G}}{q_{\lambda}^{G}} \, \beta_{\lambda} \, - \, \alpha_{\lambda},~~n\in\mathbb{N}, 
\qquad \qquad \textrm{(cf.\ \eqref{fo.anbn})}
\notag
\ee
which leads to the desired bound 
$B_{\lambda,n}^{G}  =  
\exp\Big\{
a^{(q_{\lambda}^{G})}_{n}\, \omega_{0} \, + \, 
\sum_{k=1}^{n} b^{(p_{\lambda}^{G},q_{\lambda}^{G})}_{k}
\Big\} \ $ (cf.\ part (b) of Proposition \ref{propUPa1}).
According to (p-ii), the fundamentally important sequence 
$\left(a^{(q_{\lambda}^{G})}_{n}\right)_{n\in\mathbb{N}}$ is strictly negative, strictly
decreasing and converges to the unique solution 
$x_0^{(q_{\lambda}^{G})} \in ]-\infty,a^{(q_{\lambda}^{G})}_{1}[$ of the equation 
\be
\xi^{(q_{\lambda}^{G})}_{\lambda}(x)~=~ q_{\lambda}^{G}\cdot e^x - \beta_{\lambda} = x , 
\quad x < 0 \ .
\qquad \qquad \textrm{(cf.\ \eqref{fo.zeroUPa1}, \eqref{fo.zero})}
\notag
\ee
For an upper bound of the sequence $a_{n}^{(q_{\lambda}^{G})}$ we introduce the recursion
\bea
\oa_{0}^{(q_{\lambda}^{G})}
:= 0 &;& \qquad \oa_{n}^{(q_{\lambda}^{G})}
\ := \
\xi_{\lambda}^{(q_{\lambda}^{G}),S}\hspace{-0.12cm}\left(\oa^{(q_{\lambda}^{G})}_{n-1}\right) \ + \ \orho_{n-1}^{(q_{\lambda}^{G})}
,~~n\in\mathbb{N}, \label{defoa}
\eea 
i.e.\ we replace  
the nonlinear function $\xi^{(q_{\lambda}^{G})}_{\lambda}(x)~=~q_{\lambda}^{G}\cdot e^{x}-\beta_{\lambda}$
by the secant line of $\xi^{(q_{\lambda}^{G})}_{\lambda}$ across its arguments $x_0^{(q_{\lambda}^{G})}$ and $0$, defined by
\be\label{fo.xiS} 
\xi_{\lambda}^{(q_{\lambda}^{G}),S}(x)\ := \ 
q_{\lambda}^{G}-\beta_{\lambda} \ + \ 
\frac{x_{0}^{(q_{\lambda}^{G})}-(q_{\lambda}^{G}-\beta_{\lambda})}{x_{0}^{(q_{\lambda}^{G})}}
\cdot x
\, , \qquad x\in[x_{0}^{(q_{\lambda}^{G})},0] \ ,
\ee
and reduce the error we face by adding the ``correction-term'' 
\be\label{deforho}
\orho_{n-1}^{(q_{\lambda}^{G})} \ := \ -\frac{1}{2} \cdot\left(x_{0}^{(q_{\lambda}^{G})}\right)^{2}
\cdot\left(q_{\lambda}^{G} \cdot e^{x_{0}^{(q_{\lambda}^{G})}} \right)^{n}\cdot\left(1-\left(
\frac{x_{0}^{(q_{\lambda}^{G})}-(q_{\lambda}^{G}-\beta_{\lambda})}{x_{0}^{(q_{\lambda}^{G})}}
\right)^{n-1}\right) \ < \ 0.
\ee  
In other words, by means of \eqref{fo.dTundGamma} and the function on the domain $[0,\infty[$
\be
q \ \mapsto \  d^{(q),S} := 
\frac{x_{0}^{(q)}-(q-\beta_{\lambda})}{x_{0}^{(q)}} 
\label{fo.dS}
\ee
we use \eqref{fo.defanCLO}, \eqref{fo.xiCLO}, \eqref{defrho}
with the constants  
$d:= d^{(q_{\lambda}^{G}),S}  \in ]d^{(q_{\lambda}^{G}),T},1[$,  $c:= 
q_{\lambda}^{G}-\beta_{\lambda}  \in ]-\infty,0[$, 
\ $K_{1}:= - \Gamma^{(q_{\lambda}^{G})} <0$,  
\ $\varkappa:= d^{(q_{\lambda}^{G}),T} \in ]0,d[$,
\ $K_{2}:= -K_{1}$,
\ $\nu:= d^{(q_{\lambda}^{G}),T} \cdot d^{(q_{\lambda}^{G}),S} \in ]0,\varkappa[$. 
The following fundamental properties will be proved in the appendix:

\begin{lem}\label{lemoa}
For all $\qui \in \quaset\backslash\{\quasetSPdreid\cup\quasetSPvier\}\times]0,1[$ it holds
\bea
& (a) & \oa_{n}^{(q_{\lambda}^{G})} \ \geq \ a_{n}^{(q_{\lambda}^{G})}, \qquad \textrm{for all} \ n\in\mathbb{N}, \ \textrm{with equality iff} \ n=1.\hspace{5cm}\notag\\
& (b) & \textrm{The sequence} \ \left(\oa_{n}^{(q_{\lambda}^{G})}\right)_{n\in\mathbb{N}} \ \textrm{is strictly decreasing}.\notag\\
& (c) & \lim_{n\rightarrow\infty}\oa_{n}^{(q_{\lambda}^{G})} \ = \ \lim_{n\rightarrow\infty}a_{n}^{(q_{\lambda}^{G})} \ = \ x_{0}^{(q_{\lambda}^{G})}.\notag
\eea
\end{lem}

\noindent
Applying Theorem \ref{thm2}, Lemma \ref{lemoa} as well as the formulae \eqref{fo.anbn}, \eqref{fo.repanCLO} and \eqref{fo.repsumanCLO}, one obtains



\begin{thm}\label{thm3varUP}
For all $\qui \in (\quasetNI \cup \quasetSPeins \cup \quasetSPzwei \cup \quasetSPdreiab 
\cup \quasetSPdreic) \times ]0,1[$ 
and  all initial population sizes $\omega_{0}\in\mathbb{N}$ 
the following assertions hold:\\
\noindent
(a) for all observation horizons $n \in \mathbb{N}$ the Hellinger integral can be bounded from above by the closed-form bounds $H_{\lambda}(\Pna||\Pnh) \ < \ C_{\lambda,n}^{G}$ given by
\bea\label{fo.oanCLO}
 C_{\lambda,n}^{G} && := \ \exp\left\{ x_{0}^{(q_{\lambda}^{G})}\left[\omega_{0}-\frac{p_{\lambda}^{G}}{q_{\lambda}^{G}}\frac{d^{(q_{\lambda}^{G}),S}}{1-d^{(q_{\lambda}^{G}),S}}\right]\left(1-\left(d^{(q_{\lambda}^{G}),S}\right)^{n}\right) \ + \ \left[\frac{p_{\lambda}^{G}}{q_{\lambda}^{G}}\left(\beta_{\lambda}+x_{0}^{(q_{\lambda}^{G})}\right)-\alpha_{\lambda}\right]\cdot n\right.\notag\\
&& \hspace{1.4cm}\left.- \  \overline{\zeta}^{(q_{\lambda}^{G})}_{n}\cdot\omega_{0} \ - \ \overline{\vartheta}^{(q_{\lambda}^{G})}_{n}\right\}, \qquad \textrm{where for all} \ n\in\mathbb{N} 
\eea

\be\label{def.ozeta}
\overline{\zeta}^{(q_{\lambda}^{G})}_{n} \ := \ \Gamma^{(q_{\lambda}^{G})}\left[\frac{\left(d^{(q_{\lambda}^{G}),S}\right)^{n}-\left(d^{(q_{\lambda}^{G}),T}\right)^{n}}{d^{(q_{\lambda}^{G}),S}-d^{(q_{\lambda}^{G}),T}}-\frac{\left(d^{(q_{\lambda}^{G}),S}\right)^{n}\left(1-\left(d^{(q_{\lambda}^{G}),T}\right)^{n}\right)}{d^{(q_{\lambda}^{G}),S}\left(1-d^{(q_{\lambda}^{G}),T}\right)}\right] \ > \ 0 \qquad \textrm{and}
\ee

\be\label{def.ovartheta}
\overline{\vartheta}^{(q_{\lambda}^{G})}_{n}  :=   \frac{\Gamma^{(q_{\lambda}^{G})} \cdot p_{\lambda}^{G}\cdot d^{(q_{\lambda}^{G}),T}}{q_{\lambda}^{G}(1-d^{(q_{\lambda}^{G}),T})}\left[\frac{1-\left(d^{(q_{\lambda}^{G}),S}d^{(q_{\lambda}^{G}),T}\right)^{n}}{1-d^{(q_{\lambda}^{G}),S}d^{(q_{\lambda}^{G}),T}}+\frac{1-\left(d^{(q_{\lambda}^{G}),S}\right)^{n}}{d^{(q_{\lambda}^{G}),S}-d^{(q_{\lambda}^{G}),T}}-\frac{1-\left(d^{(q_{\lambda}^{G}),T}\right)^{n}}{d^{(q_{\lambda}^{G}),S}-d^{(q_{\lambda}^{G}),T}}\right]  >  0 \ .
\ee
The parameters $0<q_{\lambda}^{G}<\beta_{\lambda}$, $0<p_{\lambda}^{G}\leq\alpha_{\lambda}$ can be chosen 
subject to the restrictions explained in the parameter-adequate Subsections \ref{secDETUP}(a1),(a2),(a3),(a4) (for $G=U$)
respectively Subsections \ref{secDETEX}(aNI),(aEF) (for $G=E$).\\
\noindent
(b)
the sequence $\left(C_{\lambda,n}^{G}\right)_{n\in\mathbb{N}}$ is strictly decreasing.\\
\noindent
(c) for all observation horizons $n \in \mathbb{N}$
\be
C_{\lambda,n}^{G} \ \geq \ \left\{
\begin{array}{ll}
B_{\lambda,n}^{U}, & \textrm{if } \qui \in (\quasetSPzwei \cup \quasetSPdreiab 
\cup \quasetSPdreic) \times ]0,1[,\\ 
V_{\lambda,n}, & \textrm{if } \qui \in (\quasetNI \cup \quasetSPeins) \times ]0,1[ \, ,
\end{array}
\right.
\notag
\ee 
with equality iff $n=1$.\\
(d) 
\be
\lim_{n\rightarrow\infty}C_{\lambda,n}^{G} \ = \ \left\{
\begin{array}{ll}
0, & \textrm{if } \qui \in \quasetSP\backslash(\quasetSPdreid\cup\quasetSPvier)
 \times ]0,1[,\\
\exp\Big\{
\omega_{0} \, x_{0}^{(q_{\lambda}^{E})} 
\Big\} > 0, & \textrm{if } \qui \in \quasetNI  \times ]0,1[,
\end{array}
\right.
\notag
\ee
\be
\hspace{-3.9cm} = \ \lim_{n \rightarrow \infty} \,H_{\lambda}(\Pna||\Pnh). \hspace{4.7cm} \ 
\notag
\ee 
\be
\hspace{-3.6cm} 
\textrm{(e)} \quad
\lim_{n \rightarrow \infty} \, \frac{1}{n} \, \log C_{\lambda,n}^{G} \ = \ 
\frac{p_{\lambda}^{G}}{q_{\lambda}^{G}} \, \left(x_0^{(q_{\lambda}^{G})} + \beta_{\lambda}\right)
- \alpha_{\lambda} \qquad \textrm{which coincides with}
\notag
\ee 
$\lim_{n \rightarrow \infty} \, \frac{1}{n} \, \log H_{\lambda}(\Pna||\Pnh)$ for all 
$\qui \in \left(\quasetNI \cup \quasetSPeins \right) \times ]0,1[$
repectively with $\lim_{n \rightarrow \infty} \, \frac{1}{n} \, \log B_{\lambda,n}^{U}$ for all 
$\qui \in \left(\quasetSPzwei \cup \quasetSPdreiab 
\cup \quasetSPdreic \right) \times ]0,1[$.
\end{thm}

\noindent 
Notice that the strict positivity in \eqref{def.ozeta} and \eqref{def.ovartheta} can be easily seen from \eqref{focntil}.

\begin{rem}
\label{remCLUP}
The formula \eqref{fo.oanCLO} simplifies in the parameter case $\qui \in \quasetSPeins $ $\times ]0,1[$
which results in 
$\frac{p_{\lambda}^{E}}{q_{\lambda}^{E}}\left(\beta_{\lambda}+x_{0}^{(q_{\lambda}^{E})}\right)-\alpha_{\lambda}
 = 
 \frac{\aal}{\bal}
 \cdot x_{0}^{(q_{\lambda}^{E})}$.
For the case $\qui \in \quasetNI \times ]0,1[$, one can even use the stronger relation $p_{\lambda}^{E} = 0 = \alpha_{\lambda}$.
\end{rem}

\noindent
In order to get an 
``explicit'' upper bound which does not rely on the implicitly given fixed point $x_0^{(q_{\lambda}^{G})}$ $(G\in\{U,E\})$, one can replace the latter by a 
close explicit upper approximate $\overline{x}_{0}^{(q_{\lambda}^{G})} > x_{0}^{(q_{\lambda}^{G})}$ and proceed completely analogously,
leading to a larger upper bound (say) $\overline{C}_{\lambda,n}^{G} > C_{\lambda,n}^{G}$ in assertion (a) of Theorem \ref{thm3varUP};
in the corresponding assertions (b), (c) and (d) one then has to replace $C_{\lambda,n}^{G}$ by $\overline{C}_{\lambda,n}^{G}$
and $x_{0}^{(q_{\lambda}^{E})}$ by $\overline{x}_{0}^{(q_{\lambda}^{E})}$.
One possibility along these lines 
is the choice
\be\label{defox0}
\overline{x}_{0}^{(q_{\lambda}^{G})} \ := \ \frac{1}{q_{\lambda}^{G}}\cdot\left[\left(1-q_{\lambda}^{G}\right)-\sqrt{\left(1-q_{\lambda}^{G}\right)^{2}-2\cdot q_{\lambda}^{G}\cdot\left(q_{\lambda}^{G}-\beta_{\lambda}\right)} \, \right],
\ee
which is exactly the unique negative solution to the quadratic equation
\be 
\overline{Q}_{\lambda}^{(q_{\lambda}^{G})}(x) \ := \ \frac{q_{\lambda}^{G}}{2}\cdot x^{2}+q_{\lambda}^{G}\cdot x+q_{\lambda}^{G}-\beta_{\lambda} \ = \ x.
\notag
\ee
By inspection of the first two derivatives, one gets $\overline{Q}_{\lambda}^{(q_{\lambda}^{G})}(x) > \xi_{\lambda}^{(q_{\lambda}^{G})}(x)$ for all $x\in[-\infty,0[$, and thus $\overline{x}_{0}^{(q_{\lambda}^{G})} > x_{0}^{(q_{\lambda}^{G})}$.
Notice that the additional fundamental requirement $\overline{x}_{0}^{(q_{\lambda}^{G})} < q_{\lambda}^{G}-\beta_{\lambda}$ 
holds for parameter constellations for which $\underline{x}_{0}^{(q_{\lambda}^{G})} > -1$ for any $\underline{x}_{0}^{(q_{\lambda}^{G})} \leq x_{0}^{(q_{\lambda}^{G})}$, since then one has
$\frac{\textrm{d}\overline{Q}_{\lambda}^{(q_{\lambda}^{G})}}{\textrm{d}x}\, (x) > 0$ for all $x \in [\underline{x}_{0}^{(q_{\lambda}^{G})},0]$.
Such a situation will be used as an auxiliary tool for the proof of Lemma \ref{lem4}(c) in the appendix. 


\section{Hellinger integral bounds in the diffusion limit}
\label{sec.diflim}
One can show that a properly rescaled Galton-Watson process with immigration (GWI) 
converges weakly to a diffusion process $\widetilde{X} := 
\left\{\widetilde{X}_t \, , t \in [0,\infty[ \right\}$ which is the unique, strong, nonnegative
-- and in case of 
$\frac{\eta}{\sigma^{2}}\geq\frac{1}{2}$
strictly positive --  solution
of the stochastic differential equation (SDE) of the form
\be
d\widetilde{X}_{t} \ = \ \left(\eta \, - \, \kappa \, \widetilde{X}_{t}\right)\, dt \, + \, \sigma\sqrt{\widetilde{X}_{t}} \, dW_{t}, \ \quad  t \in [0,\infty[, \qquad \widetilde{X}_{0} \in ]0,\infty[  \textrm{ given},
\label{def.CIR-SDE}
\ee
where $\eta \in [0,\infty[$, $\kappa \in [0,\infty[$,
$\sigma \in ]0,\infty[$ are constants 
and $W_{t}$ denotes a standard Brownian motion with respect to the underlying probability measure $P$;
see e.g.\ 
Feller \cite{Fel51}, Jirina \cite{Jir69}, Lamperti \cite{Lam67a}, \cite{Lam67b},  Lindvall \cite{Lin72}, \cite{Lin74},
Grimvall \cite{Gri74}, Jagers \cite{Jag75}, Borovkov \cite{Bor86}, Ethier and Kurtz \cite{EK86}, Durrett \cite{Dur96}
for the 
non-immigration case 
corresponding to $\eta = 0$, $\kappa \geq 0$, 
Kawazu and Watanabe \cite{Kaw71}, Wei and Winnicki \cite{Wei89}, Winnicki \cite{Win91} for the immigration case corresponding to $\eta \ne 0$, $\kappa=0$,
as well as Sriram \cite{Sri94} for the general case $\eta \in[0,\infty[$, $\kappa \in \mathbb{R}$. 
Feller-type branching processes of the form \eqref{def.CIR-SDE}, which are special cases of continuous state branching processes
with immigration (see e.g.\ Kawazu und Watanabe \cite{Kaw71}, Li \cite{Li06}, as well as Dawson and Li \cite{Daw06} for
imbeddings to affine processes) play for instance an important role in the modelling of the term structure of interest rates, 
cf.\ the seminal Cox-Ingersoll-Ross CIR model 
\cite{CIR85b} and the vast follow-up literature thereof. 
Furthermore, \eqref{def.CIR-SDE} is also prominently used as (a special case of) Cox and Ross's \cite{Co76}
constant elasticity of variance CEV asset price process, as
(part of) Heston's \cite{Hes93} stochastic asset-volatility framework, 
as a model of neuron activity (see e.g.\ Lansky and Lanska \cite{Lan87}, Giorno et al.\ \cite{Gio88},
Lanska et al.\ \cite{Lan94}, Lansky et al.\ \cite{Lan95}, Ditlevsen and Lansky \cite{Dit06}, Höpfner \cite{Hoe07},
Lansky and Ditlevsen \cite{Lan08}),
as a time-dynamic description of the nitrous oxide emission rate from the soil surface
(see e.g.\ Pedersen \cite{Ped00}),
as well as a model for the individual hazard rate in a survival analysis context (see e.g.\ Aalen and Gjessing \cite{Aal04}).\\[0.1cm]

\noindent
Along these lines 
of branching-type diffusion limits, it makes sense to consider the solutions of two SDEs \eqref{def.CIR-SDE} with different
fixed parameter sets $(\eta,\kappa_{\Al},\sigma)$ and $(\eta,\kappa_{\Hy},\sigma)$, determine for each of them 
a corresponding approximating GWI, investigate the Hellinger integral between the laws of these two GWI, and finally
calculate the Hellinger integral (bounds) limit as the GWI approach their SDE solutions.
Notice that for technicality reasons (which will be explained below), the constants
$\eta$ and $\sigma$ ought to be independent of $\Al$, $\Hy$ in our current context.
\ \\[0.05cm]

\noindent
In order to make the abovementioned limit procedure rigorous, it is reasonable to 
work with appropriate approximations such that in each convergence step $m$ one 
faces the setup $(\quasetNI \cup \quasetSPeins) \times ]0,1[$
(i.e.\ the non-immigration or the equal-fraction case), where the corresponding Hellinger integral can be calculated
exactly in a recursive way (cf.\ Section \ref{secDETEX}). Let us explain the details in the following.\\[0.1cm]
\noindent
Consider a sequence of GWI $\left(X^{(m)}\right)_{m\in\mathbb{N}}$ with probability laws $P^{(m)}_{\bullet}$ on a measurable space $(\Omega,\mathcal{A})$, 
where
as above the subscript $\bullet$ stands for either the hypothesis $\Hy$ or 
the alternative $\Al$. Analogously to \eqref{033}, we use for each fixed step $m \in \mathbb{N}$ the representation
$X^{(m)}:= \left\{ X^{(m)}_{n}, \, n  \in \mathbb{N} \right\}$ with
\vspace{-0.3cm}
\be\label{def1}
X^{(m)}_{n}~:=~\sum_{k=1}^{X^{(m)}_{n-1}}Y^{(m)}_{n-1,k}+\widetilde{Y}^{(m)}_{n}, \qquad n \in \mathbb{N},   
\qquad X^{(m)}_{0} \in \mathbb{N} \textrm{ given},
\ee
where under the law $P^{(m)}_{\bullet}$

\begin{itemize}

\item the collection $Y^{(m)}:= \left\{ Y^{(m)}_{n-1,k}, \, n  \in \mathbb{N}, k  \in \mathbb{N} \right\}$
consists of i.i.d.\ 
random variables which are
Poisson distributed with parameter $\beta^{(m)}_{\bullet} >0$,

\item the collection $\widetilde{Y}^{(m)}:= \left\{ \widetilde{Y}^{(m)}_{n}, \, n  \in \mathbb{N} \right\}$
consists of i.i.d. random variables which are
Poisson distributed with parameter $\alpha^{(m)}_{\bullet} \geq 0$,

\item $Y^{(m)}$ and $\widetilde{Y}^{(m)}$ are independent.

\end{itemize}

\noindent
From arbitrary drift-parameters $\eta \in [0,\infty[$, $\kappa_{\bullet} \in [0,\infty[$, 
and diffusion-term-parameter $\sigma>0$, we construct the 
offspring-distribution-parameter and the immigration-distribution parameter
of the sequence $\left(X^{(m)}_{n}\right)_{n\in\mathbb{N}}$ by 
\be\label{defparm}
\beta^{(m)}_{\bullet} \, := \, 1-\frac{\kappa_{\bullet}}{\sigma^{2}m} \qquad \text{and} \qquad \alpha^{(m)}_{\bullet} \, := \, \beta^{(m)}_{\bullet}\cdot\frac{\eta}{\sigma^{2}}.
\ee
Here and henceforth, we always assume that the approximation step $m$ is large enough to ensure 
that $\beta^{(m)}_{\bullet}\in]0,1]$ and at least one of $\bal^{(m)}$, $\bhy^{(m)}$ is strictly less than 1; this will be abbreviated by $m\in \largeint$. 
Let us point out that 
-- as mentioned above -- our choice entails the best-to-handle setup
$(\quasetNI \cup \quasetSPeins) \times ]0,1[$ (which does not happen if instead of $\eta$ one uses $\eta_{\bullet}$ 
with $\eta_{\Al} \ne \eta_{\Hy}$).
Based on the GWI $X^{(m)}$, let us construct the continuous-time branching process
$\widetilde{X}^{(m)}:= \left\{ \widetilde{X}^{(m)}_{t}, \, t  \in [0, \infty[ \right\}$ by
\vspace{-0.3cm}
\be\label{def2}
\widetilde{X}^{(m)}_{t}~:=~\frac{1}{m} \,  X^{(m)}_{\left\lfloor \sigma^{2}mt\right\rfloor} \ ,
\ee
living on the state space $E^{(m)}:=\frac{1}{m}\mathbb{N}_{0}$. 
From \eqref{def2} one can see immediately the necessity of having $\sigma$ to be independent of $\Al$, $\Hy$ because for the
required absolute continuity in \eqref{def.Zn} both models at stake have to ``live''
on the same time-scale $\tau_{t}^{(m)}:=\left\lfloor \sigma{^2}mt\right\rfloor$.
For this setup, one obtains the following convergence result:



\begin{thm}\label{thm4}
Let $\eta\in[0,\infty[$, $\kappa_{\bullet}\in [0,\infty[$, 
$\sigma \in]0,\infty[$ and $\widetilde{X}^{(m)}$ 
be as defined in~(\ref{def1}) to (\ref{def2}). 
Furthermore,  let us suppose that $\lim_{m\rightarrow \infty} \frac{1}{m} \, X^{(m)}_{0} = \widetilde{X}_{0} >0$ and denote by 
$D([0,\infty[, [0,\infty[)$ the space of right-continuous functions $f: [0,\infty[ \mapsto [0,\infty[$ with left limits. 
Then the sequence of processes
$\left(\widetilde{X}^{(m)}\right)_{m \in \largeint}$ convergences in distribution
in $D([0,\infty[, [0,\infty[)$ to a diffusion process $\widetilde{X}$
which is the unique strong, nonnegative -- and in case of 
$\frac{\eta}{\sigma^{2}}\geq\frac{1}{2}$
strictly positive -- solution of the SDE
\be\label{fo.diffusion}
d\widetilde{X}_{t} \ = \ \big(\eta \, - \, \kappa_{\bullet} \, \widetilde{X}_{t}\big) \, dt \ + \ \sigma\sqrt{\widetilde{X}_{t}}~dW_{t}^{\bullet},
\quad t \in [0,\infty[, \qquad \widetilde{X}_{0} \in ]0,\infty[ \textrm{ given},
\ee
where $W_{t}^{\bullet}$ denotes a standard Brownian motion with respect to the 
limit probability measure $P_{\bullet}$.
\end{thm}

\indent Notice that the condition $\frac{\eta}{\sigma^{2}}\geq\frac{1}{2}$ can be interpreted in our approximation setup \eqref{defparm}
as 
$\alpha^{(m)}_{\bullet} \geq \beta^{(m)}_{\bullet}/2$, which quantifies the intuitively reasonable 
indication  
that if
the probability $P_{\bullet}[\widetilde{Y}^{(m)}_{n}=0] = e^{-\alpha^{(m)}_{\bullet}}$ 
of having no immigration is small enough relative
to the probability $P_{\bullet}[Y^{(m)}_{n-1,k}=0]=e^{-\beta^{(m)}_{\bullet}}$ of having no offspring
($m \in \largeint$),
then the limiting diffusion $\widetilde{X}$ never hits zero almost surely.\\[0.1cm]

\noindent
The corresponding proof of Theorem \ref{thm4} -- which is outlined 
in Appendix \ref{App5} -- is  
an adaption of the proof of Theorem 9.1.3 in Ethier and Kurtz \cite{EK86}
which deals with drift-parameters $\eta=0$, $\kappa_{\bullet}=0$ in the SDE \eqref{fo.diffusion} whose solution 
is approached on a $\sigma-$independent time scale
by a sequence of (critical) Galton-Watson processes without immigration but with
general offspring distribution with mean $1$ and variance $\sigma$.
Notice that due to \eqref{defparm} the latter is inconsistent with our Poissonian setup,
but this is compensated by our chosen $\sigma-$dependent time scale. Furthermore, \eqref{defparm} is also inconsistent
with the other concrete parameter choices in the abovementioned corresponding references.
\\[0.1cm]

\noindent
For each approximation step $m$ and each observation horizon $\tau_{t}^{(m)}$, 
the corresponding Hellinger integrals $H_{\lambda}\left(\Pmat\Big|\Big|\Pmht\right)$ obey the
results of  

\begin{enumerate}

\item[(ap1)] the Propositions \ref{propNI} (for $\eta=0$) and
\ref{propPSP1} (for $\eta \in]0,\infty[$), as far as recursively computable exact values are concerned,

\item[(ap2)] Theorem \ref{thm3varLOW} 
and
Theorem \ref{thm3varUP}, 
as far as closed-form bounds are concerned;
recall that the current setup is of type $(\quasetNI \cup \quasetSPeins) \times ]0,1[$,
and thus we can use the simplifications proposed in the Remarks
\ref{remCLLOW} and \ref{remCLUP}.

\end{enumerate}

\noindent
In order to obtain the desired limits of 
the Hellinger integrals $H_{\lambda}\left(\Pmat\Big|\Big|\Pmht\right)$
respectively of their closed-form bounds as $m\rightarrow \infty$, one faces the following problems: 
in accordance with Section \ref{secDETEX}, for each fixed $m$ in (ap1) one has to choose the parameters 
$p_{\lambda}^{(m)}:=\left(\aam\right)^{\lambda}\left(\ahm\right)^{1-\lambda}$,
$q_{\lambda}^{(m)}:=\left(\bam\right)^{\lambda}\left(\bhm\right)^{1-\lambda}$,
which in particular 
determine the fundamental 
sequence $\left(a_{n}^{(m)}\right)_{n\in\mathbb{N}}:=\left(a_{n}^{(q_{\lambda}^{(m)})}\right)_{n\in\mathbb{N}}$
(cf.\ \eqref{defan}). 
This enters in the appropriate versions of part (b) in Propositions \ref{propNI} and \ref{propPSP1} respectively  
in form of  $a_{\tau_{t}^{(m)}}^{(q_{\lambda}^{(m)})}$ , and the correspondingly arising convergences (as $m\rightarrow \infty$) seem to
be not (straightforwardly) tractable due to the recursive nature of \eqref{defan}. 
In contrast, for the \textit{closed-form bounds} of Section \ref{secCFB} the desired
convergences 
are tractable, which will be worked out in the following. To begin with, let us explicitly
formulate the results of the application of 
Theorem \ref{thm3varLOW} (where Remark \ref{remCLLOW} applies) 
and Theorem \ref{thm3varUP} (where Remark \ref{remCLUP} applies) to the current setup.
For this, we use the following SDE-parameter constellations
(which are consistent with \eqref{defparm} in combination with our requirement 
to work here only on $(\quasetNI \cup \quasetSPeins) \times ]0,1[$):
let $\widetilde{\mathcal{P}}_{NI}$ be the set of all 
$(\kappa_{\Al},\kappa_{\Hy},\eta)$  for which 
$\eta=0$, $\kappa_{\Al} \in [0,\infty[$, $\kappa_{\Hy} \in [0,\infty[$
with $\kappa_{\Al} \ne \kappa_{\Hy}$; furthermore,
denote by $\widetilde{\mathcal{P}}_{SP,1}$
the set of all 
$(\kappa_{\Al},\kappa_{\Hy},\eta)$  for which 
$\eta \in ]0,\infty[$, $\kappa_{\Al} \in [0,\infty[$, $\kappa_{\Hy} \in [0,\infty[$
with $\kappa_{\Al} \ne \kappa_{\Hy}$. On $\widetilde{\mathcal{P}}_{NI} \cup \widetilde{\mathcal{P}}_{SP,1}$
there hold for $m \in \largeint$ the useful restrictions $q_{\lambda}^{(m)} \in ]0,1[$ and $\beta_{\lambda}^{(m)}\in]0,1[$. 
For the sake of brevity, let us henceforth use the abbreviations
$\alpha_{\lambda}^{(m)} := \lambda\cdot\aam+(1-\lambda)\cdot\ahm$, \   
$\beta_{\lambda}^{(m)} := \lambda\cdot\bam+(1-\lambda)\cdot\bhm$,
$x_{0}^{(m)}  := x_{0}^{(q_{\lambda}^{(m)})}$, \ 
$\Gamma^{(m)}  :=   \Gamma^{(q^{(m)}_{\lambda})}  =   \frac{q^{(m)}_{\lambda}}{2}\cdot e^{x_{0}^{(m)}}\cdot\left(x_{0}^{(m)}\right)^{2}$, \ 
$d^{(m),S} := d^{(q_{\lambda}^{(m)}),S} = 
\frac{x_{0}^{(m)}-(q_{\lambda}^{(m)}-\beta_{\lambda}^{(m)})}{x_{0}^{(m)}}$ 
and 
$ d^{(m),T} \ := \ d^{(q_{\lambda}^{(m)}),T} \ = \ q_{\lambda}^{(m)}\cdot e^{x_{0}^{(m)}}.
$
By the above considerations, the Theorems \ref{thm3varLOW} and \ref{thm3varUP}
(together with their remarks) adapt to the current setup as follows:



\begin{cor}\label{cor1} \ \hspace{-0.4cm} 
For all $(\kappa_{\Al},\kappa_{\Hy},\eta,\lambda) \in (\widetilde{\mathcal{P}}_{NI}\cup\widetilde{\mathcal{P}}_{SP,1}) \times ]0,1[$, all $t\in[0,\infty[$, all approximation steps $m \in\largeint$ and all initial population sizes 
$X_{0}^{(m)} \in\mathbb{N}$ 
the Hellinger integral can be bounded by
\bea
&& \exp\Bigg\{x_{0}^{(m)}\cdot\left[X_{0}^{(m)}-\frac{\eta}{\sigma^{2}}\frac{d^{(m),T}}{1-d^{(m),T}}\right]\left(1-\left(d^{(m),T}\right)^{\left\lfloor \sigma^{2}mt\right\rfloor}\right)\notag\\
&& \hspace{2.0cm}+ \ x_{0}^{(m)}\frac{\eta}{\sigma^{2}}\cdot\left\lfloor \sigma^{2}mt\right\rfloor \ + \ \underline{\zeta}^{(m)}_{\left\lfloor \sigma^{2}mt\right\rfloor} \cdot X_{0}^{(m)} \ + \ \underline{\vartheta}^{(m)}_{\left\lfloor \sigma^{2}mt\right\rfloor}\Bigg\}
\label{fo.bouapprox1} \\
& \leq & H_{\lambda}\left(\Pmat\Big|\Big|\Pmht\right) \notag\\
& \leq & \exp\Bigg\{x_{0}^{(m)}\cdot\left[X_{0}^{(m)}-\frac{\eta}{\sigma^{2}}\frac{d^{(m),S}}{1-d^{(m),S}}\right]\left(1-\left(d^{(m),S}\right)^{\left\lfloor \sigma^{2}mt\right\rfloor}\right)\notag\\
&& \hspace{2.0cm}+ \ x_{0}^{(m)}\frac{\eta}{\sigma^{2}}\cdot\left\lfloor \sigma^{2}mt\right\rfloor \ - \ \overline{\zeta}^{(m)}_{\left\lfloor \sigma^{2}mt\right\rfloor}\cdot X_{0}^{(m)} \ - \ \overline{\vartheta}^{(m)}_{\left\lfloor \sigma^{2}mt\right\rfloor}\Bigg\} \, ,
\label{fo.bouapprox}
\eea
where 
we define analogously to 
\eqref{def.uzeta}, \eqref{def.uvartheta}, \eqref{def.ozeta} and \eqref{def.ovartheta} 
\bea 
\underline{\zeta}^{(m)}_{n} & := & \Gamma^{(m)}\cdot\frac{\left(d^{(m),T}\right)^{n-1}\left(1-\left(d^{(m),T}\right)^{n}\right)}{1-d^{(m),T}} \  > 0 \, , 
\notag\\
\underline{\vartheta}^{(m)}_{n} & := & \frac{\eta}{\sigma^{2}}\cdot\Gamma^{(m)}\cdot\frac{\left(1-\left(d^{(m),T}\right)^{n}\right)}{\left(1-d^{(m),T}\right)^{2}}\cdot\left(1-d^{(m),T}\cdot\frac{1+\left(d^{(m),T}\right)^{n}}{1+d^{(m),T}}\right) \geq0 \, ,\notag\\
\overline{\zeta}^{(m)}_{n} & := & \Gamma^{(m)}\cdot\left[\frac{\left(d^{(m),S}\right)^{n}-\left(d^{(m),T}\right)^{n}}{d^{(m),S}-d^{(m),T}}-\frac{\left(d^{(m),S}\right)^{n-1}\left(1-\left(d^{(m),T}\right)^{n}\right)}{1-d^{(m),T}}\right] >0 \, ,\notag\\
\overline{\vartheta}^{(m)}_{n} & := & \frac{\eta}{\sigma^{2}}\cdot\Gamma^{(m)}\cdot\frac{d^{(m),T}}{1-d^{(m),T}}\cdot\left[\frac{1-\left(d^{(m),S}d^{(m),T}\right)^{n}}{1-d^{(m),S}d^{(m),T}}-\frac{\left(d^{(m),S}\right)^{n}-\left(d^{(m),T}\right)^{n}}{d^{(m),S}-d^{(m),T}}\right] \geq0.
\notag
\eea
Notice that the bounds \eqref{fo.bouapprox} simplify significantly in the case $(\kappa_{\Al},\kappa_{\Hy},\eta,\lambda) \in \widetilde{\mathcal{P}}_{NI} \times ]0,1[$ for which $\eta=0$ holds.
\end{cor}

\noindent
Let us finally present the corresponding desired limit assertions as the approximation step
$m$ tends to infinity,  
by making use of the quantities
\be\label{def.kl-Lambda}
\kappa_{\lambda} \ := \ \lambda\kappa_{\Al}+(1-\lambda)\kappa_{\Hy} \ > \ 0
 \qquad \textrm{as well as}\qquad  \Lambda_{\!\lambda} \ := \ \sqrt{\lambda\kappa_{\Al}^{2}+(1-\lambda)\kappa_{\Hy}^{2}} \  > \ \kappa_{\lambda} \ :
\ee



\begin{thm}\label{thmlimit}
Let the initial SDE-value  $\widetilde{X}_{0}\in ]0,\infty[$  be arbitrary but fixed, and suppose that $\lim_{m\rightarrow \infty} \frac{1}{m} \, X^{(m)}_{0} = \widetilde{X}_{0}$. 
Then, for all  $t\in[0,\infty[$ 
and all $(\kappa_{\Al},\kappa_{\Hy},\eta,\lambda) \in (\widetilde{\mathcal{P}}_{NI}\cup\widetilde{\mathcal{P}}_{SP,1}) \times ]0,1[$
the Hellinger integral limit can be bounded by
\bea 
&& D_{\lambda,t}^{L} \ := \  \exp\Bigg\{-\frac{\Lambda_{\!\lambda}-\kappa_{\lambda}}{\sigma^{2}}\left[\widetilde{X}_{0}-\frac{\eta}{\Lambda_{\!\lambda}}\right]\left(1-e^{-\Lambda_{\!\lambda}\cdot t}\right)-\frac{\eta}{\sigma^{2}}\left(\Lambda_{\!\lambda}-\kappa_{\lambda}\right)\cdot t
\notag\\
&& \hspace{2.5cm} 
\ + \ L_{\lambda}^{(1)}(t)\cdot\widetilde{X}_{0} \ + \ \frac{\eta}{\sigma^{2}}\cdot L_{\lambda}^{(2)}(t)\Bigg\} 
\\[-0.2cm]
& &
\leq \ \lim_{m\rightarrow\infty} H_{\lambda}\left(\Pmat\Big|\Big|\Pmht\right)\label{fo.boulim1} \notag\\
& &  \leq \  \exp\Bigg\{-\frac{\Lambda_{\!\lambda}-\kappa_{\lambda}}{\sigma^{2}}\left[\widetilde{X}_{0}-\frac{\eta}{\frac{1}{2}(\Lambda_{\!\lambda}+\kappa_{\lambda})}\right]\left(1-e^{-\frac{1}{2}(\Lambda_{\!\lambda}+\kappa_{\lambda})\cdot t}\right)-\frac{\eta}{\sigma^{2}}\left(\Lambda_{\!\lambda}-\kappa_{\lambda}\right)\cdot t\notag\\
&& \hspace{1.5cm}- \ U_{\lambda}^{(1)}(t)\cdot\widetilde{X}_{0} \ - \  \frac{\eta}{\sigma^{2}}\cdot U_{\lambda}^{(2)}(t)\Bigg\}
\ =: \ D_{\lambda,t}^{U} \ ,
\label{fo.boulim} 
\eea
\ \\[-0.7cm]
where for all $t\geq0$
\bea
L_{\lambda}^{(1)}(t) & := & \frac{\left(\Lambda_{\!\lambda}-\kappa_{\lambda}\right)^{2}}{2\sigma^{2}\cdot\Lambda_{\!\lambda}}\cdot e^{-\Lambda_{\!\lambda}\cdot t}\cdot\left(1-e^{-\Lambda_{\!\lambda}\cdot t}\right) \ > \ 0,\label{def.Leins}\\
L_{\lambda}^{(2)}(t) & := & \frac{1}{4}\cdot\left(\frac{\Lambda_{\!\lambda}-\kappa_{\lambda}}{\Lambda_{\!\lambda}}\right)^{2}\cdot\left(1-e^{-\Lambda_{\!\lambda}\cdot t}\right)^{2} \ > \ 0,\label{def.Lzwei}\\
U_{\lambda}^{(1)}(t) & := & \frac{\left(\Lambda_{\!\lambda}-\kappa_{\lambda}\right)^{2}}{\sigma^{2}}\cdot\left[\frac{e^{-\frac{1}{2}(\Lambda_{\!\lambda}+\kappa_{\lambda})\cdot t}-e^{-\Lambda_{\!\lambda}\cdot t}}{\Lambda_{\!\lambda}-\kappa_{\lambda}}-\frac{e^{-\frac{1}{2}(\Lambda_{\!\lambda}+\kappa_{\lambda})\cdot t}\left(1-e^{-\Lambda_{\!\lambda}\cdot t}\right)}{2\cdot \Lambda_{\!\lambda}}\right] \ \geq \ 0,\label{def.Ueins}\\
U_{\lambda}^{(2)}(t) & := & \frac{\left(\Lambda_{\!\lambda}-\kappa_{\lambda}\right)^{2}}{\Lambda_{\!\lambda}}\cdot\left[\frac{1-e^{-\frac{1}{2}\left(3\Lambda_{\!\lambda}+\kappa_{\lambda}\right)\cdot t}}{3\Lambda_{\!\lambda}+\kappa_{\lambda}}+\frac{e^{-\Lambda_{\!\lambda}\cdot t}-e^{-\frac{1}{2}(\Lambda_{\!\lambda}+\kappa_{\lambda})\cdot t}}{\Lambda_{\!\lambda}-\kappa_{\lambda}}\right] \ \geq \ 0.\label{def.Uzwei}
\eea
Notice that the components $L_{\lambda}^{(i)}(t)$ and $U_{\lambda}^{(i)}(t)$ ($i=1,2$) do not depend on the parameter $\eta$, and that
the bounds \eqref{fo.boulim1} and \eqref{fo.boulim} simplify significantly in the case $(\kappa_{\Al},\kappa_{\Hy},\eta,\lambda) \in \widetilde{\mathcal{P}}_{NI} \times ]0,1[$, for which $\eta=0$ holds.
\end{thm}




\section{Power divergences and relative entropy}
\label{sec.ent}

All the results of the previous sections carry correspondingly over from the Hellinger integrals 
$H_{\lambda}( \cdot || \cdot )$ ($\lambda \in ]0,1[$) to the power divergences $I_{\lambda}\left( \cdot || \cdot \right)$
by virtue of the relation (cf.\ \eqref{fo.powerdivdef})
\be
I_{\lambda}\left(\Pa||\Ph \right) \ = \ 
\frac{1-H_{\lambda}(\Pa||\Ph)}{\lambda\cdot(1-\lambda)} \ .
\notag
\ee
In particular, this leads to bounds on $I_{\lambda}\left(\Pa||\Ph \right)$ which are tighter than the 
general rudimentary bound \eqref{fo.powbound} connected with \eqref{fo.hellgenup}.
Furthermore, it is well known that in general the relative entropy 
defined by \eqref{fo.relent}
\be\label{fo.liment}
I(\Pa||\Ph) \ = \ 
\lim_{\lambda\nearrow1} \ I_{\lambda}\left(\Pa||\Ph \right) \ ,
\ee
see e.g.\ Liese and Vajda \cite{LieVa87}. 
Accordingly, for our context of GWI we can use \eqref{fo.liment}
in combination with the \textit{recursive} exact values respectively \textit{recursive} 
lower bounds of Theorem \ref{thm2} 
and Section \ref{secDETLOW}
to obtain the following \textit{closed-form} exact values respectively \textit{closed-form} upper bounds of 
the relative entropy $I(\Pna||\Pnh)$:

\begin{thm}\label{thm.entex}
(a) For all $\qua \in (\quasetNI \cup \quasetSPeins)$, all initial population sizes $\omega_{0}\in\mathbb{N}$ and  all observation horizons $n \in \mathbb{N}$ 
\be\label{fo.entex}
I(\Pna||\Pnh) \ = \ \left\{
\begin{array}{ll}
\frac{\bal\cdot\left(\log\left(\frac{\bal}{\bhy}\right)-1\right)+\bhy}{1-\bal}\cdot\left[\omega_{0}-\frac{\aal}{1-\bal}\right]\cdot\left(1-\left(\bal\right)^{n}\right) & ~\\[0.2cm]
 \quad + \ \frac{\aal\cdot\left[\bal\cdot\left(\log\left(\frac{\bal}{\bhy}\right)-1\right)+\bhy\right]}{\bal(1-\bal)}\cdot n \, , & \textrm{if }\bal\neq1,\\
 ~&~\\[-0.1cm]
\left[\bhy -\log \bhy -1 \right]\cdot\left[\frac{\aal}{2}\cdot n^{2}+\left(\omega_{0}+\frac{\aal}{2}\right)\cdot n\right] \, , & \textrm{if }\bal=1.
\end{array}
\right.
\ee
(b) For all $\qua \in \quasetSPcompvar$ \hspace{-0.10cm} , all initial population sizes $\omega_{0}\in\mathbb{N}$ and  all observation hori\-zons $n \in \mathbb{N}$ it holds $ \ I(\Pna||\Pnh) \ \leq \ E_{n}^{U}$, where
\be\label{fo.entbouup}
E_{n}^{U} := \left\{
\begin{array}{ll}
\frac{\bal\cdot\left(\log\left(\frac{\bal}{\bhy}\right)-1\right)+\bhy}{1-\bal}\cdot\left[\omega_{0}-\frac{\aal}{1-\bal}\right]\cdot\left(1-\left(\bal\right)^{n}\right) & ~\\[0.2cm]
 \quad + \ \left[\frac{\aal\cdot\left[\bal\cdot\left(\log\left(\frac{\bal}{\bhy}\right)-1\right)+\bhy\right]}{\bal(1-\bal)}+\aal\left[\log\left(\frac{\aal\bhy}{\ahy\bal}\right)-\frac{\bhy}{\bal}\right]+\ahy\right]\cdot n \, ,& \textrm{if }\bal\neq1,\\
 ~&~\\[-0.2cm]
\left[\bhy -\log \bhy -1\right]\cdot\left[\frac{\aal}{2}\cdot n^{2}+\left(\omega_{0}+\frac{\aal}{2}\right)\cdot n\right] & \\[0.1cm]
 + \ \left[\aal\left[\log\left(\frac{\aal\bhy}{\ahy}\right)-\bhy\right]+\ahy \right]\cdot n \, , &  \textrm{if }\bal=1.
\end{array}
\right.
\ee
\end{thm}

\begin{rem}
The $n-$behaviour of (the bounds of) the relative entropy $I(\Pna||\Pnh)$ in Theorem \ref{thm.entex} 
is influenced by the following facts:
\begin{itemize}
\item  $\bal\cdot\left(\log\left(\frac{\bal}{\bhy}\right)-1\right)+\bhy\geq0$ with equality iff $\bal=\bhy$.
\item In the case $\bal\neq1$ of \eqref{fo.entbouup}, there holds $\frac{\aal\cdot\left[\bal\cdot\left(\log\left(\frac{\bal}{\bhy}\right)-1\right)+\bhy\right]}{\bal(1-\bal)}+\aal\left[\log\left(\frac{\aal\bhy}{\ahy\bal}\right)-\frac{\bhy}{\bal}\right]+\ahy\geq0$, with equality iff $\aal=\ahy$ and $\bal=\bhy$. 
\end{itemize}
\end{rem}
 
\noindent In contrast, in order to derive (semi-)\textit{closed-form} lower bounds of 
the relative entropy $I(\Pna||\Pnh)$ we use \eqref{fo.liment}
in combination with the \textit{recursive} 
upper bounds of Theorem \ref{thm2}(b) 
and appropriately adapted detailed analyses along the lines of Section \ref{secDETUP}.
This amounts to


\begin{thm}\label{thm.entexUP}
For all $\qua \in  \quasetSP\backslash \quasetSPeins$, 
all initial population sizes $\omega_{0}\in\mathbb{N}$ and  all observation hori\-zons $n \in \mathbb{N}$
\be\label{fo.entboulowSP}
I(\Pna||\Pnh) \ \geq \ E^{L}_{n} \ := \ \sup_{k\in\mathbb{N}_{0}  , ~ y\in[0,\infty[}\Big\{E^{L,tan}_{y,n} \, , \, E^{L,sec}_{k,n} \ , \ E^{L,hor}_{n} \Big\} \in [0,\infty[ \ ,
\ee
where for all $y\in[0,\infty[$ we define the 
-- possibly negatively valued -- 
finite bound component 
\be\label{fo.ELyn}
E^{L,\textrm{tan}}_{y,n} \ := \ \left\{
\begin{array}{ll}
\left[\bal\ln\left(\frac{\aal+\bal y}{\ahy+\bhy y}\right)+\bhy\left(1-\frac{\aal+\bal y}{\ahy+\bhy y}\right)\right]\cdot\frac{1-\left(\bal\right)^{n}}{1-\bal}\cdot\left[\omega_{0}-\frac{\aal}{1-\bal}\right] & \\[0.2cm]
+ \ \Big[\frac{\aal}{\bal(1-\bal)}\left[\bal\ln\left(\frac{\aal+\bal y}{\ahy+\bhy y}\right)+\bhy\left(1-\frac{\aal+\bal y}{\ahy+\bhy y}\right)\right] & \\
\hspace{0.6cm}+ \ \left(\ahy-\aal\frac{\bhy}{\bal}\right)\left(1-\frac{\aal+\bal y}{\ahy+\bhy y}\right)\Big]\cdot n \, , & \textrm{if} \ \bal\neq1,\vspace{0.2cm}\\
\left[\ln\left(\frac{\aal+y}{\ahy+\bhy y}\right)+\bhy\left(1-\frac{\aal+y}{\ahy+\bhy y}\right)\right]\cdot\left[\frac{\aal}{2}\cdot n^{2}+\left(\omega_{0}+\frac{\aal}{2}\right)\cdot n\right] & \\
+ \ \left(\ahy-\aal\bhy\right)\left(1-\frac{\aal+ y}{\ahy+\bhy y}\right)\cdot n \, , & \textrm{if} \ \bal=1,
\end{array}
\right.
\ee
and for all $k\in\mathbb{N}_{0}$ the 
-- possibly negatively valued -- 
finite bound component
\be\label{fo.ELkn}
E^{L,\textrm{sec}}_{k,n} \ := \ \left\{
\begin{array}{ll}
\left[f_{\Al}(k+1)\ln\left(\frac{f_{\Al}(k+1)}{f_{\Hy}(k+1)}\right)-f_{\Al}(k)\ln\left(\frac{f_{\Al}(k)}{f_{\Hy}(k)}\right)+\bhy-\bal\right]\cdot\frac{1-\left(\bal\right)^{n}}{1-\bal}\cdot\left[\omega_{0}-\frac{\aal}{1-\bal}\right] & \vspace{0.1cm}\\
+ \ \Big[\frac{\aal}{\bal(1-\bal)}\left(f_{\Al}(k+1)\ln\left(\frac{f_{\Al}(k+1)}{f_{\Hy}(k+1)}\right)-f_{\Al}(k)\ln\left(\frac{f_{\Al}(k)}{f_{\Hy}(k)}\right)+\bhy-\bal\right) & \vspace{0.1cm}\\
\hspace{0.5cm}-\left(f_{\Al}(k+1)\ln\left(\frac{f_{\Al}(k+1)}{f_{\Hy}(k+1)}\right)-f_{\Al}(k)\ln\left(\frac{f_{\Al}(k)}{f_{\Hy}(k)}\right)\right)\cdot\left(k+\frac{\aal}{\bal}\right) & \vspace{0.1cm}\\
\hspace{0.5cm}+f_{\Al}(k)\ln\left(\frac{f_{\Al}(k)}{f_{\Hy}(k)}\right)-\frac{\aal\bhy}{\bal}+\ahy\Big]\cdot n \, , & \hspace{-1.9cm}\textrm{if} \ \bal\neq1, \vspace{0.3cm}\\
\left[f_{\Al}(k+1)\ln\left(\frac{f_{\Al}(k+1)}{f_{\Hy}(k+1)}\right)-f_{\Al}(k)\ln\left(\frac{f_{\Al}(k)}{f_{\Hy}(k)}\right)+\bhy-1\right]\cdot\left[\frac{\aal}{2}\cdot n^{2}+\left(\omega_{0}+\frac{\aal}{2}\right)\cdot n\right] & \vspace{0.1cm}\\
-\Big[\left(f_{\Al}(k+1)\ln\left(\frac{f_{\Al}(k+1)}{f_{\Hy}(k+1)}\right)-f_{\Al}(k)\ln\left(\frac{f_{\Al}(k)}{f_{\Hy}(k)}\right)\right)\left(k+\aal\right) & \vspace{0.1cm}\\
\hspace{0.5cm}-f_{\Al}(k)\ln\left(\frac{f_{\Al}(k)}{f_{\Hy}(k)}\right)+\aal\bhy-\ahy\Big]\cdot n \, , & \hspace{-1.9cm}\textrm{if} \ \bal=1.
\end{array}
\right.
\ee
Furthermore, on $\quasetSPvier$ we set $E^{L,hor}_{n}:=0$ for all $n\in \mathbb{N}$ whereas  
on $\quasetSP\backslash(\quasetSPeins\cup\quasetSPvier)$ we define
\be\label{fo.ELhorn}
E^{L,hor}_{n} \ := \ \left[\left(\aal+\bal z^{*}\right)\cdot\left[\log\left(\frac{\aal+\bal z^{*}}{\ahy+\bhy z^{*}}\right)-1\right]+\ahy+\bhy z^{*}\right]\cdot n,
\qquad, n \in \mathbb{N},
\ee
with $z^{*} \ := \ \arg\max_{x\in\mathbb{N}_{0}}\left\{(\aal+\bal x)\left[ - \log\left(\frac{\aal+\bal x}{\ahy+\bhy x}\right)+1\right]-(\ahy+\bhy x)\right\}$. 
In the subcases $\quasetSPzwei \cup \quasetSPdreiab \cup \quasetSPdreic \cup \quasetSPvier$ one gets even $E^{L}_{n}>0$
for all $\omega_{0}\in\mathbb{N}$ and all $n \in \mathbb{N}$. In the subcase $\quasetSPdreid$, one obtains for each fixed $n\in \mathbb{N}$ and each fixed $\omega_{0}\in\mathbb{N}$ the strict positivity
$E^{L}_{n}>0$ if $\left(\frac{\partial}{\partial y}E^{L,tan}_{y,n}\right)(y^{*})\neq0$, 
where $y^{*}:=\frac{\aal-\ahy}{\bhy-\bal} \in \mathbb{N}$ and hence
\bea
\label{fo.deryELtanystar}&& \hspace{0.5cm}\left(\frac{\partial}{\partial y}E^{L,tan}_{y,n}\right)(y^{*}) 
\\
&& \hspace{-0.6cm} =  \left\{
\begin{array}{ll}
\hspace{-0.1cm}
-  \frac{(\bal-\bhy)^{3}}{\aal\bhy-\ahy\bal}\cdot\frac{1-\left(\bal\right)^{n}}{1-\bal}\cdot\left[\omega_{0}-\frac{\aal}{1-\bal}\right]  -  \frac{(\bal-\bhy)^{2}}{\bal}\left(1+\frac{\aal(\bal-\bhy)}{(1-\bal)(\aal\bhy-\ahy\bal)}\right)\cdot n \, ,  & \textrm{if} \ \bal\neq1,\vspace{0.2cm}\notag\\
- \ \frac{(1-\bhy)^{3}}{\aal\bhy-\ahy}\cdot\left[\frac{\aal}{2}\cdot n^{2}+\left(\omega_{0}+\frac{\aal}{2}\right)\cdot n\right] - (1-\bhy)^{2}\cdot n \, , & \textrm{if} \ \bal=1.
\end{array}
\right.
\eea
\end{thm}
\begin{rem}\label{rem.ex.ent}
Consider the exemplary parameter setup $\qua=(\frac{1}{3},\frac{2}{3},2,1)\in\quasetSPdreid$. For initial population $\omega_{0}=3$ it holds $\left(\frac{\partial}{\partial y}E^{L,tan}_{y,n}\right)(y^{*})=0$ for all $n\in\mathbb{N}$, whereas for $\omega_{0}\neq3$ one obtains $\left(\frac{\partial}{\partial y}E^{L,tan}_{y,n}\right)(y^{*})\neq0$ for all $n\in\mathbb{N}$.
\end{rem}

\vspace{0.3cm}
\indent
It seems 
that the optimzation problem in \eqref{fo.entboulowSP} admits in general only 
an implicitly representable solution. 
Of course, as a less tight but less involved explicit lower bound of the relative entropy $I(\Pna||\Pnh)$ 
one can use any term of the form $\max\left\{E^{L,tan}_{y,n} \, , \, E^{L,sec}_{k,n}  \ , \ E^{L,hor}_{n} \right\}$
($y\in[0,\infty[$, $k\in\mathbb{N}_{0}$), as well as the following

\begin{cor}\label{cor.entLB}
(a) 
For all $\qua \in  \quasetSP\backslash \quasetSPeins$, 
all initial population sizes $\omega_{0}\in\mathbb{N}$ and  all observation hori\-zons $n \in \mathbb{N}$
\be
I(\Pna||\Pnh) \ \geq \ E^{L}_{n} \ \geq  \ \widetilde{E^{L}_{n}} :=  \max\left\{E^{L,tan}_{\infty,n} \, , \, E^{L,sec}_{0,n}  \ , \ E^{L,hor}_{n} \right\} \in [0,\infty[ \ ,
\notag
\ee
with $E^{L,hor}_{n}$ defined by \eqref{fo.ELhorn}, with
-- possibly negatively valued -- 
finite bound component $ \ E^{L,tan}_{\infty,n}  \ := \ \lim_{y\rightarrow\infty}E^{L,tan}_{y,n}$, where
\be 
E^{L,tan}_{\infty,n}  \ := \  
\left\{
\begin{array}{ll}
\frac{\bal\cdot\left(\log\left(\frac{\bal}{\bhy}\right)-1\right)+\bhy}{1-\bal}\cdot\left[\omega_{0}-\frac{\aal}{1-\bal}\right]\cdot\left(1-\left(\bal\right)^{n}\right) & ~\\[0.2cm]
 + \ \left[\frac{\aal\cdot\left[\bal\cdot\left(\log\left(\frac{\bal}{\bhy}\right)-1\right)+\bhy\right]}{\bal(1-\bal)}+\aal\left(1-\frac{\bhy}{\bal}\right)+\ahy\left(1-\frac{\bal}{\bhy}\right)\right]\cdot n \, , & \textrm{if }\bal\neq1, \\
 ~&~\\
\left[\bhy -\log \bhy -1\right]\cdot\left[\frac{\aal}{2}\cdot n^{2}+\left(\omega_{0}+\frac{\aal}{2}\right)\cdot n\right] & \\
+ \ \left[ \aal\left(1-\bhy\right)+\ahy\left(1-\frac{1}{\bhy}\right)\right]\cdot n \, , & \textrm{if }\bal=1,
\end{array}
\right.
\notag
\ee
and  -- possibly negatively valued -- finite bound component 
\be
\label{fo.entboulowSP23ab}
E^{L,sec}_{0,n} \ = \ 
\left\{
\begin{array}{ll}
\left[\left(\aal+\bal\right)\cdot\log\left(\frac{\aal+\bal}{\ahy+\bhy}\right)-\aal\cdot\log\left(\frac{\aal}{\ahy}\right)+\bhy-\bal\right]\cdot\frac{1-\left(\bal\right)^{n}}{1-\bal}\cdot\left[\omega_{0}-\frac{\aal}{1-\bal}\right]
 & ~\\ +\bigg\{\frac{\aal}{\bal(1-\bal)}\left(\left(\aal+\bal\right)\cdot\log\left(\frac{\aal+\bal}{\ahy+\bhy}\right)-\aal\cdot\log\left(\frac{\aal}{\ahy}\right)\right)-\frac{\aal}{1-\bal}\left(1-\bhy\right) & ~ \\
 \quad -\aal\left(1+\frac{\aal}{\bal}\right)\cdot\log\left(\frac{\ahy(\aal+\bal)}{\aal(\ahy+\bhy)}\right)+\ahy\bigg\}\cdot n \, , & \hspace{-1.8cm}\textrm{if }\bal\neq1,\\
 ~&~\\
\left[\left(\aal+1\right)\cdot\log\left(\frac{\aal+1}{\ahy+\bhy}\right)-\aal\cdot\log\left(\frac{\aal}{\ahy}\right)+\bhy-1\right]\cdot\left[n\cdot\omega_{0}+\frac{\aal}{2}\cdot n^{2}\right] & \\
+\Big\{\frac{\aal}{2}\left[\left(\aal+1\right)\cdot\log\left(\frac{\aal+1}{\ahy+\bhy}\right)-\aal\cdot\log\left(\frac{\aal}{\ahy}\right)-\bhy-1\right] & ~\\[0.2cm]
\quad-\aal\left(1+\aal\right)\cdot\log\left(\frac{\ahy(\aal+1)}{\aal(\ahy+\bhy)}\right)+\ahy\Big\}\cdot n \, , & \hspace{-1.8cm}\textrm{if }\bal=1.
\end{array}
\right.
\ee
For the cases $\quasetSPzwei \cup \quasetSPdreiab 
\cup \quasetSPdreic$
one gets even $\widetilde{E^{L}_{n}} > 0$ for all $\omega_{0}\in\mathbb{N}$ and  all $n \in \mathbb{N}$.
\end{cor}

\vspace{0.3cm}

\noindent For the diffusion-limit of the relative entropy we obtain \textit{closed-form} exact values:

\begin{thm}\label{thm.entdiflim} \hspace{-0.4cm}
Within the framework of Section \ref{sec.diflim}, one gets for all initial 
SDE-values $\widetilde{X}_{0} \in ]0,\infty[$,
all observation horizons $t\in[0,\infty[$ and all parameter constellations $(\kappa_{\Al},\kappa_{\Hy},\eta)\in
(
\widetilde{\mathcal{P}}_{NI}\cup\widetilde{\mathcal{P}}_{SP,1}
)
$
\bea 
&& \lim_{m\rightarrow\infty} I\left(P^{(m)}_{\Al,\left\lfloor \sigma^{2}mt\right\rfloor}\Big|\Big|P^{(m)}_{\Hy,\left\lfloor \sigma^{2}mt\right\rfloor}\right) \ = \ 
\lim_{m\rightarrow\infty}\lim_{\lambda\nearrow1}I_{\lambda}\left(P^{(m)}_{\Al,\left\lfloor \sigma^{2}mt\right\rfloor}\Big|\Big|P^{(m)}_{\Hy,\left\lfloor \sigma^{2}mt\right\rfloor}\right) \notag\\
& & = \left\{\begin{array}{ll} \frac{\left(\kappa_{\Al}-\kappa_{\Hy}\right)^{2}}{2\sigma^{2}\cdot\kappa_{\Al}}\cdot\left[\left(\widetilde{X}_{0}-\frac{\eta}{\kappa_{\Al}}\right)\cdot\left(1-e^{-\kappa_{\Al}\cdot t}\right)+\eta\cdot t\right] \, , & \textrm{if} \ \kappa_{\Al}>0,\\
 & \\
\frac{\kappa_{\Hy}^{2}}{2\sigma^{2}}\cdot\left[\frac{\eta}{2}\cdot t^{2} \ + \ \widetilde{X}_{0}\cdot t\right] \, , & \textrm{if} \ \kappa_{\Al}=0,
\end{array}
\right. \notag\\
&& \ = \ \lim_{\lambda\nearrow1}\lim_{m\rightarrow\infty}I_{\lambda}\left(P^{(m)}_{\Al,\left\lfloor \sigma^{2}mt\right\rfloor}\Big|\Big|P^{(m)}_{\Hy,\left\lfloor \sigma^{2}mt\right\rfloor}\right) \ . \notag
\eea
\end{thm}




\section{Applications}
\label{sec.dec}
As already mentioned in the introduction, there are numerous
applications of both ingredients -- power divergences resp.\ Hellinger integrals resp.\ relative entropy
on the one hand and Galton-Watson branching processes with immigration on the other hand.
In order to indicate the concrete applicability of our combinating investigations, for the sake of brevity
we confine ourselves to some issues in the context of Bayesian decision making BDM and Neyman-Pearson testing NPT.
In BDM, we decide here 
between an action $d_{\Hy}$ ``associated with'' the (say) hypothesis law $P_{\Hy}$ and an action $d_{\Al}$ ``associated with''
the (say) alternative law $P_{\Al}$,
based on the sample path observation $\mathcal{X}_{n}:=\{X_{l}: \ l\in\{0,1,\ldots,n\} \, \}
$
of the GWI-generation-sizes up to observation horizon 
$n \in \mathbb{N}$. 
Following the lines of Stummer and Vajda \cite{StuVa07} (adapted to our branching process context),
for BDM let us consider as admissible decision rules  $\delta_{n}: \Omega_{n} \mapsto \{d_{\Hy} , d_{\Al} \}$
the ones generated by all path sets $G_{n} \in \Omega_{n}$ 
through
\be  
	\begin{array}{lrll}
	& \delta_{n}(\mathcal{X}_{n}) \ := \ \delta_{G_{n}}(\mathcal{X}_{n}) &  
	:= & \left\{\begin{array}{ll}
	d_{\Al}, & \textrm{if} \ \mathcal{X}_{n}\in G_{n},
	\\
	d_{\Hy}, & \textrm{if} \ \mathcal{X}_{n}\notin G_{n},
	\end{array}\right.
	\end{array} 
	\notag
\ee
as well as loss functions of the form
	\be\left(\begin{array}{ll}  \label{lossfu}
	L(d_{\Hy},\Hy) & L(d_{\Hy},\Al)\\
	L(d_{\Al},\Hy) & L(d_{\Al},\Al)
	\end{array}\right) \ := \ 
	\left(\begin{array}{ll}
	0 & L_{\Al}\\
  L_{\Hy} & 0
	\end{array}\right)
	\ee
with pregiven constants $L_{\Al}>0$, $L_{\Hy}>0$ (e.g.\ arising as bounds from 
quantities in worst-case scenarios); notice that in \eqref{lossfu}, $d_{\Hy}$ is assumed to be
a zero-loss action under $\Hy$ and $d_{\Al}$ a zero-loss action under $\Al$.
Per definition, the Bayes decision rule $\delta_{G_{n,\text{min}}}$ minimizes -- over $G_{n}$ --
the mean decision loss
	\bea \label{meandecloss}
	\mathcal{L}(\delta_{G_{n}}) & := & p_{\Hy}^{\text{prior}}\cdot 
	\mathcal{L}_{\Hy} \cdot  
	Pr\left(\delta_{G_{n}}(\mathcal{X}_{n}) = d_{\Al} \Big| \Hy  \right)
	\ + \ p_{\Al}^{\text{prior}}\cdot 
	 \mathcal{L}_{\Al} \cdot 
Pr\left(\delta_{G_{n}}(\mathcal{X}_{n}) = d_{\Hy} \Big| \Al  \right)  
	\notag\\
	 & = & p_{\Hy}^{\text{prior}}\cdot 
	 \mathcal{L}_{\Hy} \cdot 
	\Pnh(G_{n}) \ + \ p_{\Al}^{\text{prior}}\cdot 
	\mathcal{L}_{\Al}
	\cdot 
	\Pna (\Omega_{n}- G_{n})
	\eea
for given prior probabilities $p_{\Hy}^{\text{prior}} = Pr(\Hy) \in ]0,1[$ for $\Hy$ and  $p_{\Al}^{\text{prior}} := Pr(\Al) = 1- p_{\Hy}^{\text{prior}}$ for $\Al$.
As a side remark let us mention that, in a certain sense, the involved model (parameter) uncertainty expressed by the ``superordinate'' Bernoulli-type law 
$Pr=Bin(1,p_{\Hy}^{\text{prior}})$ can also be reinterpreted as a rudimentary static random environment caused
e.g.\ by a random Bernoulli-type external static force.
By straightforward calculations, one gets with \eqref{def.Zn} the minimzing path set
$G_{n,\text{min}} = 
\left\{ Z_{n} \geq \frac{p_{\Hy}^{\text{prior}}L_{\Hy}}{p_{\Al}^{\text{prior}}L_{\Al}} \right\}$
leading to the minimal mean decision loss, i.e.\ the Bayes risk, 
	\bea \label{bayrisk1}
\mathcal{R}_{n} \ := \ \min_{G_{n}} \mathcal{L}(\delta_{G_{n}}) \ = \ 	
		\mathcal{L}(\delta_{G_{n,\text{min}}}) & = & \int_{\Omega_{n}}\min\left\{p_{\Hy}^{\text{prior}}L_{\Hy},p_{\Al}^{\text{prior}}L_{\Al}~Z_{n}\right\}d\Pnh \ .
	\eea
Notice that -- by straightforward standard arguments --  
the \textit{alternative} decision procedure
\be
\text{take action $d_{\Al}$ (resp.\ $d_{\Hy}$) \ \  if \ \
$L_{\Hy}  \cdot p_{\Hy}^{\text{post}}\hspace{-0.05cm}(\mathcal{X}_{n})  \ \leq \text{(resp.\ >)} \ \ L_{\Al} \cdot p_{\Al}^{\text{post}}\hspace{-0.05cm}(\mathcal{X}_{n})$}
\notag
\ee
with posterior probabilities $p_{\Hy}^{\text{post}}\hspace{-0.05cm}(\mathcal{X}_{n}) := \frac{p_{\Hy}^{\text{prior}}}{(1-p_{\Hy}^{\text{prior}}) \cdot Z_{n}\hspace{-0.05cm}(\mathcal{X}_{n}) \, + \, p_{\Hy}^{\text{prior}}}  =: 1- p_{\Al}^{\text{post}}\hspace{-0.05cm}(\mathcal{X}_{n})$, leads 
exactly to the same actions as 
$\delta_{G_{n,\text{min}}}$. By adapting Lemma 6.5 of Stummer and Vajda \cite{StuVa07}, one gets for 
all $L_{\Hy} >0$, $L_{\Al}>0$, $p_{\Hy}^{\text{prior}} \in ]0,1[$, $\lambda \in ]0,1[$ and $n \in \mathbb{N}$
the upper bound
\be
\label{upbayesrisk}
\mathcal{R}_{n} \ \leq \ \Lambda_{\Al}^{\lambda} \ \Lambda_{\Hy}^{1-\lambda} \ 
H_{\lambda}\left(\Pna||\Pnh\right)
\ ,
\qquad \textrm{with } \Lambda_{\Hy}:=p_{\Hy}^{prior}L_{\Hy},~\Lambda_{\Al}:=(1-p_{\Hy}^{prior})L_{\Al},
\ee
as well as the lower bound
\be
\left(\mathcal{R}_{n}\right)^{\min\{\lambda,1-\lambda\}} 
\cdot
\left( \Lambda_{\Hy} + \Lambda_{\Al}
-\mathcal{R}_{n}\right)^{\max\{\lambda,1-\lambda\}} 
\ \geq \ \Lambda_{\Al}^{\lambda} \ \Lambda_{\Hy}^{1-\lambda} \ 
H_{\lambda}\left(\Pna||\Pnh\right)
\
\notag
\ee
which implies in particular the ``direct'' lower bound
\be
\label{lowbayesrisk2}
\mathcal{R}_{n} \ \geq \ 
\frac{\Lambda_{\Al}^{\max\{1,\frac{\lambda}{1-\lambda}\}} \, \Lambda_{\Hy}^{\max\{1,\frac{1-\lambda}{\lambda}\}}}{\left(\Lambda_{\Al}+\Lambda_{\Hy}\right)^{\max\{\frac{\lambda}{1-\lambda},\frac{1-\lambda}{\lambda}\}}}\cdot\left(
H_{\lambda}\left(\Pna||\Pnh\right)
\right)^{\max\{\frac{1}{\lambda},\frac{1}{1-\lambda}\}} \ .
\ee
By using \eqref{upbayesrisk} (respectively \eqref{lowbayesrisk2}) together with the exact values and the upper (respectively lower) bounds of 
the Hellinger integrals $H_{\lambda}\left(\Pna||\Pnh\right)$
derived in the preceding sections, we end up with upper (respectively lower) bounds of the Bayes risk $\mathcal{R}_{n}$.
For different types of -- mainly parameter estimation (squared-error type loss function) concerning  --  Bayesian analyses 
based on GW(I) generation size observations, see e.g.\ 
Jagers \cite{Jag75}, Heyde \cite{Hey79}, Heyde and Johnstone \cite{HeyJoh79}, Johnson et al.\ \cite{Joh79}, Basawa and Rao \cite{Bas80}, 
Basawa and Scott \cite{Bas83}, Scott \cite{Sco87}, Guttorp \cite{Gut91}, Yanev and Tsokos \cite{Yan99},
Mendoza and Gutierrez-Pena \cite{Men00}, and the references therein. \\[0.2cm]

\noindent
Alternatively to the BDM applications above, let us now briefly deal with the corresponding NPT framework
with randomized tests $\mathcal{T}_{n}:  \Omega_{n} \mapsto [0,1]$ of the hypothesis $P_{\Hy}$ against
the alternative $P_{\Al}$,
based on the GWI-generation-size sample path observations $\mathcal{X}_{n}:=\{X_{l}: \ l\in\{0,1,\ldots,n\} \, \}
$. 
In contrast to \eqref{meandecloss}, \eqref{bayrisk1} a Neyman-Pearson test minimizes
-- over $\mathcal{T}_{n}$ -- the type II error probability
$\int_{\Omega_{n}} (1-\mathcal{T}_{n}) \, \text{d}\Pna$
in the class of the tests for which the type I error probability  
$\int_{\Omega_{n}} \mathcal{T}_{n} \, \text{d}\Pnh$
is at most
 $\varsigma \in ]0,1[$. The corresponding minimal type II error probability
\be 
 \mathcal{E}_{\varsigma}\left(\Pna||\Pnh\right) \ := \ 
\inf_{\mathcal{T}_{n}: 
\int_{\Omega_{n}} \hspace{-0.15cm}\mathcal{T}_{n} \, \text{d}\Pnh \leq  \varsigma}  
\  \int_{\Omega_{n}} (1-\mathcal{T}_{n}) \, \text{d}\Pna 
\notag
\ee
can for all $\varsigma \in ]0,1[$, $\lambda \in ]0,1[$, $n \in \mathbb{N}$ 
be bounded from above by 
\be \label{seconderrorbound}
 \mathcal{E}_{\varsigma}\left(\Pna||\Pnh\right) \ \leq \ \min \left\{ \ (1-\lambda) \cdot
\left(\frac{\lambda}{\varsigma} \right)^{\hspace{-0.05cm}\lambda/(1-\lambda)} \cdot
 \Big(\,  H_{1-\lambda}\left(\Pna||\Pnh\right) \, \Big)^{1/(1-\lambda)} \ , \ 1 \ \right\} \ 
\ee
which is an adaption of a general result of Krafft and Plachky \cite{Kra70}, see also Liese and Vajda \cite{LieVa87}
as well as Stummer and Vajda \cite{StuVa07}.
Hence, by combining \eqref{seconderrorbound} with the exact values respectively upper bounds of 
the Hellinger integrals $H_{1-\lambda}\left(\Pna||\Pnh\right)$
from the preceding sections, we obtain for our context of GWI with Poisson offspring and
Poisson immigration (including the non-immigration case) some upper bounds of 
$\mathcal{E}_{\varsigma}\left(\Pna||\Pnh\right)$, which
can also be immediately rewritten as lower bounds for the power
$1-\mathcal{E}_{\varsigma}\left(\Pna||\Pnh\right)$ of a most powerful test at level $\varsigma$.
In contrast to such finite-time-horizon results, for the (to our context)
incompatible setup of GWI with Poisson offspring but nonstochastic immigration of constant value 1,
the asymptotic rates of decrease as $n\rightarrow \infty$ of the unconstrained
type II  error probabilities as well as the
type I error probabilites 
were studied in Linkov and Lunyova \cite{Link96}
by a different approach employing also Hellinger integrals.
Some other types of GW(I) concerning Neyman-Pearson testing investigations different to ours can be found e.g.\ in 
Basawa and Scott \cite{Bas76}, Feigin \cite{Fei78}, Sweeting \cite{Swe78}, Basawa and Scott \cite{Bas83}, and the references therein.
\\[0.0cm]
       
\noindent
For the sake of brevity, a further more detailed 
discussion of GWI statistical issues along the lines of this section
as well as power-divergences-connected goodness-of-fit investigations
will appear in a forthcoming paper. \\[0.0cm]



\begin{appendix}

\section{Proofs and auxiliary lemmas}
\subsection{Tool and proof for Section \ref{secDET}}
\label{App3}
\begin{lem}\label{lem2}
\noindent
For all real numbers $x,y,z>0$ and all $\lambda\in ]0,1[$  one has
\be
x^{\lambda}y^{1-\lambda} \ - \ \left(\lambda\, x\, z^{\lambda-1}\, + \, (1-\lambda)\, y \, z^{\lambda}\right) \ \leq \ 0 \notag
\ee
with equality iff $\frac{x}{y}=z$.
\end{lem}
\prl\ref{lem2} \ For fixed $\tilde{x}:=xz^{\lambda-1}>0$, $\tilde{y}:=yz^{\lambda}>0$ with $\tilde{x} \ne \tilde{y}$
we inspect the function $g$ on $[0,1]$ defined by $g(\lambda) := \ \tilde{x}^{\lambda}\tilde{y}^{1-\lambda}-(\lambda \tilde{x} +(1-\lambda)\tilde{y})$
which satisfies $g(0)=g(1)=0$, $g'(0)=\tilde{y}\log(\tilde{x}/\tilde{y})-(\tilde{x}-\tilde{y})~<~\tilde{y}((\tilde{x}/\tilde{y})-1)-(\tilde{x}-\tilde{y})~=~0$
and which is strictly convex. Thus, the assertion follows immediately by taking into account the obvious case $\tilde{x} = \tilde{y}$.
\qed\\

\vspace{0.2cm}

\noindent
\prf\eqref{fo.hellnull3d}: \ For the parameter constellation in Subsection \ref{secDETUP}(a5) we
employ as upper bound for $\phi_{\lambda}(x)$, $x\in\mathbb{N}_{0}$ 
the function
\be
\overline{\phi_{\lambda}}(x):=\left\{
\begin{array}{ll}
\phi_{\lambda}(0), & \textrm{if } x=0,\\
0, & \textrm{if } x>0.
\end{array}
\right.
\notag
\ee 
Notice that this method is rather crude, and gives in the other cases treated in the Subsections 
\ref{secDETUP}(a1) to (a4) worse bounds than those derived there. For the calculation of the Hellinger integral, we first set $\epsilon:=1-e^{\phi_{\lambda}(0)}\in ]0,1[$. 
Hence, we obtain for all $n\in\mathbb{N}\backslash\{1\}$ 
\bea 
&& \sum_{\omega_{n-1}=0}^{\infty}\frac{\left[\varphi_{\lambda}(\omega_{n-2})\right]^{\omega_{n-1}}}{\omega_{n-1}!}\cdot \exp\{\phi_{\lambda}(\omega_{n-1})\}
~\leq~
\sum_{\omega_{n-1}=0}^{\infty}\frac{\left[\varphi_{\lambda}(\omega_{n-2})\right]^{\omega_{n-1}}}{\omega_{n-1}!}\cdot \exp\{\overline{\phi_{\lambda}}(\omega_{n-1})\}\notag\\
&& =~\exp\{\varphi_{\lambda}(\omega_{n-2})\}-\epsilon
~=~\exp\{\varphi_{\lambda}(\omega_{n-2})\} \cdot \left[1-\epsilon \cdot \exp\{-\varphi_{\lambda}(\omega_{n-2})\} \right] \notag\\
&&\leq~
\exp\left\{\varphi_{\lambda}(\omega_{n-2})-\epsilon\cdot e^{-\varphi_{\lambda}(\omega_{n-2})}\right\}.
\notag
\eea
In the current setup of Subsection \ref{secDETUP}(a5) we have $\bal\ne\bhy$, which means that $\lim_{x\rightarrow\infty}\phi_{\lambda}(x)=-\infty$
(cf.\ (p-xiii)). But this together with the nonnegativity of $\varphi_{\lambda}$ implies $\sup_{x\in\mathbb{N}_{0}}\exp\{\phi_{\lambda}(x)-\epsilon\cdot e^{-\varphi_{\lambda}(x)}\}=:\delta<1$. Incorporating these considerations as well as 
the formulae \eqref{fo.Znk} to \eqref{defflambda}, we get
for $n=1$ the relation  $H_{\lambda}\left(\Pna||\Pnh\right) = \exp\{\phi_{\lambda}(\omega_{n-1})\} < 1$
and for all $n\in\mathbb{N}\backslash\{1\}$  as a continuation of formula \eqref{fo.hellzw} (with the obvious shortcut for $n=2$)
\bea\label{fo.neuemethode}
&& H_{\lambda}\left(\Pna||\Pnh\right)
~=~\sum_{\omega_{1}=0}^{\infty}\cdots\sum_{\omega_{n}=0}^{\infty}\prod_{k=1}^{n}Z^{(\lambda)}_{n,k}(\omega) \notag \\
&&=~\sum_{\omega_{1}=0}^{\infty}\cdots\sum_{\omega_{n-1}=0}^{\infty}\prod_{k=1}^{n-1}Z^{(\lambda)}_{n,k}(\omega)
\notag\\
&&\hspace{1.5cm}\cdot\exp\Big\{\left(f_{\Al}(\omega_{n-1})\right)^{\lambda}\left(f_{\Hy}(\omega_{n-1})\right)^{(1-\lambda)}-(\lambda f_{\Al}(\omega_{n-1})+(1-\lambda)f_{\Hy}(\omega_{n-1}))\Big\}\notag\\
&& =~\sum_{\omega_{1}=0}^{\infty}\cdots\sum_{\omega_{n-2}=0}^{\infty}\prod_{k=1}^{n-2}Z^{(\lambda)}_{n,k}(\omega)\cdot\exp\left\{-f_{\lambda}(\omega_{n-2})\right\}\sum_{\omega_{n-1}=0}^{\infty}\frac{\left[\varphi_{\lambda}(\omega_{n-2})\right]^{\omega_{n-1}}}{\omega_{n-1}!}\cdot \exp\{\phi_{\lambda}(\omega_{n-1})\}\notag\\
&& \leq~\sum_{\omega_{1}=0}^{\infty}\cdots\sum_{\omega_{n-2}=0}^{\infty}\prod_{k=1}^{n-2}Z^{(\lambda)}_{n,k}(\omega)\cdot\exp\{\phi_{\lambda}(\omega_{n-2})-\epsilon\cdot e^{-\varphi_{\lambda}(\omega_{n-2})}\}\notag\\
&&\leq~\delta\cdot\sum_{\omega_{1}=0}^{\infty}\cdots\sum_{\omega_{n-2}=0}^{\infty}\prod_{k=1}^{n-2}Z^{(\lambda)}_{n,k}(\omega)~\leq\cdots\leq~\delta^{\left\lfloor n/2\right\rfloor} \ .
\eea
Hence, $H_{\lambda}\left(\Pna||\Pnh\right) ~\stackrel{n\rightarrow\infty}{\longrightarrow}~0$. \qed \\[0.0cm]

\noindent Notice that the above proof method of formula \eqref{fo.hellnull3d} does not work 
for the parameter setup in Subsection \ref{secDETUP}(a6), because there one has $\delta=\sup_{x\in\mathbb{N}_{0}}\exp\{\phi_{\lambda}(x)-\epsilon\cdot e^{-\varphi_{\lambda}(x)}\}=1$.

\subsection{Proofs of Section \ref{secCFB}}\label{App4} \ \\
\prl\ref{lemua}\ 
Recall the fundamental nonlinear recursion of $\big(a^{(q_{\lambda}^{\bigstar})}_{n}\big)_{n\in\mathbb{N}_{0}}$ (cf.\ 
\eqref{defanvar}, \eqref{fo.zeroSTAR}), the corresponding ``substitute'' inhomogeneous linear recursion
of 
$\big(\ua_{n}^{(q_{\lambda}^{\bigstar})}\big)_{n\in\mathbb{N}_{0}}$ 
(cf.\ \eqref{defua}, \eqref{fo.xiT}, \eqref{defurho})
and its homogenous linear relative 
 $\big(\ua_{n}^{(q_{\lambda}^{\bigstar}),hom}\big)_{n\in\mathbb{N}_{0}}$
 (cf.\ \eqref{fo.defanCLO}, \eqref{fo.xiCLO}) 
 which by \eqref{fo.repanCLO} and \eqref{fo.xiT} takes the form
\bea
&& \ua_{0}^{(q_{\lambda}^{\bigstar}),hom}
:= 0, \qquad
 \ua_{n}^{(q_{\lambda}^{\bigstar}),hom}
:= 
\xi_{\lambda}^{(q_{\lambda}^{\bigstar}),T}\hspace{-0.12cm}\left(\ua^{(q_{\lambda}^{\bigstar}),hom}_{n-1}\right) 
\ = \ x_{0}^{(q_{\lambda}^{\bigstar})}\left(1-\left(d^{(q_{\lambda}^{\bigstar}),T}\right)^{n}\right)
, \quad n\in\mathbb{N}, \label{defanT}
\vspace{-0.2cm}
\eea
with $d^{(q_{\lambda}^{\bigstar}),T}=q_{\lambda}^{\bigstar}\cdot e^{x_{0}^{(q_{\lambda}^{\bigstar})}} \in ]0,1[$. By construction, one has   
\be
\ua_{n}^{(q_{\lambda}^{\bigstar}),hom} \ < \ a^{(q_{\lambda}^{\bigstar})}_{n}  \quad \textrm{for all } n \in \mathbb{N},
\label{fo.anTANinequ}
\qquad \text{as well as} \qquad
\lim_{n\rightarrow\infty}\ua_{n}^{(q_{\lambda}^{\bigstar}),hom}=\lim_{n\rightarrow\infty}a_{n}^{(q_{\lambda}^{\bigstar})}=x_{0}^{(q_{\lambda}^{\bigstar})} \ .
\ee
As an auxiliary step, let us compare $x \mapsto \xi^{(q_{\lambda}^{\bigstar})}_{\lambda}(x)~=~ q_{\lambda}^{\bigstar}\cdot e^x - \beta_{\lambda}$
with the quadratic function 
\bea 
\uu_{\lambda}^{(q_{\lambda}^{\bigstar})}(x) 
& := & \frac{q_{\lambda}^{\bigstar}}{2} \, e^{x_{0}^{(q_{\lambda}^{\bigstar})}}  x^{2}  +  q_{\lambda}^{\bigstar}e^{x_{0}^{(q_{\lambda}^{\bigstar})}}\left(1-x_{0}^{(q_{\lambda}^{\bigstar})}\right)\cdot x  +  x_{0}^{(q_{\lambda}^{\bigstar})}\left(1-q_{\lambda}^{\bigstar}e^{x_{0}^{(q_{\lambda}^{\bigstar})}}+\frac{q_{\lambda}^{\bigstar}}{2}e^{x_{0}^{(q_{\lambda}^{\bigstar})}}x_{0}^{(q_{\lambda}^{\bigstar})}\right). \notag\\[-0.6cm]
&& \notag 
\eea
Clearly, we have the relations $\uu_{\lambda}^{(q_{\lambda}^{\bigstar})}(
x_{0}^{(q_{\lambda}^{\bigstar})}
) 
= x_{0}^{(q_{\lambda}^{\bigstar})}
= \xi^{(q_{\lambda}^{\bigstar})}_{\lambda}(x_{0}^{(q_{\lambda}^{\bigstar})})$,  \ 
$\frac{\partial \uu_{\lambda}^{(q_{\lambda}^{\bigstar})}}{\partial x}(x_{0}^{(q_{\lambda}^{\bigstar})})= q_{\lambda}^{\bigstar}\cdot 
e^{x_{0}^{(q_{\lambda}^{\bigstar})}}
= \frac{\partial \xi^{(q_{\lambda}^{\bigstar})}_{\lambda}}{\partial x}(x_{0}^{(q_{\lambda}^{\bigstar})})$,
and 
$\frac{\partial^2 \uu_{\lambda}^{(q_{\lambda}^{\bigstar})}}{\partial x^2}(x) < \frac{\partial^2 \xi^{(q_{\lambda}^{\bigstar})}_{\lambda}}{\partial x^2}(x)$
for all $x \in ]x_{0}^{(q_{\lambda}^{\bigstar})},0]$. Hence, $\uu_{\lambda}^{(q_{\lambda}^{\bigstar})}(\cdot)$
is on $]x_{0}^{(q_{\lambda}^{\bigstar})},0]$ a strict lower functional bound of $\xi^{(q_{\lambda}^{\bigstar})}_{\lambda}(\cdot)$.
We are now ready to prove  
part (a) by induction. For $n=1$, we easily see that
$\ua_{1}^{(q_{\lambda}^{\bigstar})} <  a_{1}^{(q_{\lambda}^{\bigstar})}$ iff  \, 
$e^{x_{0}^{(q_{\lambda}^{\bigstar})}} \cdot \big\{ \frac{(x_{0}^{(q_{\lambda}^{\bigstar})})^2}{2} - x_{0}^{(q_{\lambda}^{\bigstar})} + 1  \big\} -1 < 0$,
and the latter is obviously true. To continue, let us assume that 
$\ua_{n}^{(q_{\lambda}^{\bigstar})}\leq a_{n}^{(q_{\lambda}^{\bigstar})}$ 
holds. From this, \eqref{defurho}, \eqref{defanT} and \eqref{fo.anTANinequ} we obtain\\[-0.4cm]
\bea  
&& 0 \ < \ \urho_{n}^{(q_{\lambda}^{\bigstar})} \ = \ 
 \frac{q_{\lambda}^{\bigstar}}{2} \, e^{x_{0}^{(q_{\lambda}^{\bigstar})}}
 \left(
 x_{0}^{(q_{\lambda}^{\bigstar})} \cdot
 \left(q_{\lambda}^{\bigstar}\cdot e^{x_{0}^{(q_{\lambda}^{\bigstar})}} \right)^{n} \, 
 \right)^{2} 
 \ = \ 
 \frac{q_{\lambda}^{\bigstar}}{2} \ e^{x_{0}^{(q_{\lambda}^{\bigstar})}}\left(\ua_{n}^{(q_{\lambda}^{\bigstar}),hom}-x_{0}^{(q_{\lambda}^{\bigstar})}\right)^{2}
 \notag \\
 && < \ \frac{q_{\lambda}^{\bigstar}}{2} \ e^{x_{0}^{(q_{\lambda}^{\bigstar})}}\left(a_{n}^{(q_{\lambda}^{\bigstar})}-x_{0}^{(q_{\lambda}^{\bigstar})}\right)^{2}
 \ = \ \uu_{\lambda}^{(q_{\lambda}^{\bigstar})}\left(a_{n}^{(q_{\lambda}^{\bigstar})}\right) - d^{(q_{\lambda}^{\bigstar}),T} \cdot a_{n}^{(q_{\lambda}^{\bigstar})}
 - x_{0}^{(q_{\lambda}^{\bigstar})} \cdot \left(1-d^{(q_{\lambda}^{\bigstar}),T}\right)
 \notag \\
  && < \ \xi^{(q_{\lambda}^{\bigstar})}_{\lambda}\left(a_{n}^{(q_{\lambda}^{\bigstar})}\right) - d^{(q_{\lambda}^{\bigstar}),T} \cdot a_{n}^{(q_{\lambda}^{\bigstar})}
 - x_{0}^{(q_{\lambda}^{\bigstar})} \cdot \left(1-d^{(q_{\lambda}^{\bigstar}),T}\right) \notag\\
 && < \ a_{n+1}^{(q_{\lambda}^{\bigstar})} - d^{(q_{\lambda}^{\bigstar}),T} \cdot \ua_{n}^{(q_{\lambda}^{\bigstar})}
 - x_{0}^{(q_{\lambda}^{\bigstar})} \cdot \left(1-d^{(q_{\lambda}^{\bigstar}),T}\right) \ .
 \notag
\eea
\vspace{-0.5cm}

Thus, $\ua_{n+1}^{(q_{\lambda}^{\bigstar})}\leq a_{n+1}^{(q_{\lambda}^{\bigstar})}$ holds. 
\noindent
In order to show (b), we make use of the straightforward representation
\be 
\ua_{n}^{(q_{\lambda}^{\bigstar})} \ = \ \sum_{k=0}^{n-1}\left(d^{(q_{\lambda}^{\bigstar}),T}\right)^{n-1-k}\cdot
\left( \urho_{k}^{(q_{\lambda}^{\bigstar})} + x_{0}^{(q_{\lambda}^{\bigstar})}\cdot (1- d^{(q_{\lambda}^{\bigstar}),T}) \right) \ 
\notag
\ee
which implies that the sequence $\left(\ua_{n}^{(q_{\lambda}^{\bigstar})}\right)_{n\in\mathbb{N}}$ is strictly decreasing 
since for all $k\in\mathbb{N}_{0}$ there holds by \eqref{defurho}
\vspace{-0.2cm}
\be
\urho_{k}^{(q_{\lambda}^{\bigstar})} + x_{0}^{(q_{\lambda}^{\bigstar})}\cdot (1- d^{(q_{\lambda}^{\bigstar}),T}) 
\ \leq \ \uu_{\lambda}^{(q_{\lambda}^{\bigstar})}(0) \ < \ \xi^{(q_{\lambda}^{\bigstar})}_{\lambda}(0) 
\ = \ q_{\lambda}^{\bigstar} - \beta_{\lambda} \ < \ 0 \ .
\notag
\ee

\noindent
The final assertion  follows immediately from \eqref{fo.anTANinequ} and the closed-form representation \eqref{fo.repanCLO}
with 
the choices $K_{1}$ $K_{2}$, $\varkappa$, $\nu$, $c$  given just right after \eqref{fo.dTundGamma}.
\qed\\[0.05cm]


\noindent
\prl\ref{lemoa}\ \ For $\quasetSP\backslash(\quasetSPdreid\cup\quasetSPvier)$ we deal with the fundamental nonlinear recursion 
of $\big(a^{(q_{\lambda}^{G})}_{n}\big)_{n\in\mathbb{N}_{0}}, \ G\in\{E,U\}$ (cf.\ 
\eqref{defanUP}, \eqref{fo.zeroUPa1}), the corresponding ``substitute'' inhomogeneous linear recursion
of 
$\big(\oa_{n}^{(q_{\lambda}^{G})}\big)_{n\in\mathbb{N}_{0}}$ 
(cf.\ \eqref{defoa}, \eqref{fo.xiS}, \eqref{deforho})
and its homogenous linear counterpart 
 $\big(\oa_{n}^{(q_{\lambda}^{G}),hom}\big)_{n\in\mathbb{N}_{0}}$
 (cf.\ \eqref{fo.defanCLO}, \eqref{fo.xiCLO}) 
 which by \eqref{fo.repanCLO} and \eqref{fo.xiS} takes the form
\bea
 && \oa_{0}^{(q_{\lambda}^{G}),hom}
:= 0, \qquad
 \oa_{n}^{(q_{\lambda}^{G}),hom}
\ := \
\xi_{\lambda}^{(q_{\lambda}^{G}),S}\hspace{-0.12cm}\left(\oa^{(q_{\lambda}^{G}),hom}_{n-1}\right) 
\ = \ x_{0}^{(q_{\lambda}^{G})}\left(1-\left(d^{(q_{\lambda}^{G}),S}\right)^{n}\right), \quad n\in\mathbb{N}, 
\label{defanS}
\eea 
with 
$d^{(q_{\lambda}^{G}),S}=1-\frac{q_{\lambda}^{G}-\beta_{\lambda}}{x_{0}^{(q_{\lambda}^{G})}} \in \big]d^{(q_{\lambda}^{G}),T},1\big[$. 
By construction, we obatin
\be
\oa_{1}^{(q_{\lambda}^{G}),hom}  =  a^{(q_{\lambda}^{G})}_{1},  \ \oa_{n}^{(q_{\lambda}^{G}),hom}  >  a^{(q_{\lambda}^{G})}_{n}  \ \textrm{for all } n \in \mathbb{N}\backslash\{1\},
\ \text{and} \
\lim_{n\rightarrow\infty}\oa_{n}^{(q_{\lambda}^{G}),hom}=\lim_{n\rightarrow\infty}a_{n}^{(q_{\lambda}^{G})}=x_{0}^{(q_{\lambda}^{G})} \ .
\label{fo.anSEKinequ}
\ee

\noindent
In analogy to the Proof of Lemma \ref{lemua}, we use the quadratic function 
\bea 
&& \ou_{\lambda}^{(q_{\lambda}^{G})}(x) 
\ := \ \frac{q_{\lambda}^{G}}{2} \, e^{x_{0}^{(q_{\lambda}^{G})}}\cdot x^{2} \ + \ \left(1-\frac{q_{\lambda}^{G}}{2} \, e^{x_{0}^{(q_{\lambda}^{G})}}x_{0}^{(q_{\lambda}^{G})}-\frac{q_{\lambda}^{G}-\beta_{\lambda}}{x_{0}^{(q_{\lambda}^{G})}}\right)\cdot x \ + \ q_{\lambda}^{G}-\beta_{\lambda}
\notag\\[-0.5cm]
&&\notag
\eea
which satisfies $\ou_{\lambda}^{(q_{\lambda}^{G})}(x_{0}^{(q_{\lambda}^{G})}) = x_{0}^{(q_{\lambda}^{G})} = \xi^{(q_{\lambda}^{G})}_{\lambda}(x_{0}^{(q_{\lambda}^{G})})$,  \ 
$\ou_{\lambda}^{(q_{\lambda}^{G})}(0)=q_{\lambda}^{G}-\beta_{\lambda} = \xi^{(q_{\lambda}^{G})}_{\lambda}(0)$, 
and 
$\frac{\partial^2 \ou_{\lambda}^{(q_{\lambda}^{G})}}{\partial x^2}(x) < \frac{\partial^2 \xi^{(q_{\lambda}^{G})}_{\lambda}}{\partial x^2}(x)$
for all $x \in ]x_{0}^{(q_{\lambda}^{G})},0]$. Hence, $\ou_{\lambda}^{(q_{\lambda}^{G})}(\cdot)$
is on $]x_{0}^{(q_{\lambda}^{G})},0]$ a strict upper functional bound of $\xi^{(q_{\lambda}^{G})}_{\lambda}(\cdot)$.
To start with the proof of part (a), 
let us first observe for $n=1$ the obvious relation
$\oa_{1}^{(q_{\lambda}^{G})}
=q_{\lambda}^{G}-\beta_{\lambda}=a_{1}^{(q_{\lambda}^{G})}=0$.
Furthermore, let us assume that 
$\oa_{n}^{(q_{\lambda}^{G})}\geq a_{n}^{(q_{\lambda}^{G})}$ ($n \in \mathbb{N}$)
holds. From this, 
\eqref{deforho}, \eqref{defanS}, \eqref{fo.anSEKinequ} and the appropriately adapted
version of $\ua_{n}^{(\cdot),hom}$ we obtain the desired inequality $\oa_{n+1}^{(q_{\lambda}^{G})} > a_{n+1}^{(q_{\lambda}^{G})}$
by estimating
\bea  
&& \hspace{-0.7cm} 0  >  \orho_{n}^{(q_{\lambda}^{G})}  =  
-\frac{\left(x_{0}^{(q_{\lambda}^{G})}\right)^{2}}{2}
\cdot\left(q_{\lambda}^{G} \cdot e^{x_{0}^{(q_{\lambda}^{G})}} \right)^{n+1}\hspace{-0.2cm}\cdot
\frac{\oa_{n}^{(q_{\lambda}^{G}),hom}}{x_{0}^{(q_{\lambda}^{G})}} 
 \ = \ 
 \frac{q_{\lambda}^{G}}{2} \ e^{x_{0}^{(q_{\lambda}^{G})}}\left(\ua_{n}^{(q_{\lambda}^{G}),hom}-x_{0}^{(q_{\lambda}^{G})}\right)
 \cdot \oa_{n}^{(q_{\lambda}^{G}),hom}
 \notag \\
 && \hspace{-0.4cm}\geq \frac{q_{\lambda}^{G}}{2} \ e^{x_{0}^{(q_{\lambda}^{G})}}\left(a_{n}^{(q_{\lambda}^{G})}-x_{0}^{(q_{\lambda}^{G})}\right)
 \cdot a_{n}^{(q_{\lambda}^{G})}
\ = \ \ou_{\lambda}^{(q_{\lambda}^{G})}\left(a_{n}^{(q_{\lambda}^{G})}\right) - d^{(q_{\lambda}^{G}),S} \cdot a_{n}^{(q_{\lambda}^{G})}
 - (q_{\lambda}^{G}-\beta_{\lambda})
 \notag \\
  && \hspace{-0.4cm}> \ \xi^{(q_{\lambda}^{G})}_{\lambda}\left(a_{n}^{(q_{\lambda}^{G})}\right) - d^{(q_{\lambda}^{G}),S} \cdot a_{n}^{(q_{\lambda}^{G})}
 - (q_{\lambda}^{G}-\beta_{\lambda}) \ \geq \ a_{n+1}^{(q_{\lambda}^{G})} - d^{(q_{\lambda}^{G}),S} \cdot \oa_{n}^{(q_{\lambda}^{G})}
 - (q_{\lambda}^{G}-\beta_{\lambda}) \ .
\notag
\eea
Moreover, the property (b) follows from the representation
\be 
\oa_{n}^{(q_{\lambda}^{G})} \ = \ \sum_{k=0}^{n-1}\left(d^{(q_{\lambda}^{G}),S}\right)^{n-1-k}\cdot
\left( \orho_{k}^{(q_{\lambda}^{G})} + (q_{\lambda}^{G}-\beta_{\lambda}) \right) \ 
\notag
\ee
which implies that the sequence $\left(\oa_{n}^{(q_{\lambda}^{G})}\right)_{n\in\mathbb{N}}$ is strictly decreasing 
since for all $k\in\mathbb{N}_{0}$ one has $\orho_{k}^{(q_{\lambda}^{G})} + (q_{\lambda}^{G}-\beta_{\lambda}) < 0$.
Finally, 
part (c) follows immediately from \eqref{fo.anSEKinequ} and the closed-form representation \eqref{fo.repanCLO}
with the choices $K_{1}$ $K_{2}$, $\varkappa$, $\nu$, $c$  given just right after \eqref{fo.dS}.
\qed

\subsection{Proofs of Section \ref{sec.diflim}}\label{App5} \ \\
\prt \ref{thm4}
As already mentioned above, one can adapt the proof of Theorem 9.1.3 in Ethier-Kurtz~\cite{EK86} who deal 
with drift-parameters $\eta=0$, $\kappa_{\bullet}=0$, and the different setup of \textit{$\sigma-$independent time-scale} and
a sequence of \textit{critical} Galton-Watson processes \textit{without immigration} with
\textit{general} offspring distribution. For the sake of brevity, we basically outline here only the main
differences to their proof;
for similar limit investigations involving offspring/immigration distributions and parametrizations
which are incompatble to ours, see e.g.\ Sriram \cite{Sri94}.  

\vspace{0.2cm}
\noindent As a first step, let us define 
the generator
\be
A_{\bullet}f(x)~:=~\big(\eta-\kappa_{\bullet}\cdot x\big)~f'(x)+\frac{\sigma^{2}}{2}\cdot x\cdot f''(x), \quad f\in C^{\infty}_{c}\big([0,\infty)\big) \ ,
\notag
\ee
which corresponds to the diffusion process $\widetilde{X}$ governed by \eqref{fo.diffusion}.
In connection with \eqref{def1}, we study 
\be
T^{(m)}_{\bullet}f(x)~:=~EP_{\bullet}\left[f\left(\frac{1}{m}\left(\sum_{k=1}^{mx}
Y^{(m)}_{0,k}+\widetilde{Y}^{(m)}_{0} 
\right)\right)\right],
\quad x \in E^{(m)}:=\frac{1}{m}\mathbb{N}_{0}, \quad f\in C^{\infty}_{c}\big([0,\infty), 
\notag
\ee
where the $Y^{(m)}_{0,k}$, $\widetilde{Y}^{(m)}_{0}$ are independent and (Poisson-$\beta^{(m)}_{\bullet}$  respectively 
Poisson-$\alpha^{(m)}_{\bullet}$) distributed as the members of the
collection $Y^{(m)}$ respectively $\widetilde{Y}^{(m)}$.
By the Theorems 8.2.1 and 1.6.5 as well as Corollary 4.8.9 of \cite{EK86} it is sufficient to show
\be
\label{oper}
\lim_{m\rightarrow\infty}\sup_{x\in E^{(m)}}\left|\sigma^{2}m\Big(T^{(m)}_{\bullet}f(x)-f(x)\Big)-A_{\bullet}f(x)\right|~=~0,~~f\in C^{\infty}_{c}\big([0,\infty)\big)
\ .
\ee
But \eqref{oper} follows mainly from the next

\begin{lem}\label{lem1}
Let
\be
S^{(m)}_{n} \ := \ \frac{1}{\sqrt{n}}\left(\sum_{k=1}^{n}\left(Y^{(m)}_{0,k}-\beta^{(m)}_{\bullet}\right)+\widetilde{Y}^{(m)}_{0}-\alpha^{(m)}_{\bullet}\right) \ , \quad n \in \mathbb{N}, \ m \in \largeint,  
\notag
\ee
with the usual convention $S^{(m)}_{0}:=0$. Then for all $m \in \largeint$, $x\in E^{(m)}$ and all $f\in C^{\infty}_{c}\big([0,\infty)\big)$
\bea\label{eqn2}
&& \epsilon^{(m)}(x) \ := \ EP_{\bullet}\left[\int_{0}^{1}\left(S^{(m)}_{mx}\right)^{2}x(1-v)\left(f''\left(\beta^{(m)}_{\bullet}x+\frac{\alpha^{(m)}_{\bullet}}{m}+v\sqrt{\frac{x}{m}}S^{(m)}_{mx}\right)-f''(x)\right)dv\right]\notag\\
&& = \ \frac{1}{\sigma^{2}}\cdot\left[\sigma^{2}m\cdot\left( T^{(m)}_{\bullet}f(x)-f(x)\right)-A_{\bullet}f(x)\right] \ + \ R^{(m)}, \qquad \textrm{where} \ \lim_{m\rightarrow\infty}R^{(m)}=0.
\eea
\end{lem}
\prl \ref{lem1} \ Let us fix $f\in C^{\infty}_{c}\big([0,\infty)\big)$. 
From the involved Poissonian expectations it is easy to see that
\begin{equation*}
\label{oper1}
\lim_{m\rightarrow\infty} \left|\sigma^{2}m\Big(T^{(m)}_{\bullet}f(0)-f(0)\Big)-A_{\bullet}f(0)\right|~=~0
\ ,
\end{equation*}
and thus \eqref{eqn2} holds for $x=0$. Accordingly, we next consider the case $x\in E^{(m)}\backslash\{0\}$,
with fixed $m\in \largeint$. 
From \, 
$EP_{\bullet}\left[\left(S^{(m)}_{mx}\right)^{2}\right] =
\beta^{(m)}_{\bullet}+\frac{\alpha^{(m)}_{\bullet}}{mx} $ \ 
we obtain
\be\label{eqn3}
EP_{\bullet}\left[\left(S^{(m)}_{mx}\right)^{2}xf''(x)\int_{0}^{1}(1-v)dv\right] \ = \ \frac{1}{2}\left(\beta^{(m)}_{\bullet}\cdot x+\frac{\alpha^{(m)}_{\bullet}}{m}\right)
f''(x) \ =: \ a_{mx} \, \frac{f''(x)}{2} \ =: \ a \, \frac{f''(x)}{2} \ .
\ee
Furthermore, with $ \ b_{mx}  :=\ b  :=  
a +\sqrt{x/m}\cdot S^{(m)}_{mx}  =  
\frac{1}{m}\left(\sum_{k=1}^{mx}Y^{(m)}_{0,k} +\widetilde{Y}^{(m)}_{0}\right) $ 
we get on $\{ S^{(m)}_{mx} \ne 0  \}$ 
\be\label{eqn4}
\int_{0}^{1}f''\left(\beta^{(m)}_{\bullet}x+\frac{\alpha^{(m)}_{\bullet}}{m}+v\sqrt{\frac{x}{m}} \ S^{(m)}_{mx}\right)dv \ = \ \sqrt{\frac{m}{x}}\cdot\frac{1}{S^{(m)}_{mx}}\int_{a}^{b}f''(y)dy \ = \ \sqrt{\frac{m}{x}}\cdot\frac{f'(b)-f'(a)}{S^{(m)}_{mx}}
\ee
as well as
\begin{eqnarray}\label{eqn5}
&& \int_{0}^{1}vf''\left(\beta^{(m)}_{\bullet}x+\frac{\alpha^{(m)}_{\bullet}}{m}+v\sqrt{\frac{x}{m}} \ S^{(m)}_{mx}\right)dv 
\ = \ \frac{m}{x\left(S^{(m)}_{mx}\right)^{2}}\Bigg[\int_{a}^{b}yf''(y) \, dy-
a \int_{a}^{b}f''(y) \, dy\Bigg]\notag\\
&& = \ \sqrt{\frac{m}{x}}\cdot\frac{f'(b)}{S^{(m)}_{mx}}
\ + \ \frac{m}{x} \cdot \frac{f(a)-f(b)}{\left(S^{(m)}_{mx}\right)^{2}} \ .
\end{eqnarray}
With our 
choice $\beta^{(m)}_{\bullet}=1-\frac{\kappa_{\bullet}}{\sigma^{2}m}$ and $\alpha^{(m)}_{\bullet}=\beta^{(m)}_{\bullet}\cdot\frac{\eta}{\sigma^{2}}$, a Taylor expansion of $f$ at $x$ gives
\bea\label{eqn6}
f(a) 
\ = \ f(x)\ + \ \frac{1}{\sigma^{2}m}\cdot f'(x)\left(\beta^{(m)}_{\bullet}\cdot\eta-\kappa_{\bullet}\cdot x\right) \ + \ o\left(\frac{1}{m}\right),
\eea
where for the 
case $\eta=\kappa=0$ we use the convention $o\left(\frac{1}{m}\right) \equiv 0$.
Combining \eqref{eqn3} to \eqref{eqn6} and the centering $EP_{\bullet}\left[S^{(m)}_{mx}\right]=0$, the left hand side of equation~(\ref{eqn2}) becomes
\begin{eqnarray}
&& EP_{\bullet}\left[\int_{0}^{1}\left(S^{(m)}_{mx}\right)^{2}x(1-v)\left(f''\left(\beta^{(m)}_{\bullet}x+\frac{\alpha^{(m)}_{\bullet}}{m}+v\sqrt{\frac{x}{m}}
\ S^{(m)}_{mx}\right)-f''(x)\right)dv\right]\notag\\
& = & EP_{\bullet}\left[\sqrt{mx}\cdot S^{(m)}_{mx}\cdot\Big(f'(b)-f'(a)\Big)\right] \ - \ EP_{\bullet}\left[\sqrt{mx}\cdot S^{(m)}_{mx}\cdot f'(b)
+ m\cdot (f(a)-f(b))
\right]\notag\\
&& - \ \frac{1}{2}\left(\beta^{(m)}_{\bullet}\cdot x+\frac{\alpha^{(m)}_{\bullet}}{m}\right)\cdot f''(x)\notag\\
& = & m\cdot \left( EP_{\bullet}\Big[f(b)\Big]-f(a) \right)  \ - \ \frac{1}{2}\left(\beta^{(m)}_{\bullet}\cdot x+\frac{\alpha^{(m)}_{\bullet}}{m}\right)\cdot f''(x)\notag\\
&= & m\cdot\left\{ EP_{\bullet}\left[f\left(\frac{1}{m}\left(\sum_{k=1}^{mx}Y^{(m)}_{0,k}+\widetilde{Y}_{0}\right)\right)\right]-f(x)\right\}-\frac{1}{\sigma^{2}}A_{\bullet}f(x)\notag\\
&& +  \frac{1}{\sigma^{2}}\left[\left(\eta-\kappa_{\bullet}\cdot x\right)-\beta^{(m)}_{\bullet}\cdot\eta+\kappa_{\bullet}\cdot x\right]\cdot f'(x)  + \frac{x}{2}\left[1-\beta^{(m)}_{\bullet}-\frac{\alpha^{(m)}_{\bullet}}{m}\right]\cdot f''(x)  -  m\cdot o\left(\frac{1}{m}\right)\notag\\
\notag
\end{eqnarray}

\vspace{-0.6cm}
which immediately leads to the right hand side of ~(\ref{eqn2}).
\qed\\

\noindent To proceed with the proof of Theorem \ref{thm4}, we obtain for $m\geq2\kappa_{\bullet}/\sigma^{2}$ the inequality $\beta_{\bullet}^{(m)}\geq1/2$ 
and accordingly for all $v \in ]0,1[$, $x\in E^{(m)}$ 
\be
\label{absch}
\beta^{(m)}_{\bullet}x+\frac{\alpha^{(m)}_{\bullet}}{m}+v\sqrt{\frac{x}{m}}\ S^{(m)}_{mx} 
 =  (1-v)\cdot x\cdot\beta_{\bullet}^{(m)}+(1-v)\frac{\alpha_{\bullet}^{(m)}}{m}+v\left(\sum_{k=1}^{mx}Y_{0,k}^{(m)}+\widetilde{Y}_{0}\right)
\ \geq \ x \cdot \frac{1-v}{2} \, . \notag\\
\ee
Suppose that the support of $f$ is contained in the interval $[0,c]$. Correspondingly, for $v\leq1-2c/x$  
the integrand in $\epsilon^{(m)}(x)$
is zero and hence with \eqref{absch} we can estimate
\bea
&& \left|\int_{0}^{1}\left(S^{(m)}_{mx}\right)^{2}x(1-v)\left(f''\left(\beta^{(m)}_{\bullet}x+\frac{\alpha^{(m)}_{\bullet}}{m}+v\sqrt{\frac{x}{m}} \ S^{(m)}_{mx}\right)-f''(x)\right)dv\right|\notag\\
& \leq & \int_{0\vee(1-2c/x)}^{1}\left(S^{(m)}_{mx}\right)^{2}x(1-v)\cdot2\left\|f''\right\|_{\infty}dv \ \leq \ x\cdot\left(S^{(m)}_{mx}\right)^{2}\left(1\wedge\frac{2c}{x}\right)^{2}\left\|f''\right\|_{\infty}.
\notag
\eea
From this, one can deduce $\lim_{m\rightarrow\infty}\sup_{x\in E^{(m)}}\epsilon^{(m)}(x)=0$ -- and thus \eqref{oper} --
in the same manner as at the end of the proof of Theorem 9.1.3 in~\cite{EK86}
(by means of the dominated convergence theorem).
\qed\\[0.1cm]

\noindent
The following lemma is the main tool for the proof of Theorem \ref{thmlimit} below.
\begin{lem}\label{lem4}
Let $(\kappa_{\Al},\kappa_{\Hy},\eta,\lambda) \in (\widetilde{\mathcal{P}}_{NI}\cup\widetilde{\mathcal{P}}_{SP,1}) \times ]0,1[$.
By using the quantities $\kappa_{\lambda}:=\lambda\kappa_{\Al}+(1-\lambda)\kappa_{\Hy} \, > \, 0$ 
 and $\Lambda_{\!\lambda}:=\sqrt{\lambda\kappa_{\Al}^{2}+(1-\lambda)\kappa_{\Hy}^{2}}  \, > \, \kappa_{\lambda}$ from \eqref{def.kl-Lambda}, 
one gets for all $t>0$ 
\bea
(a) && \hspace{-0.6cm} \lim_{m\rightarrow\infty}m\cdot(1-q_{\lambda}^{(m)})~=~\frac{\kappa_{\lambda}}{\sigma^{2}} \ > \ 0 \, .
	\notag\\
 (b) && \hspace{-0.6cm} \lim_{m\rightarrow\infty}m^{2}\cdot a^{(m)}_{1} = -\frac{\lambda(1-\lambda)\left(\kappa_{\Al}-\kappa_{\Hy}\right)^{2}}{2\sigma^{4}}
	= - \frac{\Lambda_{\!\lambda}^2-\kappa_{\lambda}^2}{2\sigma^{4}}  <  0 \ ; \quad
	\lim_{m\rightarrow\infty}m\cdot(1-\beta_{\lambda}^{(m)}) = \frac{\kappa_{\lambda}}{\sigma^{2}}   >   0   .\notag\\
 (c) && \hspace{-0.6cm} \lim_{m\rightarrow\infty}m\cdot x_{0}^{(m)}~=~-\frac{\Lambda_{\!\lambda}-\kappa_{\lambda}}{\sigma^{2}} \ < \ 0 \ ; \qquad
	\lim_{m\rightarrow\infty}m^2\cdot\Gamma^{(m)}~=~\frac{(\Lambda_{\!\lambda}-\kappa_{\lambda})^2}{2\sigma^4}  \ > \ 0 \, .\notag\\
 (d) && \hspace{-0.6cm} \lim_{m\rightarrow\infty}m\cdot(1-d^{(m),S}) \ = \ \frac{\Lambda_{\!\lambda}+\kappa_{\lambda}}{2\sigma^{2}}  \ > \ 0 \, .	\notag\\
 (e) && \hspace{-0.6cm} \lim_{m\rightarrow\infty}m\cdot(1-d^{(m),T})  =  \frac{\Lambda_{\!\lambda}}{\sigma^{2}}   >  0 \, ;
	\quad
	\lim_{m\rightarrow\infty}m^2\cdot x_{0}^{(m)}\cdot(1-d^{(m),T})~=~-\frac{\Lambda_{\!\lambda}\cdot(\Lambda_{\!\lambda}-\kappa_{\lambda})}{\sigma^4 } 
	\ < \ 0  \, .\notag\\
 (f) && \hspace{-0.6cm} \lim_{m\rightarrow\infty}m\cdot(1-d^{(m),S}d^{(m),T}) \ = \ \frac{3\Lambda_{\!\lambda}+\kappa_{\lambda}}{2\sigma^{2}}  \ > \ 0 \, . 	\notag\\
 (g) && \hspace{-0.6cm} \lim_{m\rightarrow\infty}\left(d^{(m),S}\right)^{\sigma^{2}mt} \ = \ \exp\left\{-\frac{\Lambda_{\!\lambda}+\kappa_{\lambda}}{2}\cdot t\right\}
	 \ < \ 1 \, . 	\notag\\
 (h) && \hspace{-0.6cm} \lim_{m\rightarrow\infty}\left(d^{(m),T}\right)^{\sigma^{2}mt} \ = \ \exp\left\{-\Lambda_{\!\lambda}\cdot t\right\} 	 \ < \ 1 \, . 	\notag\\
 (i) && \hspace{-0.6cm} \lim_{m\rightarrow\infty}\left(d^{(m),S}d^{(m),T}\right)^{\sigma^{2}mt} \ = \ \exp\left\{-\frac{3\Lambda_{\!\lambda}+\kappa_{\lambda}}{2}\cdot t\right\}
		 \ < \ 1 \, . 	\notag\\
 (j) && \hspace{-0.6cm} \lim_{m\rightarrow\infty}m\cdot\underline{\zeta}^{(m)}_{\left\lfloor \sigma^{2}mt\right\rfloor} \ = \ \frac{\left(\Lambda_{\!\lambda}-\kappa_{\lambda}\right)^{2}}{2\sigma^{2}\cdot\Lambda_{\!\lambda}}\cdot e^{-\Lambda_{\!\lambda}\cdot t}\cdot\left(1-e^{-\Lambda_{\!\lambda}\cdot t}\right)
	 \ > \ 0 \, . 	\notag\\
 (k) && \hspace{-0.6cm} \lim_{m\rightarrow\infty}\underline{\vartheta}^{(m)}_{\left\lfloor \sigma^{2}mt\right\rfloor} \ = \ \frac{\eta}{4\sigma^{2}}\cdot\left(\frac{\Lambda_{\!\lambda}-\kappa_{\lambda}}{\Lambda_{\!\lambda}}\right)^{2}\cdot\left(1-e^{-\Lambda_{\!\lambda}\cdot t}\right)^{2}
	 \ \geq \ 0 \, . 	\notag\\
 (l) && \hspace{-0.6cm} \lim_{m\rightarrow\infty}m\cdot\overline{\zeta}^{(m)}_{\left\lfloor \sigma^{2}mt\right\rfloor} \ = \ \frac{\left(\Lambda_{\!\lambda}-\kappa_{\lambda}\right)^{2}}{\sigma^{2}}\cdot\left[\frac{e^{-\frac{1}{2}(\Lambda_{\!\lambda}+\kappa_{\lambda})\cdot t}-e^{-\Lambda_{\!\lambda}\cdot t}}{\Lambda_{\!\lambda}-\kappa_{\lambda}}-\frac{e^{-\frac{1}{2}(\Lambda_{\!\lambda}+\kappa_{\lambda})\cdot t}\left(1-e^{-\Lambda_{\!\lambda}\cdot t}\right)}{2\cdot \Lambda_{\!\lambda}}\right]
		 \ \geq \ 0 \, . 	\notag\\
 (m) && \hspace{-0.6cm} \lim_{m\rightarrow\infty}\overline{\vartheta}^{(m)}_{\left\lfloor \sigma^{2}mt\right\rfloor} \ = \ \frac{\eta}{\sigma^{2}}\frac{\left(\Lambda_{\!\lambda}-\kappa_{\lambda}\right)^{2}}{\Lambda_{\!\lambda}}\cdot\left[\frac{1-e^{-\frac{1}{2}\left(3\Lambda_{\!\lambda}+\kappa_{\lambda}\right)\cdot t}}{3\Lambda_{\!\lambda}+\kappa_{\lambda}}+\frac{e^{-\Lambda_{\!\lambda}\cdot t}-e^{-\frac{1}{2}(\Lambda_{\!\lambda}+\kappa_{\lambda})\cdot t}}{\Lambda_{\!\lambda}-\kappa_{\lambda}}\right]
			 \ \geq \ 0 \, . 	\notag
\eea
\end{lem}
\prl\ref{lem4} \ 
For each of the assertions (a) to (m), we will make use of l'Hospital's rule. To begin with,
we obtain for arbitrary $\mu,\nu\in\mathbb{R}$ 
\bea 
&& \lim_{m\rightarrow\infty}m\cdot\left[1-(\bam)^{\mu}(\bhm)^{\nu}\right]\notag\\
&  = & \lim_{m\rightarrow\infty}m^{2}\cdot\left[\mu \cdot (\bam)^{\mu-1}(\bhm)^{\nu}\frac{\kappa_{\Al}}{\sigma^2\, m^{2}}+\nu \cdot (\bam)^{\mu}(\bhm)^{\nu-1}\frac{\kappa_{\Hy}}{\sigma^2 \, m^{2}}\right] \ = \ \mu~\frac{\kappa_{\Al}}{\sigma^{2}}+\nu~\frac{\kappa_{\Hy}}{\sigma^{2}} \ .
\notag
\eea
From this, (a) follows immediately and (b) can be deduced by
\bea
 \lim_{m\rightarrow\infty}m^{2}\cdot a^{(m)}_{1} & = & \lim_{m\rightarrow\infty}\frac{m}{2\sigma^2}\cdot\Big[\lambda \cdot \kappa_{\Al}\left(1-(\bam)^{\lambda-1}(\bhm)^{1-\lambda}\right)\notag\\
&& + \ (1-\lambda)\cdot \kappa_{\Hy}\left(1-(\bam)^{\lambda}(\bhm)^{-\lambda}\right)\Big] \ = \ -\frac{\lambda(1-\lambda)(\kappa_{\Al}-\kappa_{\Hy})^{2}}{2\sigma^4}
\ .
\notag
\eea
For the proof of the first part of (c), we rely on the inequalities 
$\underline{x}_{0}^{(m)}\leq x_{0}^{(m)}\leq\overline{x}_{0}^{(m)}$ ($m\in\mathbb{N}$), where $\underline{x}_{0}^{(m)}$ and $\overline{x}_{0}^{(m)}$ 
are the obvious notational adaptions of \eqref{defux0} and \eqref{defox0}, respectively. By using (a) and (b), one can calculate
\bea
\lim_{m\rightarrow\infty}m\cdot\overline{x}_{0}^{(m)} & = & \lim_{m\rightarrow\infty}\left(q_{\lambda}^{(m)}\right)^{-1}\cdot\left[m\cdot(1-q_{\lambda}^{(m)})-\sqrt{\left(m\cdot(1-q_{\lambda}^{(m)})\right)^{2}-2\cdot q_{\lambda}^{(m)}\cdot m^{2}\cdot a_{1}^{(m)}}\right]\notag\\
& = & -\frac{\Lambda_{\!\lambda}-\kappa_{\lambda}}{\sigma^{2}} \ .
\notag
\eea
From \eqref{defux0}, (a), (b) and $\lim_{m\rightarrow\infty}\beta_{\lambda}^{(m)}=1$ we obtain   
the limit
\be
\lim_{m\rightarrow\infty}h\left(q_{\lambda}^{(m)}\right) \ = \ \lim_{m\rightarrow\infty}\max\left\{-\beta_{\lambda}^{(m)} \ ; \ \frac{a_{1}^{(m)}}{1-q_{\lambda}^{(m)}}\right\} \ = \ \lim_{m\rightarrow\infty}\frac{1}{m}\cdot\frac{m^{2}\cdot a_{1}^{(m)}}{m\cdot\left(1-q_{\lambda}^{(m)}\right)} \ = \ 0, 
\notag
\ee
which implies\vspace{-0.4cm}
\bea 
\lim_{m\rightarrow\infty}m\cdot\underline{x}_{0}^{(m)} & = & \lim_{m\rightarrow\infty}\frac{e^{-h\left(q_{\lambda}^{(m)}\right)}}{q_{\lambda}^{(m)}} \cdot\Bigg[m\cdot(1-q_{\lambda}^{(m)})\notag\\
&& - \ \sqrt{\left(m\cdot(1-q_{\lambda}^{(m)})\right)^{2}-2 e^{h\left(q_{\lambda}^{(m)}\right)} q_{\lambda}^{(m)}\cdot m^{2}\cdot a_{1}^{(m)}}\Bigg] \ = \ -\frac{\Lambda_{\!\lambda}-\kappa_{\lambda}}{\sigma^{2}}
\notag
\eea
and thus the first part of (c). The second part is an immediate consequence thereof.
Assertion (d) follows from (b) and (c) by 
\vspace{-0.1cm}
\be
\lim_{m\rightarrow\infty}m\cdot(1-d^{(m),S}) \ = \ \lim_{m\rightarrow\infty}
\frac{m^2 \cdot a_{1}^{(m)}}{m \cdot x_{0}^{(m)}}
\ = \ \frac{\Lambda_{\!\lambda}+\kappa_{\lambda}}{2\sigma^{2}}.
\notag
\ee
For the first part of (e), we use the general limit 
$ \ \lim_{x\rightarrow0}\frac{e^{x}-1}{x} \ = \ 1, \notag $ 
to get with (a) and (c)
\be
\lim_{m\rightarrow\infty}m\cdot(1-d^{(m),T}) \ = \ \lim_{m\rightarrow\infty}\left(m\cdot\left(1-q_{\lambda}^{(m)}\right)-q_{\lambda}^{(m)}\cdot m\cdot x_{0}^{(m)}\cdot\frac{e^{x_{0}^{(m)}}-1}{x_{0}^{(m)}}\right) \ = \ \frac{\Lambda_{\!\lambda}}{\sigma^{2}} \ .
\notag
\ee
From this and (c), the second part of (e) is obvious.
The limit (f) can be obtained  from (d) and (e). The assertions (g) respectively (h) respectively (i) follow from (d) respectively (e) respectively (f) 
by using the general relation $\lim_{m\rightarrow\infty}\left(1+\frac{x_m}{m}\right)^{m}=e^{\lim_{m\rightarrow \infty} x_m}$.
The last four limits (j) to (m) are straightforward implications of (a) to (i).
\qed\\

\noindent
\prt\ref{thmlimit} \ 
It suffices to compute the limits of the bounds given in Corollary \ref{cor1} as $m$ tends to infinity. 
This is done by applying Lemma \ref{lem4} 
which provides corresponding limits of various involved quantities.
Accordingly, for all $t>0$ the lower bound \eqref{fo.boulim1} can be obtained from \eqref{fo.bouapprox1} by
\bea
&& \lim_{m\rightarrow\infty}\exp\Bigg\{x_{0}^{(m)}\cdot\left[X_{0}^{(m)}-\frac{\eta}{\sigma^{2}}\cdot\frac{d^{(m),T}}{1-d^{(m),T}}\right]\left(1-\left(d^{(m),T}\right)^{\left\lfloor \sigma^{2}mt\right\rfloor}\right)\notag\\
&& \hspace{2.0cm}+ \ x_{0}^{(m)}\frac{\eta}{\sigma^{2}}\cdot\left\lfloor \sigma^{2}mt\right\rfloor \ + \ \underline{\zeta}^{(m)}_{\left\lfloor \sigma^{2}mt\right\rfloor} \cdot X_{0}^{(m)} \ + \ \underline{\vartheta}^{(m)}_{\left\lfloor \sigma^{2}mt\right\rfloor}\Bigg\}\notag\\
& = & \lim_{m\rightarrow\infty}\exp\Bigg\{m\cdot x_{0}^{(m)}\cdot\left[\frac{X_{0}^{(m)}}{m}-\frac{\eta}{\sigma^{2}}\cdot\frac{d^{(m),T}}{m\cdot\left(1-d^{(m),T}\right)}\right]\left(1-\left(d^{(m),T}\right)^{\left\lfloor \sigma^{2}mt\right\rfloor}\right)\notag\\
&& \hspace{2.0cm}+ \ m\cdot x_{0}^{(m)}\frac{\eta}{\sigma^{2}}\cdot\frac{\left\lfloor \sigma^{2}mt\right\rfloor}{m} \ + \ m\cdot\underline{\zeta}^{(m)}_{\left\lfloor \sigma^{2}mt\right\rfloor} \cdot \frac{X_{0}^{(m)}}{m} \ + \ \underline{\vartheta}^{(m)}_{\left\lfloor \sigma^{2}mt\right\rfloor}\Bigg\}\notag\\
& = & \exp\Bigg\{-\frac{\Lambda_{\!\lambda}-\kappa_{\lambda}}{\sigma^{2}}\cdot\left[\widetilde{X}_{0}-\frac{\eta}{\sigma^{2}}\cdot\frac{\sigma^{2}}{\Lambda_{\!\lambda}}\right]\left(1-e^{-\Lambda_{\!\lambda} t}\right)-\frac{\Lambda_{\!\lambda}-\kappa_{\lambda}}{\sigma^{2}}\cdot\frac{\eta}{\sigma^{2}}\cdot \sigma^{2}t\notag\\
&& \hspace{1cm}+ \ \frac{\left(\Lambda_{\!\lambda}-\kappa_{\lambda}\right)^{2}}{2\sigma^{2}\cdot\Lambda_{\!\lambda}}\cdot e^{-\Lambda_{\!\lambda}\cdot t}\cdot\left(1-e^{-\Lambda_{\!\lambda}\cdot t}\right)\cdot\widetilde{X}_{0}+\frac{\eta}{4\sigma^{2}}\cdot\left(\frac{\Lambda_{\!\lambda}-\kappa_{\lambda}}{\Lambda_{\!\lambda}}\right)^{2}\cdot\left(1-e^{-\Lambda_{\!\lambda}\cdot t}\right)^{2}\Bigg\}\notag\\
& = & \exp\left\{-\frac{\Lambda_{\!\lambda}-\kappa_{\lambda}}{\sigma^{2}}\left[\widetilde{X}_{0}-\frac{\eta}{\Lambda_{\!\lambda}}\right]\left(1-e^{-\Lambda_{\!\lambda}\cdot t}\right)-\frac{\eta}{\sigma^{2}}\left(\Lambda_{\!\lambda}-\kappa_{\lambda}\right)\cdot t \ + \ L_{\lambda}^{(1)}(t)\cdot\widetilde{X}_{0} \ + \ \frac{\eta}{\sigma^{2}}\cdot L_{\lambda}^{(2)}(t)\right\}.
\notag
\eea
For all $t>0$, the upper bound \eqref{fo.boulim} follows analogously  from \eqref{fo.bouapprox}  by 
\bea
&& \lim_{m\rightarrow\infty}\exp\Bigg\{x_{0}^{(m)}\cdot\left[X_{0}^{(m)}-\frac{\eta}{\sigma^{2}}\cdot\frac{d^{(m),S}}{1-d^{(m),S}}\right]\left(1-\left(d^{(m),S}\right)^{\left\lfloor \sigma^{2}mt\right\rfloor}\right)\notag\\
&& \hspace{2.0cm}+ \ x_{0}^{(m)}\frac{\eta}{\sigma^{2}}\cdot\left\lfloor \sigma^{2}mt\right\rfloor \ - \ \overline{\zeta}^{(m)}_{\left\lfloor \sigma^{2}mt\right\rfloor}\cdot X_{0}^{(m)} \ - \ \overline{\vartheta}^{(m)}_{\left\lfloor \sigma^{2}mt\right\rfloor}\Bigg\}\notag\\
& = & \lim_{m\rightarrow\infty}\exp\Bigg\{m\cdot x_{0}^{(m)}\cdot\left[\frac{X_{0}^{(m)}}{m}-\frac{\eta}{\sigma^{2}}\cdot\frac{d^{(m),S}}{m\cdot\left(1-d^{(m),S}\right)}\right]\left(1-\left(d^{(m),S}\right)^{\left\lfloor \sigma^{2}mt\right\rfloor}\right)\notag\\
&& \hspace{2.0cm}+ \ m\cdot x_{0}^{(m)}\frac{\eta}{\sigma^{2}}\cdot\frac{\left\lfloor \sigma^{2}mt\right\rfloor}{m} \ - \ m\cdot\overline{\zeta}^{(m)}_{\left\lfloor \sigma^{2}mt\right\rfloor}\cdot \frac{X_{0}^{(m)}}{m} \ - \ \overline{\vartheta}^{(m)}_{\left\lfloor \sigma^{2}mt\right\rfloor}\Bigg\}\notag
\eea
\bea
& = & \exp\Bigg\{-\frac{\Lambda_{\!\lambda}-\kappa_{\lambda}}{\sigma^{2}}\left[\widetilde{X}_{0}-\frac{\eta}{\sigma^{2}}\cdot\frac{8\sigma^{2}}{\Lambda_{\!\lambda}+\kappa_{\lambda}}\right]\left(1-\left(e^{-\frac{1}{2}(\Lambda_{\!\lambda}+\kappa_{\lambda})t}\right)\right)-\frac{\Lambda_{\!\lambda}-\kappa_{\lambda}}{\sigma^{2}}\cdot\frac{\eta}{\sigma^{2}}\cdot\sigma^{2}t\notag\\
&& \hspace{1cm} + \ \frac{\left(\Lambda_{\!\lambda}-\kappa_{\lambda}\right)^{2}}{\sigma^{2}}\cdot\left[\frac{e^{-\frac{1}{2}(\Lambda_{\!\lambda}+\kappa_{\lambda})\cdot t}-e^{-\Lambda_{\!\lambda}\cdot t}}{\Lambda_{\!\lambda}-\kappa_{\lambda}}-\frac{e^{-\frac{1}{2}(\Lambda_{\!\lambda}+\kappa_{\lambda})\cdot t}\left(1-e^{-\Lambda_{\!\lambda}\cdot t}\right)}{2\cdot \Lambda_{\!\lambda}}\right]\cdot\widetilde{X}_{0}\notag\\
&& \hspace{1cm} + \ \frac{\eta}{\sigma^{2}}\frac{\left(\Lambda_{\!\lambda}-\kappa_{\lambda}\right)^{2}}{\Lambda_{\!\lambda}}\cdot\left[\frac{1-e^{-\frac{1}{2}\left(3\Lambda_{\!\lambda}+\kappa_{\lambda}\right)\cdot t}}{3\Lambda_{\!\lambda}+\kappa_{\lambda}}+\frac{e^{-\Lambda_{\!\lambda}\cdot t}-e^{-\frac{1}{2}(\Lambda_{\!\lambda}+\kappa_{\lambda})\cdot t}}{\Lambda_{\!\lambda}-\kappa_{\lambda}}\right]\Bigg\}\notag\\
& = & \exp\Bigg\{-\frac{\Lambda_{\!\lambda}-\kappa_{\lambda}}{\sigma^{2}}\left[\widetilde{X}_{0}-\frac{\eta}{\frac{1}{2}(\Lambda_{\!\lambda}+\kappa_{\lambda})}\right]\left(1-e^{-\frac{1}{2}(\Lambda_{\!\lambda}+\kappa_{\lambda})\cdot t}\right)-\frac{\eta}{\sigma^{2}}\left(\Lambda_{\!\lambda}-\kappa_{\lambda}\right)\cdot t\notag\\
&& \hspace{1.0cm}- \ U_{\lambda}^{(1)}(t)\cdot\widetilde{X}_{0} \ - \  \frac{\eta}{\sigma^{2}}\cdot U_{\lambda}^{(2)}(t)\Bigg\}.
\hspace{6.0cm} \square \notag
\eea

\subsection{Proofs of Section \ref{sec.ent}}\label{App6}

We start with two lemmas which will be useful for the proof of Theorem \ref{thm.entex},
and which can be easily seen by induction. They deal with
the sequence $\left(a_{n}^{(q_{\lambda})}\right)_{n\in\mathbb{N}}$ from \eqref{defan}. 

\begin{lem}\label{lem.anqlambda}
For arbitrarily fixed parameter constellation $\qua\in\quaset$, suppose that $q_{\lambda} >0$ $(\lambda \in ]0,1[)$
and $\lim_{\lambda\nearrow1} \ q_{\lambda} \ = \ \bal$ holds. Then one gets the limit
\be\label{fo.anzero}
\forall \ n\in\mathbb{N}: \quad \lim_{\lambda\nearrow1} \ a_{n}^{(q_{\lambda})} \ = 0.
\ee
\end{lem}

\begin{lem}\label{lem.derliman}
In addition to the assumptions of Lemma \ref{lem.anqlambda}, 
suppose that $\lambda \mapsto q_{\lambda}$ is continuously differentiable 
on $]0,1[$
and that the limit $l :=\lim_{\lambda\nearrow1}\frac{\partial\, q_{\lambda}}{\partial\lambda}$
is finite. 
Then one gets the limit
\be  
\forall \ n\in\mathbb{N}: \quad \lim_{\lambda\nearrow1}\frac{\partial\, a_{n}^{(q_{\lambda})}}{\partial\lambda} \ = \ u_{n} \ := \ \left\{
\begin{array}{ll}
\frac{l + \bhy - \bal}{1-\bal}\cdot\left(1-\left(\bal\right)^{n}\right) \, , & \textrm{if }\, \bal\neq1,\\
~&~\\
n\cdot\left(l + \bhy - 1 \right) \, , & \textrm{if }\, \bal=1,
\end{array}
\right.
\notag
\ee 
which is the unique solution of the linear recursion equation 
\be 
u_{n}  \ = \ l + \bhy - \bal \ + \ \bal\cdot u_{n-1}  \ , \qquad
u_{0}  \ = \ 0 \ . 
\notag
\ee
Furthermore,
\be 
\forall \ n\in\mathbb{N}: \quad \sum_{k=1}^{n} \lim_{\lambda\nearrow1}\frac{\partial\, a_{k}^{(q_{\lambda})}}{\partial\lambda} \ = \  \sum_{k=1}^{n} u_{k} \ = \ \left\{
\begin{array}{ll}
\frac{l + \bhy - \bal}{1-\bal}\cdot\left[n-\frac{\bal}{1-\bal}\left(1-\left(\bal\right)^{n}\right)\right] \, , & \textrm{if }\, \bal\neq1,\\
~&~\\
\frac{n\cdot(n+1)}{2}\cdot\left(l + \bhy - 1 \right) \, , & \textrm{if }\, \bal=1.
\end{array}
\right.
\notag
\ee
\end{lem}




\vspace{0.3cm}

\noindent We are now ready to give the

\vspace{0.3cm}

\noindent
\prt\ref{thm.entex} \\
(a) Recall that for the setup $\qua\in(\quasetNI\cup\quasetSPeins)$ we chose 
the intercept as $p_{\lambda} := p_{\lambda}^{E}:=\aal^{\lambda}\ahy^{1-\lambda}$ and the slope as $q_{\lambda}:= q_{\lambda}^{E}:=\bal^{\lambda}\bhy^{1-\lambda}$,
which in \eqref{fo.genequality} lead to the exact value $V_{\lambda,n}$ 
of the Hellinger integral.  
Because of  $\frac{p_{\lambda}}{q_{\lambda}}\beta_{\lambda}-\alpha_{\lambda}=0$ 
as well as $\lim_{\lambda\nearrow1}q_{\lambda}=\bal$, we obtain by using \eqref{fo.anbn} and Lemma \ref{lem.anqlambda} for all $n\in\mathbb{N}$
\be 
\lim_{\lambda\nearrow1}V_{\lambda,n}
\ := \ \lim_{\lambda\nearrow1}\exp\left\{a_{n}^{(q_{\lambda})}\cdot\omega_{0}+\sum_{k=1}^{n}b_{k}^{(p_{\lambda},q_{\lambda})}\right\} \ = \ \lim_{\lambda\nearrow1}\exp\left\{a_{n}^{(q_{\lambda})}\cdot\omega_{0}+\frac{\aal}{\bal}\sum_{k=1}^{n}a_{k}^{(q_{\lambda})}\right\} \ = \ 1, 
 n\in\mathbb{N},
\notag
\ee
which leads 
by \eqref{fo.liment} to \vspace{-0.3cm}
\bea\label{fo.entder}
I(\Pna||\Pnh) & = & \lim_{\lambda\nearrow1} \ \frac{1-H_{\lambda}(\Pna||\Pnh)}{\lambda\cdot(1-\lambda)} \ = \  \lim_{\lambda\nearrow1} \ \frac{1-V_{\lambda,n}
}{\lambda\cdot(1-\lambda)}\notag\\
& = & \lim_{\lambda\nearrow1}\frac{-V_{\lambda,n}
}{1-2\lambda}\cdot\frac{\partial}{\partial\lambda}\left[a_{n}^{(q_{\lambda})}\cdot\omega_{0}+\frac{p_{\lambda}}{q_{\lambda}}\sum_{k=1}^{n}a_{k}^{(q_{\lambda})}\right]\notag\\
& = & \lim_{\lambda\nearrow1}\left[\frac{\partial \, a_{n}^{(q_{\lambda})}}{\partial\lambda}\cdot\omega_{0}+\left(\frac{\partial}{\partial\lambda}\frac{p_{\lambda}}{q_{\lambda}}\right)\cdot\sum_{k=1}^{n}a_{k}^{(q_{\lambda})}+\frac{p_{\lambda}}{q_{\lambda}}\cdot\sum_{k=1}^{n}\frac{\partial \, a_{k}^{(q_{\lambda})}}{\partial\lambda}\right].
\eea
For further analysis, we use the obvious derivatives
\be\label{fo.dera}
\frac{\partial \, p_{\lambda}}{\partial\lambda} \ = \ p_{\lambda} \, \log\left(\frac{\aal}{\ahy}\right),\qquad 
\frac{\partial}{\partial\lambda}\frac{p_{\lambda}}{q_{\lambda}} \ = \ \frac{p_{\lambda}}{q_{\lambda}}  \, \log\left(\frac{\aal\bhy}{\ahy\bal}\right),
\qquad\frac{\partial \, q_{\lambda}}{\partial\lambda} \ = \ q_{\lambda} \, \log\left(\frac{\bal}{\bhy}\right),
\ee
where the subcase $\qua\in\quasetNI$ (with $p_{\lambda}\equiv 0$) 
is consistently covered.
From \eqref{fo.dera}
and Lemma \ref{lem.derliman} 
we deduce
\be 
\lim_{\lambda\nearrow1}\frac{\partial \, a_{n}^{(q_{\lambda})}}{\partial\lambda}\cdot\omega_{0} \ = \ \left\{
\begin{array}{ll}
\left(\bal\log\left(\frac{\bal}{\bhy}\right)-(\bal-\bhy)\right)\cdot\frac{1-\left(\bal\right)^{n}}{1-\bal}\cdot\omega_{0} & \textrm{if } \bal\neq1,\\
n\cdot\left(\bal\log\left(\frac{\bal}{\bhy}\right)-(\bal-\bhy)\right)\cdot\omega_{0} & \textrm{if } \bal=1,
\end{array}
\right.
\notag
\ee
and by means of \eqref{fo.anzero} 
\be 
\forall \ n\in\mathbb{N}: \quad \lim_{\lambda\nearrow1} \left[\left(\frac{\partial}{\partial\lambda}\frac{p_{\lambda}}{q_{\lambda}}\right)\cdot\sum_{k=1}^{n}a_{k}^{(q_{\lambda})} \right] \ = \ 0.
\notag
\ee
For the last expression in \eqref{fo.entder} we again apply Lemma \ref{lem.derliman} 
to end up with
\be\label{fo.entderproofc}
\lim_{\lambda\nearrow1}\frac{p_{\lambda}}{q_{\lambda}}\cdot\sum_{k=1}^{n}\frac{\partial}{\partial\lambda}a_{k}^{(q_{\lambda})} \ = \ \left\{
\begin{array}{ll}
\frac{\aal\cdot\left[\bal\log\left(\frac{\bal}{\bhy}\right)-(\bal-\bhy)\right]}{\bal(1-\bal)}\cdot\left[n-\frac{\bal}{1-\bal}\left(1-\left(\bal\right)^{n}\right)\right] & \textrm{if }\bal\neq1,\\
n\cdot(n+1)\frac{\aal}{2\bal}\cdot\left[\bal\log\left(\frac{\bal}{\bhy}\right)-(\bal-\bhy)\right] & \textrm{if }\bal=1,
\end{array}
\right.
\ee
which finishes the proof of part (a).
To show part (b), for the corresponding setup $\qua$ $\in \quasetSPcompvar$ let us first choose
-- according to the Section \ref{secDETLOW} --
the intercept as $p_{\lambda} := p_{\lambda}^{L}:=\aal^{\lambda}\ahy^{1-\lambda}$ and the slope as $q_{\lambda}:= q_{\lambda}^{L}:=\bal^{\lambda}\bhy^{1-\lambda}$,
which in part (b) of Proposition \ref{propLOW} lead to the lower bounds $B_{\lambda,n}^{L}$ 
of the Hellinger integral. 
This is formally the same choice as in part (a) satisfying $\lim_{\lambda\nearrow1}p_{\lambda}=\aal$, $\lim_{\lambda\nearrow1}q_{\lambda}=\bal$ but 
in contrast to (a) we now have $\frac{p_{\lambda}}{q_{\lambda}} \, \beta_{\lambda}-\alpha_{\lambda}\neq0$ but nevertheless
\be
\lim_{\lambda\nearrow1}\frac{p_{\lambda}}{q_{\lambda}} \, \beta_{\lambda}-\alpha_{\lambda} \ = \ 0.
\notag
\ee
From this, \eqref{fo.anbn}, part (b) of Proposition \ref{propLOW} and Lemma \ref{lem.anqlambda} we obtain 
\be
\label{limBL}
\lim_{\lambda\nearrow1}B_{\lambda,n}^{L} \ = \ \lim_{\lambda\nearrow1}\exp\left\{a_{n}^{(q_{\lambda})}\cdot\omega_{0}+\frac{p_{\lambda}}{q_{\lambda}}\, \sum_{k=1}^{n}a_{k}^{(q_{\lambda})}+n\cdot\left(\frac{p_{\lambda}}{q_{\lambda}}\beta_{\lambda}-\alpha_{\lambda}\right)\right\} \ = \ 1
\ee
and hence 
\bea\label{fo.proofub}
I(\Pna||\Pnh) & \leq & \hspace{-0.1cm}\lim_{\lambda\nearrow1}  \frac{1-B_{\lambda,n}^{L}}{\lambda\cdot(1-\lambda)}  =  \lim_{\lambda\nearrow1}\frac{-B_{\lambda,n}^{L}}{1-2\lambda}\cdot\frac{\partial}{\partial\lambda}\left[a_{n}^{(q_{\lambda})}\omega_{0}+\frac{p_{\lambda}}{q_{\lambda}}\sum_{k=1}^{n}a_{k}^{(q_{\lambda})}+n\left(\frac{p_{\lambda}}{q_{\lambda}}\, \beta_{\lambda}-\alpha_{\lambda}\right)\right]\notag\\
& = & \hspace{-0.1cm}\lim_{\lambda\nearrow1}\left[\frac{\partial \, a_{n}^{(q_{\lambda})}}{\partial\lambda}\omega_{0}+\left(\frac{\partial}{\partial\lambda}\frac{p_{\lambda}}{q_{\lambda}}\right)\sum_{k=1}^{n}a_{k}^{(q_{\lambda})}+\frac{p_{\lambda}}{q_{\lambda}}\sum_{k=1}^{n}\frac{\partial \, a_{k}^{(q_{\lambda})}}{\partial\lambda}+n\frac{\partial}{\partial\lambda}\left(\frac{p_{\lambda}}{q_{\lambda}}\beta_{\lambda}-\alpha_{\lambda}\right)\right].
\eea
In the current setup, the first three expressions in \eqref{fo.proofub} can be evaluated in exactly the same way as 
in \eqref{fo.dera} to \eqref{fo.entderproofc}, and for the last expression one has the limit
\bea 
\frac{\partial}{\partial\lambda}\left(\frac{p_{\lambda}}{q_{\lambda}}\, \beta_{\lambda}-\alpha_{\lambda}\right) & = &  \frac{p_{\lambda}}{q_{\lambda}} \, \log\left(\frac{\aal\bhy}{\ahy\bal}\right)\cdot\beta_{\lambda} \ + \ \frac{p_{\lambda}}{q_{\lambda}}\cdot\left(\bal-\bhy\right) \ - \ \left(\aal-\ahy\right)\notag\\
& \stackrel{\lambda\nearrow1}{\longrightarrow} & \aal\left[\log\left(\frac{\aal\bhy}{\ahy\bal}\right)-\frac{\bhy}{\bal}\right]+\ahy \ ,
\notag\\[-0.6cm]
&&\notag
\eea
which finishes the proof of part (b).
\qed\\

\noindent
\prt \ref{thm.entexUP} \ 
Let us fix $\qua \in  \quasetSP\backslash \quasetSPeins$, $\omega_{0}\in\mathbb{N}$, $n \in \mathbb{N}$ and $y\in[0,\infty[$.
The lower bound $E^{L,tan}_{y,n}$ of the relative entropy is derived by using as a linear upper bound 
$\phi_{\lambda}^{U}$ (cf.~\eqref{fo.phibou}) for $\phi_{\lambda}$  ($\lambda \in ]0,1[$)  
the tangent line of $\phi_{\lambda}$ at $y$. This corresponds to 
$\phi_{\lambda}^{U}(x) := (p_{\lambda}^{U} - \alpha_{\lambda}) + (q_{\lambda}^{U} - \beta_{\lambda}) \, x$ 
($x\in[0,\infty[$) with $p_{\lambda}:= p_{\lambda}(y) := \phi_{\lambda}(y)-y\phi_{\lambda}'(y)+\alpha_{\lambda}$ 
and $q_{\lambda}:=q_{\lambda}(y):= \phi_{\lambda}'(y)+\beta_{\lambda}$, implying $q_{\lambda}>0$ because of 
(p-xii). 
As a side remark, notice that this $\phi_{\lambda}^{U}(x)$ may become negative for some $x\in[0,\infty[$ (which is not always
consistent with goal (Gc) for fixed $\lambda$, but leads to a tractable limit bound as $\lambda$ tends to 1). 
Analogously to \eqref{limBL} and \eqref{fo.proofub}, 
we obtain from \eqref{fo.anbn} and \eqref{fo.genbounds} the convergence 
$\lim_{\lambda\nearrow1}B_{\lambda,n}^{U}=1$ and thus
\be\label{pr.thmentup0}
I(\Pna||\Pnh) \ \geq \ \lim_{\lambda\nearrow1}\left[\frac{\partial \, a_{n}^{(q_{\lambda})}}{\partial\lambda}\omega_{0}+\left(\frac{\partial}{\partial\lambda}\frac{p_{\lambda}}{q_{\lambda}}\right)\sum_{k=1}^{n}a_{k}^{(q_{\lambda})}+\frac{p_{\lambda}}{q_{\lambda}}\sum_{k=1}^{n}\frac{\partial \, a_{k}^{(q_{\lambda})}}{\partial\lambda}+n\frac{\partial}{\partial\lambda}\left(\frac{p_{\lambda}}{q_{\lambda}}\beta_{\lambda}-\alpha_{\lambda}\right)\right].
\ee
As before, we compute the involved derivatives. From \eqref{defphi} to \eqref{defflambda} as well as (p-xii) we get
\bea
&& \hspace{-0.3cm}\frac{\partial p_{\lambda}}{\partial\lambda} \ = \ \left(\frac{f_{\Al}(y)}{f_{\Hy}(y)}\right)^{\lambda}f_{\Hy}(y)\log\left(\frac{f_{\Al}(y)}{f_{\Hy}(y)}\right)-\bal y\left(\frac{f_{\Al}(y)}{f_{\Hy}(y)}\right)^{\lambda-1}\hspace{-0.3cm}-\lambda\bal y\left(\frac{f_{\Al}(y)}{f_{\Hy}(y)}\right)^{\lambda-1}\hspace{-0.1cm}\log\left(\frac{f_{\Al}(y)}{f_{\Hy}(y)}\right)\notag\\
&&  \hspace{1.2cm} + \ \bhy y\left(\frac{f_{\Al}(y)}{f_{\Hy}(y)}\right)^{\lambda}-(1-\lambda)\bhy y\left(\frac{f_{\Al}(y)}{f_{\Hy}(y)}\right)^{\lambda}\log\left(\frac{f_{\Al}(y)}{f_{\Hy}(y)}\right)\notag\\
&& \hspace{0.5cm}\stackrel{\lambda\nearrow1}{\longrightarrow} \ \aal\log\left(\frac{f_{\Al}(y)}{f_{\Hy}(y)}\right)+\frac{y\cdot(\aal\bhy-\ahy\bal)}{f_{\Hy}(y)} \ ,
\label{pr.thmentup1}
\eea
and \vspace{-0.5cm}
\bea
\frac{\partial q_{\lambda}}{\partial\lambda} & = & \bal\left(\frac{f_{\Al}(y)}{f_{\Hy}(y)}\right)^{\lambda-1}+\lambda\bal\left(\frac{f_{\Al}(y)}{f_{\Hy}(y)}\right)^{\lambda-1}\log\left(\frac{f_{\Al}(y)}{f_{\Hy}(y)}\right)-\bhy\left(\frac{f_{\Al}(y)}{f_{\Hy}(y)}\right)^{\lambda}\notag\\
&& + \ (1-\lambda)\bhy\left(\frac{f_{\Al}(y)}{f_{\Hy}(y)}\right)^{\lambda}\log\left(\frac{f_{\Al}(y)}{f_{\Hy}(y)}\right)\notag\\
& \stackrel{\lambda\nearrow1}{\longrightarrow} & \bal\left(1+\log\left(\frac{f_{\Al}(y)}{f_{\Hy}(y)}\right)\right)-\bhy\frac{f_{\Al}(y)}{f_{\Hy}(y)} \quad =: \quad l.
\label{pr.thmentup2}
\eea
Combining these two limits we get
\bea
\frac{\partial}{\partial\lambda}\left(\frac{p_{\lambda}}{q_{\lambda}}\beta_{\lambda}-\alpha_{\lambda}\right) & = & \frac{q_{\lambda}\left(\frac{\partial p_{\lambda}}{\partial\lambda}\right)-p_{\lambda}\left(\frac{\partial q_{\lambda}}{\partial\lambda}\right)}{(q_{\lambda})^{2}}\cdot\beta_{\lambda}+\frac{p_{\lambda}}{q_{\lambda}}\cdot\left(\bal-\bhy\right)-\left(\aal-\ahy\right)\notag\\
& \stackrel{\lambda\nearrow1}{\longrightarrow} & \left[\frac{y\cdot(\aal\bhy-\ahy\bal)}{f_{\Hy}(y)}-\aal\left(1-\frac{\bhy f_{\Al}(y)}{\bal f_{\Hy}(y)}\right)\right]+\ahy-\frac{\aal\bhy}{\bal}.
\label{pr.thmentup3}
\eea
The above calculation also implies that $\lim_{\lambda\nearrow1}\left(\frac{\partial}{\partial\lambda}\frac{p_{\lambda}}{q_{\lambda}}\right)$ is finite 
 and thus  
 $\lim_{\lambda\nearrow1}\left(\frac{\partial}{\partial\lambda}\frac{p_{\lambda}}{q_{\lambda}}\right)\sum_{k=1}^{n}a_{k}^{(q_{\lambda})}=0$
by means of Lemma \ref{lem.anqlambda}. The proof of \, $I(\Pna||\Pnh) \geq E^{L,\textrm{tan}}_{y,n}$ \, 
is finished by using Lemma \ref{lem.derliman} with $l$ defined in \eqref{pr.thmentup2} and by plugging the limits \eqref{pr.thmentup1} to \eqref{pr.thmentup3} into \eqref{pr.thmentup0}.\\[0.1cm]
\indent To derive the lower bound $E^{L,\textrm{sec}}_{k,n}$ (cf.~\eqref{fo.ELkn}) for fixed $k \in \mathbb{N}_{0}$,
we use as a linear upper bound $\phi_{\lambda}^{U}$ for $\phi_{\lambda}(\cdot)$  ($\lambda \in ]0,1[$)
the secant line of $\phi_{\lambda}$ through the points $k$ and $k+1$, corresponding to 
the choices $p_{\lambda} := p_{\lambda}(k):=(k+1)\cdot\phi_{\lambda}(k)-k\cdot\phi_{\lambda}(k+1)+\alpha_{\lambda}$ and $q_{\lambda}:=q_{\lambda}(k):=\phi_{\lambda}(k+1)-\phi_{\lambda}(k)+\beta_{\lambda}$, 
implying $q_{\lambda}>0$ because of 
(p-xiii) and (p-iv). 
As a side remark, notice that this $\phi_{\lambda}^{U}(x)$ may become negative for some $x\in[0,\infty[$ (which is not always
consistent with goal (Gc) for fixed $\lambda$, but leads to a tractable limit bound as $\lambda$ tends to 1). 
Analogously to \eqref{limBL} and \eqref{fo.proofub} we get again $\lim_{\lambda\nearrow1}B_{\lambda,n}^{U}=1$, which leads 
to the lower bound 
given in \eqref{pr.thmentup0} with appropriately plugged-in quantities. 
As in the above proof of the lower bound $E^{L,tan}_{y,n}$, the inequality 
 \, $I(\Pna||\Pnh) \geq E^{L,\textrm{sec}}_{k,n}$ \, 
follows straightforwardly from Lemma \ref{lem.anqlambda}, Lemma \ref{lem.derliman} and the 
three limits
\bea
& &\hspace{-0.4cm} \frac{\partial p_{\lambda}}{\partial\lambda} \ = \ \left(\frac{f_{\Al}(k)}{f_{\Hy}(k)}\right)^{\lambda}f_{\Hy}(k)\cdot(k\hspace{-0.1cm}+\hspace{-0.1cm}1)\log\left(\frac{f_{\Al}(k)}{f_{\Hy}(k)}\right)-\left(\frac{f_{\Al}(k\hspace{-0.08cm}+\hspace{-0.08cm}1)}{f_{\Hy}(k\hspace{-0.08cm}+\hspace{-0.08cm}1)}\right)^{\lambda}f_{\Hy}(k\hspace{-0.08cm}+\hspace{-0.08cm}1)\cdot k\log\left(\frac{f_{\Al}(k\hspace{-0.08cm}+\hspace{-0.08cm}1)}{f_{\Hy}(k\hspace{-0.08cm}+\hspace{-0.08cm}1)}\right)\notag\\
& & \hspace{-0.4cm} \stackrel{\lambda\nearrow1}{\longrightarrow}  f_{\Al}(k)(k\hspace{-0.08cm}+\hspace{-0.08cm}1)\log\left(\frac{f_{\Al}(k)}{f_{\Hy}(k)}\right)-f_{\Al}(k\hspace{-0.08cm}+\hspace{-0.08cm}1)k\log\left(\frac{f_{\Al}(k\hspace{-0.08cm}+\hspace{-0.08cm}1)}{f_{\Hy}(k\hspace{-0.08cm}+\hspace{-0.08cm}1)}\right),
\notag\\
& & \hspace{-0.4cm} \frac{\partial q_{\lambda}}{\partial\lambda} \ = \  \left(\frac{f_{\Al}(k\hspace{-0.08cm}+\hspace{-0.08cm}1)}{f_{\Hy}(k\hspace{-0.08cm}+\hspace{-0.08cm}1)}\right)^{\lambda}f_{\Hy}(k\hspace{-0.08cm}+\hspace{-0.08cm}1)\log\left(\frac{f_{\Al}(k\hspace{-0.08cm}+\hspace{-0.08cm}1)}{f_{\Hy}(k\hspace{-0.08cm}+\hspace{-0.08cm}1)}\right)-\left(\frac{f_{\Al}(k)}{f_{\Hy}(k)}\right)^{\lambda}f_{\Hy}(k)\log\left(\frac{f_{\Al}(k)}{f_{\Hy}(k)}\right)\notag\\
& & \hspace{-0.4cm} \stackrel{\lambda\nearrow1}{\longrightarrow} f_{\Al}(k\hspace{-0.08cm}+\hspace{-0.08cm}1)\log\left(\frac{f_{\Al}(k\hspace{-0.08cm}+\hspace{-0.08cm}1)}{f_{\Hy}(k\hspace{-0.08cm}+\hspace{-0.08cm}1)}\right)-f_{\Al}(k)\log\left(\frac{f_{\Al}(k)}{f_{\Hy}(k)}\right) \quad =: \quad l \ , \qquad \quad \textrm{and}
\notag\\
&& \hspace{-0.4cm} \frac{\partial}{\partial\lambda}\left(\frac{p_{\lambda}}{q_{\lambda}}\beta_{\lambda}-\alpha_{\lambda}\right) 
\ = \ 
\frac{q_{\lambda}\left(\frac{\partial p_{\lambda}}{\partial\lambda}\right)-p_{\lambda}\left(\frac{\partial q_{\lambda}}{\partial\lambda}\right)}{(q_{\lambda})^{2}}\cdot\beta_{\lambda}+\frac{p_{\lambda}}{q_{\lambda}}\cdot\left(\bal-\bhy\right)-\left(\aal-\ahy\right)\notag\\
& & \hspace{-0.2cm} \stackrel{\lambda\nearrow1}{\longrightarrow}  f_{\Al}(k)\log\left(\frac{f_{\Al}(k)}{f_{\Hy}(k)}\right)\left(k\hspace{-0.08cm}+\hspace{-0.08cm}1+\frac{\aal}{\bal}\right)-f_{\Al}(k\hspace{-0.08cm}+\hspace{-0.08cm}1)\log\left(\frac{f_{\Al}(k\hspace{-0.08cm}+\hspace{-0.08cm}1)}{f_{\Hy}(k\hspace{-0.08cm}+\hspace{-0.08cm}1)}\right)\left(k+\frac{\aal}{\bal}\right) -  \frac{\aal\bhy}{\bal}  +  \ahy.
\notag
\eea
To construct the third lower bound 
$E^{L,hor}_{n}$ (cf.~\eqref{fo.ELhorn}), we start by using for each fixed $\lambda \in ]0,1[$ as an upper bound of $\phi_{\lambda}$ the horizontal line through the intercept $\sup_{x\in\mathbb{N}_{0}}\phi_{\lambda}(x)$. For $\quasetSPdreiab \cup \quasetSPdreic$, this supremum is achived
at the finite integer point $z_{\lambda}^{*}:=\arg\max_{x\in\mathbb{N}_{0}}\phi_{\lambda}(x)$ 
(since $\lim_{x\rightarrow\infty}\phi_{\lambda}(x)=-\infty$) and one has $\phi_{\lambda}(z_{\lambda}^{*}) < 0$ which
leads with the setup $q_{\lambda}=\beta_{\lambda}$, $p_{\lambda}=\phi_{\lambda}(z_{\lambda}^{*})+\alpha_{\lambda}$ to the Hellinger integral upper bound $B^{U}_{\lambda,n}=\exp\left\{\phi_{\lambda}(z_{\lambda}^{*})\cdot n\right\}<1$ (cf.\ \eqref{fo.HI.horizontal}). 
To compute from this the required 
$\lim_{\lambda\nearrow1}\frac{1-B^{U}_{\lambda,n}}{\lambda(1-\lambda)}$
is not straightforward since in general it seems to be intractable to express $z_{\lambda}^{*}$ explicitly in terms of $\lambda$.
However, since $\lim_{\lambda\nearrow1}\phi_{\lambda}(x)=0$ for all $x\in [0,\infty[$, we obtain by l'Hospital's rule
\be
\lim_{\lambda\nearrow1}\frac{\phi_{\lambda}(x)}{1-\lambda} \ = \ (\aal+\bal x)\left[- \log\left(\frac{\aal+\bal x}{\ahy+\bhy x}\right)+1\right]-(\ahy+\bhy x).\notag
\ee
Accordingly, let us define $z^{*}:=\arg\max_{x\in\mathbb{N}_{0}}\left\{(\aal+\bal x)\left[-\log\left(\frac{\aal+\bal x}{\ahy+\bhy x}\right)+1\right] -(\ahy+\bhy x)\right\}$ (note that the maximum exists since $\lim_{x\rightarrow\infty}\left\{(\aal+\bal x)\left[-\log\left(\frac{\aal+\bal x}{\ahy+\bhy x}\right)+1\right] -(\ahy+\bhy x)\right\}$ $=-\infty$). 
Due to continuity of the function $(\lambda,x) \mapsto  \frac{\phi_{\lambda}(x)}{1-\lambda}$,
there exists an $\epsilon>0$ such that for all $\lambda\in]1-\epsilon,1[$ it holds $z_{\lambda}^{*}=z^{*}$. Applying these considerations, we get with l'Hospital's rule 
\be\label{fo.ELhorn2}
I(\Pna||\Pnh) \ \geq \ \lim_{\lambda\nearrow1}\frac{1-\exp\left\{\phi_{\lambda}(z^{*})\cdot n\right\}}{\lambda(1-\lambda)} \ = \ \left[ f_{\Al}(z^{*})\cdot\left[\log\left(\frac{f_{\Al}(z^{*})}{f_{\Hy}(z^{*})}\right)-1\right]+f_{\Hy}(z^{*})\right]  \cdot \, n 
\ \geq \ 0.
\ee
In fact, in the current parameter constellation $\quasetSPdreiab\cup\quasetSPdreic$ 
we have 
$\phi_{\lambda}(x) < 0$ for all $\lambda \in ]0,1[$ and all $x \in \mathbb{N}_{0}$
which implies $f_{\Al}(z^{*})\neq f_{\Hy}(z^{*})$ by Lemma \ref{lem2}; thus,   
we even get $E^{L,hor}_{n}>0$ for all $n \in \mathbb{N}$ by virtue of the inequality $-\log\left(\frac{f_{\Hy}(z^{*})}{f_{\Al}(z^{*})}\right)>-\frac{f_{\Hy}(z^{*})}{f_{\Al}(z^{*})}+1$.\\[0.2cm]
\indent For the case $\quasetSPzwei$, the abovementioned procedure leads to $z_{\lambda}^{*}=0=z^{*}$ ($\lambda \in ]0,1[$)
which implies $\phi_{\lambda}(z_{\lambda}^{*})=0$,
$B^{U}_{\lambda,n}=1$ 
and thus  the trivial lower bound
$E^{L,hor}_{n} = \lim_{\lambda\nearrow1}\frac{1-B^{U}_{\lambda,n}}{\lambda(1-\lambda)}=0$ follows. 
In contrast, for the case $\quasetSPdreid$ one gets $z_{\lambda}^{*}= 
\frac{\aal-\ahy}{\bhy-\bal} =z^{*}$ ($\lambda \in ]0,1[$) which nevertheless
also implies $\phi_{\lambda}(z_{\lambda}^{*})=0$ and hence $E^{L,hor}_{n}=0$.
On $\quasetSPvier$, we have $\sup_{x\in\mathbb{N}_{0}}\phi_{\lambda}(x) = \phi_{\lambda}(\infty) =0$
and hence we set $E^{L,hor}_{n}:=0$.\\[0.2cm]
\indent To show the strict positivity $E^{L}_{n} >0$ in the parameter case $\quasetSPzwei$, we inspect the bound $E^{L,sec}_{0,n}$.
With the help of $\alpha := \alpha_{\bullet}:=\aal=\ahy$ (the bullet will be omitted in this proof) 
and the auxiliary variable $x:=\frac{\bhy}{\bal} >0$, the definition \eqref{fo.ELkn}
respectively its special case \eqref{fo.entboulowSP23ab} rewrites as  
\be\label{fo.h}
E^{L,sec}_{0,n} \ = \ E^{L,sec}_{0,n}(x) \ = \ \left\{
\begin{array}{ll}
\left[-(\alpha+\bal) \cdot \log\left(\frac{\alpha+\bal x}{\alpha+\bal}\right)+\bal(x-1)\right]\cdot\frac{1-(\bal)^{n}}{1-\bal}\cdot\left[\omega_{0}-\frac{\alpha}{1-\bal}\right] & \vspace{0.1cm}\\
+ \ \Big[\frac{\alpha}{\bal(1-\bal)}\left(-(\alpha+\bal) \cdot \log\left(\frac{\alpha+\bal x}{\alpha+\bal}\right)+\bal(x-1)\right) & \vspace{0.1cm}\\
\hspace{0.5cm}+\frac{\alpha}{\bal}\left(\alpha+\bal\right) \cdot \log\left(\frac{\alpha+\bal x}{\alpha+\bal}\right)-\alpha(x-1)\Big]\cdot n \, , & \hspace{-1.2cm}\textrm{if} \ \bal\neq1, \vspace{0.3cm}\\
\left[-(\alpha+1) \cdot \log\left(\frac{\alpha+x}{\alpha+1}\right)+x-1\right]\cdot\left[\frac{\alpha}{2}\cdot n^{2}+\left(\omega_{0}+\frac{\alpha}{2}\right)\cdot n\right] & \vspace{0.1cm}\\
+\left[(\alpha+1) \cdot \log\left(\frac{\alpha+x}{\alpha+1}\right)-x+1\right]\cdot\alpha\cdot n \, , & \hspace{-1.2cm}\textrm{if} \ \bal=1.
\end{array}
\right.
\ee
To prove that $E^{L,sec}_{0,n}>0$ for all 
$\omega_{0} \in \mathbb{N}$ and all $n \in \mathbb{N}$ it suffices to show that $E^{L,sec}_{0,n}(1)=\left(\frac{\partial}{\partial x}E^{L,sec}_{0,n}\right)(1)=0$ and $\left(\frac{\partial^2}{\partial x^2}E^{L,sec}_{0,n}\right)(x)>0$ for all $x\in ]0,\infty[\backslash\{1\}$. The assertion $E^{L,sec}_{0,n}(1)=0$ is trivial from \eqref{fo.h}. Moreover, 
we obtain
\be
\left(\frac{\partial}{\partial x}E^{L,sec}_{0,n}\right)(x) 
\ = \ 
\left\{
\begin{array}{ll}
\bal\cdot\left[
1-\frac{\alpha+\bal}{\alpha+\bal x}\right]\cdot\frac{1-(\bal)^{n}}{1-\bal}\cdot\left[\omega_{0}-\frac{\alpha}{1-\bal}
\right] & \vspace{0.1cm}\\
+ \ \alpha\cdot 
\left(1-\frac{\alpha+\bal}{\alpha+\bal x}\right) \cdot \frac{\bal}{1-\bal}
\cdot n \, , & \textrm{if} \ \bal\neq1, \vspace{0.3cm}\\
\left[1-\frac{\alpha+1}{\alpha+x}\right]\cdot\left[\frac{\alpha}{2}\cdot n^{2}+\left(\omega_{0}-\frac{\alpha}{2}\right)\cdot n\right] \, ,  & \textrm{if} \ \bal=1,
\end{array}
\right.
\ee
which immediately yields $\left(\frac{\partial}{\partial x}E^{L,sec}_{0,n}\right)(1)=0$. For the second derivative we get
\be\label{fo.htwoprime}
\left(\frac{\partial^2}{\partial x^2}E^{L,sec}_{0,n}\right)(x) \ = \ \left\{
\begin{array}{ll}
\frac{(\alpha+\bal)\cdot\bal^{2}}{(\alpha+\bal x)^{2}}\cdot\frac{1-(\bal)^{n}}{1-\bal}\cdot\left[\omega_{0}-\frac{\alpha}{1-\bal}\right] & \vspace{0.1cm}\\
+ \ \alpha \frac{\alpha+\bal}{(\alpha+\bal x)^{2}}\cdot\frac{\bal^2}{1-\bal}\cdot n \, >0, & \textrm{if} \ \bal\neq1, \vspace{0.3cm}\\
\frac{\alpha+1}{(\alpha+x)^{2}}\cdot\left[\frac{\alpha}{2}\cdot n^{2}+\left(\omega_{0}-\frac{\alpha}{2}\right)\cdot n\right] \ > 0, & \textrm{if} \ \bal=1,
\end{array}
\right.
\ee
where the strict positivity in the case $\bal\neq1$ follows immediately by replacing $\omega_{0}$ with $1$ and by using
the obvious relation $\frac{1}{1-\bal}\cdot\left[n-\frac{1-\bal^{n}}{1-\bal}\right]=\frac{1}{1-\bal}\sum_{k=0}^{n-1}\left(1-\bal^{k}\right)>0$.\\[0.2cm]
\indent For the constellation $\quasetSPvier$ with parameters $\beta := \beta_{\bullet} := \bal=\bhy$, $\aal \ne \ahy$,
the strict positivity of $E^{L}_{n} >0$ follows by showing 
that 
$E^{L,tan}_{y,n}$ converges from above to zero as $y$ tends to infinity. In fact,
there holds $\lim_{y\rightarrow\infty}y\cdot E^{L,tan}_{y,n} \in ]0,\infty[$. To see this, let us first observe that by
l'Hospital's rule we get 
\be
\lim_{y\rightarrow\infty}y\cdot\log\left(\frac{\aal+\beta y}{\ahy+\beta y}\right)\ = \ \frac{\aal-\ahy}{\beta}
\qquad \textrm{as well as} \qquad 
\lim_{y\rightarrow\infty}y\cdot\left(1-\frac{\aal+\beta y}{\ahy+\beta y}\right) \ = \ -\frac{\aal-\ahy}{\beta} \, .
\notag
\ee
From this and \eqref{fo.ELyn}, we obtain
\ $
\lim_{y\rightarrow\infty}y\cdot E^{L,tan}_{y,n} 
= \frac{(\aal-\ahy)^{2}}{\beta}\cdot n \ > \ 0
$ \ 
in both cases $\beta\neq1$ and $\beta=1$.\\[0.2cm]
\indent Finally, in the parameter case $\quasetSPdreid$ we consider the bound $E^{L,tan}_{y^{*},n}$, with $y^{*}=\frac{\aal-\ahy}{\bhy-\bal}$. Since $\aal+\bal y^{*}=\ahy+\bhy y^{*}$ holds, it is easy to see that $E^{L,tan}_{y^{*},n}=0$ for all $n\in\mathbb{N}$. However, 
 the condition $\left(\frac{\partial}{\partial y}E^{L,tan}_{y,n}\right)(y^{*})\neq0$ implies that $\sup_{y\geq0}E^{L,tan}_{y,n}>0$. The explicit form \eqref{fo.deryELtanystar} of this condition
follows from 
\be\label{fo.deryELtan}
\left(\frac{\partial}{\partial y}E^{L,tan}_{y,n}\right)(y) \ = \ \left\{
\begin{array}{ll}
 \frac{(\aal\bhy-\ahy\bal)^{2}}{f_{\Al}(y)\left(f_{\Hy}(y)\right)^{2}}\cdot\frac{1-\left(\bal\right)^{n}}{1-\bal}\cdot\left[\omega_{0}-\frac{\aal}{1-\bal}\right] & \\
+ \ \frac{\aal\bhy-\ahy\bal}{\left(f_{\Hy}(y)\right)^{2}}\cdot\left[\frac{\aal}{\bal(1-\bal)f_{\Al}(y)}-\frac{\aal\bhy-\ahy\bal}{\bal}\right]\cdot n 
\, , & \textrm{if} \ \bal\neq1,\vspace{0.2cm}\\ 
\frac{(\aal\bhy-\ahy)^{2}}{f_{\Al}(y)\left(f_{\Hy}(y)\right)^{2}}\cdot\left[\frac{\aal}{2}\cdot n^{2}+\left(\omega_{0}+\frac{\aal}{2}\right)\cdot n\right] \ - \ \frac{(\aal\bhy-\ahy)^{2}}{\left(f_{\Hy}(y)\right)^{2}}\cdot n \, , & \textrm{if} \ \bal=1,
\end{array}
\right.\notag
\ee
$y\geq0$, by using the particular choice $y=y^{*}$ together with $f_{\Al}(y^{*})=f_{\Hy}(y^{*})=-\frac{\aal\bhy-\ahy\bal}{\bal-\bhy}$  . 
\qed\\

\noindent
The next lemma (and parts of its proof) will be useful for the verification of Theorem \ref{thm.entdiflim}:

\begin{lem}\label{lem.entlim}
Recall the bounds on the Hellinger integral $m-$limit given in \eqref{fo.boulim1} and \eqref{fo.boulim} of Theorem \ref{thmlimit}, 
in terms of $L_{\lambda}^{(i)}(t)$ and $U_{\lambda}^{(i)}(t)$ ($i=1,2$) defined by \eqref{def.Leins} to \eqref{def.Uzwei}.
Correspondingly, one gets the following $\lambda-$limits for all $t\in[0,\infty[$:
\begin{itemize}
	\item[(a)] for all $\kappa_{\Al} \in ]0,\infty[$ and all $\kappa_{\Hy} \in [0,\infty[$ with $\kappa_{\Al} \ne \kappa_{\Hy}$
	\be \label{part.all}
	\lim_{\lambda\nearrow1}\frac{\partial  L_{\lambda}^{(1)}(t)}{\partial\lambda} 
	\ = \ \lim_{\lambda\nearrow1}\frac{\partial L_{\lambda}^{(2)}(t)}{\partial\lambda}  
	\ = \ \lim_{\lambda\nearrow1}\frac{\partial \, U_{\lambda}^{(1)}(t)}{\partial\lambda} 
	\ = \ \lim_{\lambda\nearrow1}\frac{\partial \, U_{\lambda}^{(2)}(t)}{\partial\lambda}  \ = \ 0 \, .
	\ee 
	\item[(b)] for $\kappa_{\Al}=0$  and all $\kappa_{\Hy} \in ]0,\infty[$ 
		\bea	
		&& \lim_{\lambda\nearrow1}\frac{\partial L_{\lambda}^{(1)}(t)}{\partial\lambda}  \ = \ -\frac{\kappa_{\Hy}^{2}\cdot t}{2\sigma^{2}} \  .
		\label{part.low1}\\ 
		&&  
    \lim_{\lambda\nearrow1}\frac{\partial L_{\lambda}^{(2)}(t)}{\partial\lambda}  \ = \ -\frac{\kappa_{\Hy}^{2}\cdot t^{2}}{4} \ .
    \label{part.low2} \\
		&& 	
		\lim_{\lambda\nearrow1}\frac{\partial \, U_{\lambda}^{(1)}(t)}{\partial\lambda} 
		\ = \ 
		\lim_{\lambda\nearrow1}\frac{\partial \, U_{\lambda}^{(2)}(t)}{\partial\lambda} \ = \ 0 \, .
		\label{part.up}
		\eea
\end{itemize}
\end{lem}

\noindent 
\prl\ref{lem.entlim} \ For
all $\kappa_{\Al},  \kappa_{\Hy} \in [0,\infty[$ with $\kappa_{\Al} \ne \kappa_{\Hy}$ one can deduce from 
\eqref{def.kl-Lambda} as well as \eqref{def.Leins} to \eqref{def.Uzwei}
the following derivatives:
\bea\label{prlem.entlim1}
&& \hspace{-0.6cm}\frac{\partial  L_{\lambda}^{(1)}(t)}{\partial\lambda} \ = \ 
\frac{1}{2\sigma^2}
\Bigg\{\frac{t}{2}\left(\frac{\Lambda_{\!\lambda}-\kappa_{\lambda}}{\Lambda_{\!\lambda}}\right)^2 \left(\kappa_{\Al}^{2}-\kappa_{\Hy}^{2}\right)\left[2e^{-2\Lambda_{\!\lambda} t}-e^{-\Lambda_{\!\lambda} t}\right]\notag\\
&& \hspace{1.2cm}+ \ e^{-\Lambda_{\!\lambda} t}  \frac{1-e^{-\Lambda_{\!\lambda} t}}{\Lambda_{\!\lambda}} \left[\frac{\Lambda_{\!\lambda}-\kappa_{\lambda}}{\Lambda_{\!\lambda}}\left(\kappa_{\Al}^{2}-\kappa_{\Hy}^{2}-2\Lambda_{\!\lambda}(\kappa_{\Al}-\kappa_{\Hy})\right)-
\left(\frac{\Lambda_{\!\lambda}-\kappa_{\lambda}}{\Lambda_{\!\lambda}}\right)^2 \frac{\kappa_{\Al}^{2}-\kappa_{\Hy}^{2}}{2} \right] \Bigg\} \ ,
\eea
\vspace{-0.5cm}
\bea\label{prlem.entlim2}
\frac{\partial L_{\lambda}^{(2)}(t)}{\partial\lambda} & = &  
\frac{1}{4} \Bigg\{\frac{\Lambda_{\!\lambda}-\kappa_{\lambda}}{\Lambda_{\!\lambda}}\cdot\left(\frac{1-e^{-\Lambda_{\!\lambda} t}}{\Lambda_{\!\lambda}}\right)^{2}\cdot\left(\kappa_{\Al}^{2}-\kappa_{\Hy}^{2}-2\Lambda_{\!\lambda}(\kappa_{\Al}-\kappa_{\Hy})-\frac{\Lambda_{\!\lambda}-\kappa_{\lambda}}{\Lambda_{\!\lambda}}\left(\kappa_{\Al}^{2}-\kappa_{\Hy}^{2}\right)\right)\notag\\
&& \hspace{1.0cm} + \ t \cdot e^{-\Lambda_{\!\lambda} t}\cdot\left(\frac{\Lambda_{\!\lambda}-\kappa_{\lambda}}{\Lambda_{\!\lambda}}\right)^{2}\cdot\frac{1-e^{-\Lambda_{\!\lambda} t}}{\Lambda_{\!\lambda}}\cdot\left(\kappa_{\Al}^{2}-\kappa_{\Hy}^{2}\right)\Bigg\} \ ,\notag\\[-1.1cm]
&& \notag
\eea
\vspace{-0.5cm}
\bea\label{prlem.entlim3}
&& 
\frac{\partial \, U_{\lambda}^{(1)}(t)}{\partial\lambda} \ = \ 
\frac{1}{\sigma^2} \Bigg\{\frac{1}{2} \cdot \frac{\Lambda_{\!\lambda}-\kappa_{\lambda}}{\Lambda_{\!\lambda}}
\bigg[
t\, e^{-\Lambda_{\!\lambda} t}\left(\kappa_{\Al}^{2}-\kappa_{\Hy}^{2}\right)-\frac{t}{2} \  e^{-\frac{1}{2}(\Lambda_{\!\lambda}+\kappa_{\lambda})t}\left(\kappa_{\Al}^{2}-\kappa_{\Hy}^{2}+2\Lambda_{\!\lambda}(\kappa_{\Al}-\kappa_{\Hy})\right)
\notag\\
&& \hspace{5.2cm} - \   
 e^{-\frac{1}{2}(\Lambda_{\!\lambda}+\kappa_{\lambda})t} \cdot 
 \frac{1-e^{-\Lambda_{\!\lambda} t}}{\Lambda_{\!\lambda}} \cdot 
\left(\kappa_{\Al}^{2}-\kappa_{\Hy}^{2}+2\Lambda_{\!\lambda}(\kappa_{\Al}-\kappa_{\Hy})\right)
\bigg] \notag\\
&& \hspace{0.5cm}+ \ \frac{1}{2} \ \frac{e^{-\frac{1}{2}(\Lambda_{\!\lambda}+\kappa_{\lambda})t}-e^{-\Lambda_{\!\lambda} t}}{\Lambda_{\!\lambda}}\left(\kappa_{\Al}^{2}-\kappa_{\Hy}^{2}-2\Lambda_{\!\lambda}(\kappa_{\Al}-\kappa_{\Hy})\right)\notag\\
&& \hspace{0.5cm}+ \ \frac{1}{4}\left(\frac{\Lambda_{\!\lambda}-\kappa_{\lambda}}{\Lambda_{\!\lambda}}\right)^{2}\Bigg[\frac{t}{2} \ e^{-\frac{1}{2}(\Lambda_{\!\lambda}+\kappa_{\lambda})t}\left(\kappa_{\Al}^{2}-\kappa_{\Hy}^{2}+2\Lambda_{\!\lambda}(\kappa_{\Al}-\kappa_{\Hy})\right)\notag\\
&& \hspace{0.5cm} -\frac{t}{2} \ e^{-\frac{1}{2}(3\Lambda_{\!\lambda}+\kappa_{\lambda})t}\left(3\left(\kappa_{\Al}^{2}-\kappa_{\Hy}^{2}\right)+2\Lambda_{\!\lambda}(\kappa_{\Al}-\kappa_{\Hy})\right)
+ \ e^{-\frac{1}{2}(\Lambda_{\!\lambda}+\kappa_{\lambda})t} \cdot \frac{1-e^{-\Lambda_{\!\lambda} t}}{\Lambda_{\!\lambda}} \cdot \left(\kappa_{\Al}^{2}-\kappa_{\Hy}^{2}\right)\Bigg] \Bigg\}
\ , \notag\\[-0.5cm]
&&\notag
\eea
\vspace{-0.7cm}
\bea\label{prlem.entlim4}
\frac{\partial \, U_{\lambda}^{(2)}(t)}{\partial\lambda}& = &   
\Bigg\{\frac{\left(\Lambda_{\!\lambda}-\kappa_{\lambda}\right)^{2}}{\Lambda_{\!\lambda}(3\Lambda_{\!\lambda}+\kappa_{\lambda})}\Bigg[\frac{t}{2} \ 
e^{-\frac{1}{2}(3\Lambda_{\!\lambda}+\kappa_{\lambda})t}\left(3\frac{\kappa_{\Al}^{2}-\kappa_{\Hy}^{2}}{2\Lambda_{\!\lambda}}+\kappa_{\Al}-\kappa_{\Hy}\right)\notag\\
&& \hspace{3.0cm} - \ \frac{1-e^{-\frac{1}{2}(3\Lambda_{\!\lambda}+\kappa_{\lambda})t}}{3\Lambda_{\!\lambda}+\kappa_{\lambda}}
\cdot
\left(3 \, \frac{\kappa_{\Al}^{2}-\kappa_{\Hy}^{2}}{2\Lambda_{\!\lambda}}+\kappa_{\Al}-\kappa_{\Hy}\right)\Bigg]\notag\\
&& \hspace{1.0cm}+ \ \frac{\Lambda_{\!\lambda}-\kappa_{\lambda}}{\Lambda_{\!\lambda}}\left[\frac{t}{2} \ e^{-\frac{1}{2}(\Lambda_{\!\lambda}+\kappa_{\lambda})t}\left(\frac{\kappa_{\Al}^{2}-\kappa_{\Hy}^{2}}{2\Lambda_{\!\lambda}}+\kappa_{\Al}-\kappa_{\Hy}\right) - t \, e^{-\Lambda_{\!\lambda} t} \, \frac{\kappa_{\Al}^{2}-\kappa_{\Hy}^{2}}{2\Lambda_{\!\lambda}}\right]\notag\\
&& \hspace{1cm}+ \ \frac{e^{-\frac{1}{2}(\Lambda_{\!\lambda}+\kappa_{\lambda})t}-e^{-\Lambda_{\!\lambda} t}}{\Lambda_{\!\lambda}}\left(\frac{\kappa_{\Al}^{2}-\kappa_{\Hy}^{2}}{2\Lambda_{\!\lambda}}-\kappa_{\Al}+\kappa_{\Hy}\right)\notag\\
&& \hspace{1.0cm} + \ \left(2\left(\frac{\kappa_{\Al}^{2}-\kappa_{\Hy}^{2}}{2\Lambda_{\!\lambda}}-\kappa_{\Al}+\kappa_{\Hy}\right)-\frac{\Lambda_{\!\lambda}-\kappa_{\lambda}}{\Lambda_{\!\lambda}^{2}}\cdot\frac{\kappa_{\Al}^{2}-\kappa_{\Hy}^{2}}{2}\right)\notag\\
&& \hspace{1.6cm}\cdot \ \frac{1}{\Lambda_{\!\lambda}}\left[\frac{\Lambda_{\!\lambda}-\kappa_{\lambda}}{3\Lambda_{\!\lambda}+\kappa_{\lambda}}\left(1-e^{-\frac{1}{2}(3\Lambda_{\!\lambda}+\kappa_{\lambda})t}\right)-e^{-\frac{1}{2}(\Lambda_{\!\lambda}+\kappa_{\lambda})t}+e^{-\Lambda_{\!\lambda} t}\right] \ .
\eea
If $\kappa_{\Al} \in ]0,\infty[$ and $\kappa_{\Hy} \in [0,\infty[$ with $\kappa_{\Al} \ne \kappa_{\Hy}$,
then one gets $\lim_{\lambda\nearrow1} \Lambda_{\!\lambda} = \lim_{\lambda\nearrow1} \kappa_{\lambda} = \kappa_{\Al} >0$ which implies  
\eqref{part.all} from \eqref{prlem.entlim1} to \eqref{prlem.entlim4}. For the proof of part (b), let us
correspondingly assume $\kappa_{\Al}=0$ and $\kappa_{\Hy} \in ]0,\infty[$,
which by \eqref{def.kl-Lambda} leads to $\kappa_{\lambda} = \kappa_{\Hy}\cdot (1-\lambda)$,
$\Lambda_{\!\lambda}=  \kappa_{\Hy}\cdot \sqrt{1-\lambda}$ and the convergences 
$\lim_{\lambda\nearrow1}\Lambda_{\!\lambda}=\lim_{\lambda\nearrow1}\kappa_{\lambda}=0$. 
From this, the assertions \eqref{part.low1}, \eqref{part.low2}, \eqref{part.up}
follow in a straightforward manner from \eqref{prlem.entlim1}, \eqref{prlem.entlim2}, \eqref{prlem.entlim3} -- respectively --
by using (parts of) the obvious relations 
	\be \label{limhilf1}
	\lim_{\lambda\nearrow1} \frac{\kappa_{\lambda}}{\Lambda_{\!\lambda}} =0, \qquad  
	\lim_{\lambda\nearrow1}\frac{\Lambda_{\!\lambda}\pm\kappa_{\lambda}}{\Lambda_{\!\lambda}}    
	\ = \ \lim_{\lambda\nearrow1}\frac{\Lambda_{\!\lambda}-\kappa_{\lambda}}{\Lambda_{\!\lambda}+\kappa_{\lambda}} \ = \ 1 \ , 
	\ee
	\be \label{limhilf2}
	\lim_{\lambda\nearrow1}\frac{1-e^{-c_{\lambda} \cdot t}}{c_{\lambda}} 
	 \ = \ t \qquad \text{for all } c_{\lambda} \in \left\{\Lambda_{\!\lambda}, \frac{\Lambda_{\!\lambda}+\kappa_{\lambda}}{2}, 
	 \frac{3\, \Lambda_{\!\lambda}+\kappa_{\lambda}}{2}  \right\} \ .
	\ee

\noindent In order to get the last assertion \eqref{part.up} we make use of the following limits
	\be  \label{limhilf3}
\lim_{\lambda\nearrow1}\frac{1}{\Lambda_{\!\lambda}-\kappa_{\lambda}}
-\frac{3}{3\Lambda_{\!\lambda}+\kappa_{\lambda}}  \ = \ \lim_{\lambda\nearrow1}
\frac{4 \, \kappa_{\Hy}}{(\kappa_{\Hy}-\kappa_{\Hy}\cdot \sqrt{1-\lambda}) \cdot 
(3\, \kappa_{\Hy}+\kappa_{\Hy}\cdot \sqrt{1-\lambda}) }
\ = \ \frac{4}{3\, \kappa_{\Hy}} \ 
	\ee 
and
	\be	\label{limhilf4}
\lim_{\lambda\nearrow1}\frac{1}{\Lambda_{\!\lambda}}\left[\frac{1-e^{-\frac{1}{2}(3\Lambda_{\!\lambda}+\kappa_{\lambda})t}}{3\Lambda_{\!\lambda}+\kappa_{\lambda}}-\frac{1-e^{-\Lambda_{\!\lambda} t}}{\Lambda_{\!\lambda}-\kappa_{\lambda}}+\frac{1-e^{-\frac{1}{2}(\Lambda_{\!\lambda}+\kappa_{\lambda})t}}{\Lambda_{\!\lambda}-\kappa_{\lambda}}\right] \ = \ 0 \, .
	\ee		

\noindent
To see \eqref{limhilf4}, let us first observe that the involved limit can be rewritten as
	\bea
& & \lim_{\lambda\nearrow1}\Bigg\{\frac{1}{\Lambda_{\!\lambda}(\Lambda_{\!\lambda}-\kappa_{\lambda})}\left[\frac{1}{3}-\frac{1}{3}~e^{-\frac{1}{2}(3\Lambda_{\!\lambda}+\kappa_{\lambda})t}+e^{-\Lambda_{\!\lambda} t}-e^{-\frac{1}{2}(\Lambda_{\!\lambda}+\kappa_{\lambda})t}\right]\label{prlem.entlim6}\\
	&& \hspace{1cm} + \ \frac{1-e^{-\frac{1}{2}(3\Lambda_{\!\lambda}+\kappa_{\lambda})t}}{\Lambda_{\!\lambda}}\left[\frac{1}{3\Lambda_{\!\lambda}+\kappa_{\lambda}}-\frac{1}{3(\Lambda_{\!\lambda}-\kappa_{\lambda})}\right]\Bigg\} \ .
	\label{prlem.entlim6b}
	\eea
Substituting $x:=\sqrt{1-\lambda}$ and applying l'Hospital's rule twice, we get for the first limit \eqref{prlem.entlim6}
\bea
&  &\hspace{-0.7cm} ~~  \lim_{x\searrow0}\frac{\frac{1}{3}-\frac{1}{3} \ e^{-\frac{\kappa_{\Hy}t}{2}(3x+x^{2})}+e^{-\kappa_{\Hy}tx}-e^{-\frac{\kappa_{\Hy}t}{2}(x+x^{2})}}{\kappa_{\Hy}^{2}\cdot\left(x^{2}-x^{3}\right)}\notag\\
&&\hspace{-0.7cm} =  \lim_{x\searrow0}\frac{\frac{\kappa_{\Hy}t}{6}(3+2x) \ e^{-\frac{\kappa_{\Hy}t}{2}(3x+x^{2})}-\kappa_{\Hy}\, t\, e^{-\kappa_{\Hy}tx}+\frac{\kappa_{\Hy}t}{2} \, (1+2x) \, e^{-\frac{\kappa_{\Hy}t}{2}(x+x^{2})}}{\kappa_{\Hy}^{2}\cdot\left(2x-3x^{2}\right)}\notag\\
&&\hspace{-0.7cm} =  \lim_{x\searrow0}\frac{\left[-\frac{\kappa_{\Hy}^{2}t^{2}}{12}(3+2x)^{2}+\frac{\kappa_{\Hy}t}{3}\right]e^{-\frac{\kappa_{\Hy}t}{2}(3x+x^{2})}+\kappa_{\Hy}^{2} \, t^{2}\, e^{-\kappa_{\Hy}tx}\hspace{-0.1cm}-\hspace{-0.1cm}\left[\frac{\kappa_{\Hy}^{2}t^{2}}{4}(1+2x)^{2}-\kappa_{\Hy}\, t\right]e^{-\frac{\kappa_{\Hy}t}{2}(x+x^{2})}}{\kappa_{\Hy}^{2}\cdot\left(2-6x\right)}\notag\\
&&\hspace{-0.7cm} =  \frac{1}{2\kappa_{\Hy}^{2}}\left[-\frac{3\kappa_{\Hy}^{2}t^{2}}{4}+\frac{\kappa_{\Hy}t}{3}+\kappa_{\Hy}^{2}t^{2}-\frac{\kappa_{\Hy}^{2}t^{2}}{4}+\kappa_{\Hy}t\right] \ = \ \frac{2t}{3\, \kappa_{\Hy}}. 
\notag
\eea
The second limit \eqref{prlem.entlim6b} becomes
\bea
& & \lim_{\lambda\nearrow1}\frac{1-e^{-\frac{1}{2}(3\Lambda_{\!\lambda}+\kappa_{\lambda})t}}{3\Lambda_{\!\lambda}+\kappa_{\lambda}}
\cdot
\frac{3\Lambda_{\!\lambda}+\kappa_{\lambda}}{\Lambda_{\!\lambda}}
\cdot\frac{-4\kappa_{\Hy}}{(3\kappa_{\Hy}+\sqrt{1-\lambda}\kappa_{\Hy})(3\kappa_{\Hy}+3\sqrt{1-\lambda}\kappa_{\Hy})}
 \ = \ -\frac{2t}{3\, \kappa_{\Hy}} \ , \notag
\eea
and consequently \eqref{limhilf4} follows. To proceed with the proof of \eqref{part.up},
we rearrange
\bea
&& \hspace{-0.5cm}\lim_{\lambda\nearrow1}\frac{\partial \, U_{\lambda}^{(2)}(t)}{\partial\lambda}
\ = \ \lim_{\lambda\nearrow1}\Bigg\{\left(\frac{\Lambda_{\!\lambda}-\kappa_{\lambda}}{\Lambda_{\!\lambda}}\right)^{2}\Bigg[\frac{\Lambda_{\!\lambda}}{3\Lambda_{\!\lambda}+\kappa_{\lambda}}\left(\frac{t}{2} \ e^{-\frac{1}{2}(3\Lambda_{\!\lambda}+\kappa_{\lambda})t}\left(-\frac{3\kappa_{\Hy}^{2}}{2\Lambda_{\!\lambda}}-\kappa_{\Hy}\right)\right)\notag\\
&& \hspace{-0.5cm} 
- \ \frac{\Lambda_{\!\lambda}}{3\Lambda_{\!\lambda}+\kappa_{\lambda}} \cdot \frac{1-e^{-\frac{1}{2}(3\Lambda_{\!\lambda}+\kappa_{\lambda})t}}{3\Lambda_{\!\lambda}+\kappa_{\lambda}}\left(-\frac{3\kappa_{\Hy}^{2}}{2\Lambda_{\!\lambda}}-\kappa_{\Hy}\right) \ + \ \frac{\Lambda_{\!\lambda}}{\Lambda_{\!\lambda}-\kappa_{\lambda}}\frac{e^{-\frac{1}{2}(\Lambda_{\!\lambda}+\kappa_{\lambda})t}-e^{-\Lambda_{\!\lambda} t}}{\Lambda_{\!\lambda}-\kappa_{\lambda}}\left(-\frac{\kappa_{\Hy}^{2}}{2\Lambda_{\!\lambda}}+\kappa_{\Hy}\right)\notag\\
&& \hspace{-0.5cm} 
- \ \frac{\Lambda_{\!\lambda}}{\Lambda_{\!\lambda}-\kappa_{\lambda}}\left(-\frac{t}{2} \ e^{-\frac{1}{2}(\Lambda_{\!\lambda}+\kappa_{\lambda})t}\left(-\frac{\kappa_{\Hy}^{2}}{2\Lambda_{\!\lambda}}-\kappa_{\Hy}\right)-t\ e^{-\Lambda_{\!\lambda} t}\frac{\kappa_{\Hy}^{2}}{2\Lambda_{\!\lambda}}\right)\Bigg]\notag\\
&& \hspace{-0.5cm} 
 +  \left[\frac{\Lambda_{\!\lambda}-\kappa_{\lambda}}{\Lambda_{\!\lambda}}\left(-\kappa_{\Hy}^{2}+2\Lambda_{\!\lambda}\kappa_{\Hy}\right)+\left(\frac{\Lambda_{\!\lambda}-\kappa_{\lambda}}{\Lambda_{\!\lambda}}\right)^{2}\frac{\kappa_{\Hy}^{2}}{2}\right]\cdot\left[\frac{1-e^{-\frac{1}{2}(3\Lambda_{\!\lambda}+\kappa_{\lambda})t}}{\Lambda_{\!\lambda}(3\Lambda_{\!\lambda}+\kappa_{\lambda})}
-\frac{e^{-\frac{1}{2}(\Lambda_{\!\lambda}+\kappa_{\lambda})t}-e^{-\Lambda_{\!\lambda} t}}{\Lambda_{\!\lambda}(\Lambda_{\!\lambda}-\kappa_{\lambda})}\right]\Bigg\}\notag
\eea
\be
  = \  \lim_{\lambda\nearrow1}\Bigg\{\left(\frac{\Lambda_{\!\lambda}-\kappa_{\lambda}}{\Lambda_{\!\lambda}}\right)^{2}\Bigg[
\frac{\kappa_{\Hy}^{2}\, t}{4} \left(-\frac{3\, e^{-\frac{1}{2}(3\Lambda_{\!\lambda}+\kappa_{\lambda})t}}{3\Lambda_{\!\lambda}+\kappa_{\lambda}}-\frac{e^{-\frac{1}{2}(\Lambda_{\!\lambda}+\kappa_{\lambda})t}}{\Lambda_{\!\lambda}-\kappa_{\lambda}}+
\frac{2\, e^{-\Lambda_{\!\lambda} t}}{\Lambda_{\!\lambda}-\kappa_{\lambda}}\right)\hspace{2cm}\label{prlem.entlim12}
\ee
\be
+ \ \frac{\kappa_{\Hy}^{2}}{2} \left(\frac{3\left(1-e^{-\frac{1}{2}(3\Lambda_{\!\lambda}+\kappa_{\lambda})t}\right)}{\left(3\Lambda_{\!\lambda}+\kappa_{\lambda}\right)^{2}}-\frac{1-e^{-\Lambda_{\!\lambda} t}}{\left(\Lambda_{\!\lambda}-\kappa_{\lambda}\right)^{2}}+\frac{1-e^{-\frac{1}{2}(\Lambda_{\!\lambda}+\kappa_{\lambda})t}}{\left(\Lambda_{\!\lambda}-\kappa_{\lambda}\right)^{2}}\right)\hspace{2cm}\label{prlem.entlim13}
\ee
\bea
&& \hspace{-0.6cm}
+ \ \kappa_{\Hy}\Bigg(
- \frac{\Lambda_{\!\lambda}}{3\Lambda_{\!\lambda}+\kappa_{\lambda}} \cdot \frac{t \ e^{-\frac{1}{2}(3\Lambda_{\!\lambda}+\kappa_{\lambda})t}}{2}+
\frac{\Lambda_{\!\lambda}}{3\Lambda_{\!\lambda}+\kappa_{\lambda}} \cdot 
\frac{1-e^{-\frac{1}{2}(3\Lambda_{\!\lambda}+\kappa_{\lambda})t}}{3\Lambda_{\!\lambda}+\kappa_{\lambda}}-
\frac{\Lambda_{\!\lambda}}{\Lambda_{\!\lambda}-\kappa_{\lambda}} \cdot
\frac{t\ e^{-\frac{1}{2}(\Lambda_{\!\lambda}+\kappa_{\lambda})t}}{2}\notag\\
&& \hspace{-0.6cm}
+ \ \frac{\Lambda_{\!\lambda}}{\Lambda_{\!\lambda}-\kappa_{\lambda}} \cdot \frac{1-e^{-\Lambda_{\!\lambda} t}}{\Lambda_{\!\lambda}-\kappa_{\lambda}}-
\frac{\Lambda_{\!\lambda}}{\Lambda_{\!\lambda}-\kappa_{\lambda}} \cdot
\frac{1-e^{-\frac{1}{2}(\Lambda_{\!\lambda}+\kappa_{\lambda})t}}{\Lambda_{\!\lambda}-\kappa_{\lambda}}\Bigg)\Bigg]\notag\\
&& \hspace{-0.6cm}
+   \left[\frac{\Lambda_{\!\lambda}-\kappa_{\lambda}}{\Lambda_{\!\lambda}}\left(-\kappa_{\Hy}^{2}+2\Lambda_{\!\lambda}\kappa_{\Hy}\right)+\left(\frac{\Lambda_{\!\lambda}-\kappa_{\lambda}}{\Lambda_{\!\lambda}}\right)^{2}\frac{\kappa_{\Hy}^{2}}{2}\right]\cdot\left[\frac{1-e^{-\frac{1}{2}(3\Lambda_{\!\lambda}+\kappa_{\lambda})t}}{\Lambda_{\!\lambda}(3\Lambda_{\!\lambda}+\kappa_{\lambda})}-\frac{e^{-\frac{1}{2}(\Lambda_{\!\lambda}+\kappa_{\lambda})t}-e^{-\Lambda_{\!\lambda} t}}{\Lambda_{\!\lambda}(\Lambda_{\!\lambda}-\kappa_{\lambda})}\right]\Bigg\}.\notag\\
&& \label{prlem.entlim14}
\eea
By means of \eqref{limhilf1} to \eqref{limhilf3}, the limit of the 
expression after the squared brackets in \eqref{prlem.entlim12} becomes
\be
 \lim_{\lambda\nearrow1}\Bigg\{\frac{\kappa_{\Hy}^{2}\, t}{4}
\left[ \frac{1-e^{-\frac{1}{2}(\Lambda_{\!\lambda}+\kappa_{\lambda})t}}{\Lambda_{\!\lambda}-\kappa_{\lambda}}
- 2\, \frac{1-e^{-\Lambda_{\!\lambda} t}}{\Lambda_{\!\lambda}-\kappa_{\lambda}} + 
3 \, \frac{1-e^{-\frac{1}{2}(3\Lambda_{\!\lambda}+\kappa_{\lambda})t}}{3\Lambda_{\!\lambda}+\kappa_{\lambda}}
+ \frac{1}{\Lambda_{\!\lambda}-\kappa_{\lambda}}
-\frac{3}{3\Lambda_{\!\lambda}+\kappa_{\lambda}} \right]
 =  \frac{\kappa_{\Hy}\, t}{3} \, , 
\label{limhilf7}
\ee
and the limit of the 
expression in \eqref{prlem.entlim13} becomes with \eqref{limhilf4}
\bea
&& \lim_{\lambda\nearrow1}\Bigg\{\frac{\Lambda_{\!\lambda}}{\Lambda_{\!\lambda}-\kappa_{\lambda}}\cdot\frac{\kappa_{\Hy}^{2}}{2\Lambda_{\!\lambda}}\cdot\left[\frac{1-e^{-\frac{1}{2}(3\Lambda_{\!\lambda}+\kappa_{\lambda})t}}{3\Lambda_{\!\lambda}+\kappa_{\lambda}}-\frac{1-e^{-\Lambda_{\!\lambda} t}}{\Lambda_{\!\lambda}-\kappa_{\lambda}}+\frac{1-e^{-\frac{1}{2}(\Lambda_{\!\lambda}+\kappa_{\lambda})t}}{\Lambda_{\!\lambda}-\kappa_{\lambda}}\right]\notag\\
&& \hspace{1cm} - \ \frac{\kappa_{\Hy}^{2}}{2}\cdot\frac{1-e^{-\frac{1}{2}(3\Lambda_{\!\lambda}+\kappa_{\lambda})t}}{3\Lambda_{\!\lambda}+\kappa_{\lambda}}
\cdot
\left[\frac{1}{\Lambda_{\!\lambda}-\kappa_{\lambda}}
-\frac{3}{3\Lambda_{\!\lambda}+\kappa_{\lambda}}\right]
\ = \ -\frac{\kappa_{\Hy}t}{3} \ .
\label{limhilf8}
\eea
By putting \eqref{prlem.entlim14}, \eqref{limhilf7}, \eqref{limhilf8} together with \eqref{limhilf4} we finally end up with 
\be
\lim_{\lambda\nearrow1}\frac{\partial \, U_{\lambda}^{(2)}(t)}{\partial\lambda} \ = \ \left[\frac{\kappa_{\Hy}t}{3}-\frac{\kappa_{\Hy}t}{3}\right]+\kappa_{\Hy}\left(-\frac{t}{6}+\frac{t}{6}-\frac{t}{2}+t-\frac{t}{2}\right)+\left[-\kappa_{\Hy}^{2}+\frac{\kappa_{\Hy}^{2}}{2}\right]\cdot 0 \ = \ 0 \, ,
\notag
\ee
which finishes the proof of Lemma \ref{lem.entlim}.
\qed\\

\noindent
\prt \ref{thm.entdiflim} \ 
Recall from \eqref{defparm} the approximative offspring-distribution-parameter $\beta^{(m)}_{\bullet} \, := \, 1-\frac{\kappa_{\bullet}}{\sigma^{2}m}$ and immigration-distribution parameter $\alpha^{(m)}_{\bullet} \, := \, \beta^{(m)}_{\bullet}\cdot\frac{\eta}{\sigma^{2}}$,
which is a special case of $\quasetNI\cup\quasetSPeins$.
Let us first calculate
 $\lim_{m\rightarrow\infty} I\left(P^{(m)}_{\Al,\left\lfloor \sigma^{2}mt\right\rfloor}\Big|\Big|P^{(m)}_{\Hy,\left\lfloor \sigma^{2}mt\right\rfloor}\right)$ by starting from Theorem \ref{thm.entex}(a). Correspondingly, we evaluate for all $\kappa_{\Al}\geq0$,
 $\kappa_{\Hy}\geq0$ with $\kappa_{\Al} \ne \kappa_{\Hy}$
\bea\label{fo.prthmentlim5}
&& \hspace{-0.7cm} \lim_{m\rightarrow\infty}m^{2}\cdot\left[\bam\cdot\left(\log\left(\frac{\bam}{\bhm}\right)-1\right)+\bhm\right]  =  \lim_{m\rightarrow\infty}\frac{-m}{2\sigma^{2}}\left[\kappa_{\Al}\log\left(\frac{\bam}{\bhm}\right)+\kappa_{\Hy}\left(1-\frac{\bam}{\bhm}\right)\right]\notag\\
&&  \quad = \ \frac{1}{2\sigma^{4}} \cdot \lim_{m\rightarrow\infty}\frac{\bhm\cdot\kappa_{\Al}-\bam\cdot\kappa_{\Hy}}{\left(\bhm\right)^{2}}\cdot\left(\kappa_{\Al}\cdot\frac{\bhm}{\bam}-\kappa_{\Hy}\right) \ = \ \frac{\left(\kappa_{\Al}-\kappa_{\Hy}\right)^{2}}{2\sigma^{4}} \, . \\[-0.5cm]
&& \notag
\eea
Additionally there holds
\be\label{fo.prthmentlim6}
\lim_{m\rightarrow\infty}m\cdot(1-\bam)  =  \frac{\kappa_{\Al}}{\sigma^{2}} \quad\textrm{and}\quad\lim_{m\rightarrow\infty}\left(\bam\right)^{\left\lfloor \sigma^{2}mt\right\rfloor}  =  \lim_{m\rightarrow\infty}\left[\left(1-\frac{\kappa_{\Al}}{\sigma^{2}m}\right)^{m}\right]^{\left\lfloor \sigma^{2}mt\right\rfloor/m}  =  e^{-\kappa_{\Al}\cdot t} \ .
\ee
For $\kappa_{\Al}>0$, we apply the upper part of formula \eqref{fo.entex} as well as
\eqref{fo.prthmentlim5}, \eqref{fo.prthmentlim6} to derive
\bea 
&& \lim_{m\rightarrow\infty}I_{\lambda}\left(P^{(m)}_{\Al,\left\lfloor \sigma^{2}mt\right\rfloor}\Big|\Big|P^{(m)}_{\Hy,\left\lfloor \sigma^{2}mt\right\rfloor}\right) \ = \ \lim_{m\rightarrow\infty}\left[\frac{m^{2}\cdot\left[\bam\cdot\left(\log\left(\frac{\bam}{\bhm}\right)-1\right)+\bhm\right]}{m\cdot(1-\bam)}\right.\notag\\
&& \ \cdot\left[\frac{X_{0}^{(m)}}{m}-\frac{\aam}{m\cdot(1-\bam)}\right]\cdot\left(1-\left(\bam\right)^{\left\lfloor \sigma^{2}mt\right\rfloor}\right)\notag\\
&& \ +\frac{\aam}{\bam\cdot m\cdot(1-\bam)}\cdot m^{2}\cdot\left[\bam\cdot\left(\log\left(\frac{\bam}{\bhm}\right)-1\right)+\bhm\right]\cdot\left.\frac{\left\lfloor \sigma^{2}mt\right\rfloor}{m}\right]\notag\\
&& = \ \frac{\left(\kappa_{\Al}-\kappa_{\Hy}\right)^{2}}{2\sigma^{2}\cdot\kappa_{\Al}}\cdot\left[\left(\widetilde{X}_{0}-\frac{\eta}{\kappa_{\Al}}\right)\cdot\left(1-e^{-\kappa_{\Al}\cdot t}\right)+\eta\cdot t\right].
\notag
\eea
For $\kappa_{\Al}=0$ (and thus $\kappa_{\Hy}>0$, $\bam \equiv 1$, $\aam \equiv \eta/\sigma^2$), we apply the lower part of formula \eqref{fo.entex}  
as well as
\eqref{fo.prthmentlim5}, \eqref{fo.prthmentlim6}
to obtain\vspace{-0.3cm}
\bea 
&& \lim_{m\rightarrow\infty}I_{\lambda}\left(P^{(m)}_{\Al,\left\lfloor \sigma^{2}mt\right\rfloor}\Big|\Big|P^{(m)}_{\Hy,\left\lfloor 
\sigma^{2}mt\right\rfloor}\right)  \ = \  \Bigg\{\lim_{m\rightarrow\infty}m^{2}\cdot\left[
\bhm - \log \bhm -1\right]\notag\\
&&
\cdot
\left[\frac{\eta}{2\sigma^{2}}\cdot\frac{\left(\left\lfloor \sigma^{2}mt\right\rfloor\right)^{2}}{m^{2}} \ + \ \left(\frac{X_{0}^{(m)}}{m}+\frac{\eta}{2\sigma^{2}\cdot m}\right)\cdot\frac{\left\lfloor \sigma^{2}mt\right\rfloor}{m}\right]\Bigg\} \ = \ \frac{\kappa_{\Hy}^{2}}{2\sigma^{2}}\cdot\left[\frac{\eta}{2}\cdot t^{2} \ + \ \widetilde{X}_{0}\cdot t\right].
\notag
\eea

\noindent
Let us now calculate the ``converse'' double limit \vspace{-0.2cm}
\be
\lim_{\lambda\nearrow1}\lim_{m\rightarrow\infty}I_{\lambda}\left(P^{(m)}_{\Al,\left\lfloor \sigma^{2}mt\right\rfloor}\Big|\Big|P^{(m)}_{\Hy,\left\lfloor \sigma^{2}mt\right\rfloor}\right) \ = \ \lim_{\lambda\nearrow1}\lim_{m\rightarrow\infty}\frac{1-H_{\lambda}\left(P^{(m)}_{\Al,\left\lfloor \sigma^{2}mt\right\rfloor}\Big|\Big|P^{(m)}_{\Hy,\left\lfloor \sigma^{2}mt\right\rfloor}\right)}{\lambda\cdot(1-\lambda)} \ .
\notag
\ee
This will be achieved by evaluating for each $t>0$ the two limits
\be
\lim_{\lambda\nearrow1} \frac{1-
D_{\lambda,t}^{L}
}{\lambda\cdot(1-\lambda)} \qquad \textrm{and} \qquad \lim_{\lambda\nearrow1} \frac{1-
D_{\lambda,t}^{U}
}{\lambda\cdot(1-\lambda)} \label{fo.prthmentlim1c}
\ee
which will turn out to coincide; the involved lower and upper bound $D_{\lambda,t}^{L}$\, , $D_{\lambda,t}^{U}$
defined by \eqref{fo.boulim1} and \eqref{fo.boulim} 
satisfy $\lim_{\lambda\nearrow1}D_{\lambda,t}^{L}=\lim_{\lambda\nearrow1}D_{\lambda,t}^{U}=1$
as an easy consequence of the limits (cf.\ \ref{def.kl-Lambda})
\be
\label{fo.prthmentlim1d}
\lim_{\lambda\nearrow1}\Lambda_{\!\lambda}=\kappa_{\Al} \geq 0 \qquad\textrm{and}\qquad \lim_{\lambda\nearrow1}\kappa_{\lambda}=\kappa_{\Al} \geq 0 \ ,
\ee 
as well as the formulae \eqref{limhilf1}, \eqref{limhilf2} for the case $\kappa_{\Al}=0$.
Accordingly, we compute
\bea\label{fo.prthmentlim2}
&& 
\lim_{\lambda\nearrow1} \frac{1-
D_{\lambda,t}^{L}
}{\lambda\cdot(1-\lambda)} 
\ = \ 
\lim_{\lambda\nearrow1}\frac{-D_{\lambda,t}^{L}}{1-2\lambda} \ \frac{\partial}{\partial\lambda} \Bigg[-\frac{\Lambda_{\!\lambda}-\kappa_{\lambda}}{\sigma^{2}}\cdot\left[\widetilde{X}_{0}-\frac{\eta}{\Lambda_{\!\lambda}}\right]\cdot\left(1-e^{-\Lambda_{\!\lambda}\cdot t}\right)-\frac{\eta}{\sigma^{2}}\cdot\left(\Lambda_{\!\lambda}-\kappa_{\lambda}\right)\cdot t \notag\\
&& \hspace{5.7cm}
+ \ L_{\lambda}^{(1)}(t)\cdot\widetilde{X}_{0}+\frac{\eta}{\sigma^{2}}\cdot L_{\lambda}^{(2)}(t)\Bigg]\notag\\
& & = \lim_{\lambda\nearrow1}\Bigg\{-\frac{\Lambda_{\!\lambda}-\kappa_{\lambda}}{\sigma^{2}}\left[\left(\widetilde{X}_{0}-\frac{\eta}{\Lambda_{\!\lambda}}\right)\cdot te^{-\Lambda_{\!\lambda}\cdot t} \cdot \frac{\partial \, \Lambda_{\!\lambda}}{\partial\lambda}+\left(1-e^{-\Lambda_{\!\lambda}\cdot t}\right)\cdot\frac{\eta}{\Lambda_{\!\lambda}^{2}}\cdot\frac{\partial \, \Lambda_{\!\lambda}}{\partial\lambda}\right]
\notag\\
&& \hspace{1.5cm}
- \, \frac{1}{\sigma^{2}}\cdot\frac{\partial \left(\Lambda_{\!\lambda}-\kappa_{\lambda}\right)}{\partial\lambda}\cdot\left(\widetilde{X}_{0}-\frac{\eta}{\Lambda_{\!\lambda}}\right)\cdot\left(1-e^{-\Lambda_{\!\lambda}\cdot t}\right)-\frac{\eta\, t}{\sigma^{2}}\cdot\frac{\partial \left(\Lambda_{\!\lambda}-\kappa_{\lambda}\right)}{\partial\lambda}
\notag\\
&&\hspace{1.5cm}
+ \, \widetilde{X}_{0} \ \frac{\partial L_{\lambda}^{(1)}(t)}{\partial\lambda}  \ + \ \frac{\eta}{\sigma^{2}} \ \frac{\partial L_{\lambda}^{(2)}(t)}{\partial\lambda}\Bigg\} \ , \qquad \textrm{with}
\eea
\be\label{fo.prthmentlim4} 
\frac{\partial \, \Lambda_{\!\lambda}}{\partial\lambda} \ = \ \frac{\kappa_{\Al}^{2}-\kappa_{\Hy}^{2}}{2\, \Lambda_{\!\lambda}} 
\qquad \textrm{and}\qquad\frac{\partial \, \kappa_{\lambda} }{\partial\lambda} \ = \ \kappa_{\Al}-\kappa_{\Hy} \ . 
\ee
For the case $\kappa_{\Al}>0$, one can combine  
this with \eqref{fo.prthmentlim1d} and 
\eqref{part.all}
to end up with
\be\label{fo.prthmentlim11}
\lim_{\lambda\nearrow1} \frac{1-
D_{\lambda,t}^{L}
}{\lambda\cdot(1-\lambda)} 
\ = \  \frac{\left(\kappa_{\Al}-\kappa_{\Hy}\right)^{2}}{2\sigma^{2}\cdot\kappa_{\Al}}\cdot\left[\left(\widetilde{X}_{0}-\frac{\eta}{\kappa_{\Al}}\right)\cdot\left(1-e^{-\kappa_{\Al}\cdot t}\right)+\eta\cdot t\right].
\ee
For the case $\kappa_{\Al}=0$, we continue the calculation \eqref{fo.prthmentlim2}
by rearranging terms and by employing the formulae  \eqref{part.low1}, \eqref{part.low2}, \eqref{limhilf1}, \eqref{limhilf2}
as well as the obvious 
relation $\frac{1}{\Lambda}-\frac{\Lambda-\kappa_{\lambda}}{\Lambda^{2}}  = \frac{1}{\kappa_{\Hy}}$	to obtain 
\bea\label{fo.prthmentlim9}
& & 
\hspace{-0.7cm} \lim_{\lambda\nearrow1} \frac{1-
D_{\lambda,t}^{L}
}{\lambda\cdot(1-\lambda)} \ = \ \lim_{\lambda\nearrow1}\Bigg\{\frac{\kappa_{\Hy}^{2}\cdot\widetilde{X}_{0}}{2\sigma^{2}}\left[\frac{\Lambda_{\!\lambda}-\kappa_{\lambda}}{\Lambda_{\!\lambda}}\cdot t\cdot e^{-\Lambda_{\!\lambda} t}+\frac{1-e^{-\Lambda_{\!\lambda} t}}{\Lambda_{\!\lambda}}\right]\notag\\
&&  \ + \ \frac{\eta\cdot\kappa_{\Hy}^{2}\cdot t}{2\sigma^{2}}\left[\frac{1}{\Lambda_{\!\lambda}}-\frac{\Lambda_{\!\lambda}-\kappa_{\lambda}}{\Lambda_{\!\lambda}^{2}}+\frac{\Lambda_{\!\lambda}-\kappa_{\lambda}}{\Lambda_{\!\lambda}}\cdot\frac{1-e^{-\Lambda_{\!\lambda} t}}{\Lambda_{\!\lambda}}\right]
- \ \frac{\eta\cdot\kappa_{\Hy}^{2}}{2\sigma^{2}}\cdot\frac{1-e^{-\Lambda_{\!\lambda} t}}{\Lambda_{\!\lambda}}\left[\frac{1}{\Lambda_{\!\lambda}}-\frac{\Lambda_{\!\lambda}-\kappa_{\lambda}}{\Lambda_{\!\lambda}^{2}}\right] 
\notag\\
&&  
\ - \ \frac{\kappa_{\Hy}\cdot\widetilde{X}_{0}}{\sigma^{2}}\left(1-e^{-\Lambda_{\!\lambda} t}\right) \ + \ \frac{\eta\cdot\kappa_{\Hy}}{\sigma^{2}}\left[\frac{1-e^{-\Lambda_{\!\lambda} t}}{\Lambda_{\!\lambda}}-t\right]
 + \ \frac{\partial L_{\lambda}^{(1)}(t)}{\partial\lambda}\cdot\widetilde{X}_{0} \ + \ \frac{\eta}{\sigma^{2}}\cdot\frac{\partial L_{\lambda}^{(2)}(t)}{\partial\lambda}\Bigg\}\notag\\
& & 
 = \frac{\kappa_{\Hy}^{2}\, \widetilde{X}_{0} \, t}{\sigma^{2}}  +  \frac{\eta \, \kappa_{\Hy}^{2}\, t}{2\sigma^{2}}\left[\frac{1}{\kappa_{\Hy}}+t\right]  -  \frac{\eta\, \kappa_{\Hy}\, t}{2\sigma^{2}}  -  \frac{\kappa_{\Hy}^{2}\, 
\widetilde{X}_{0} \, t}{2\sigma^{2}}  -  \frac{\eta\, \kappa_{\Hy}^{2}\, t^{2}}{4\sigma^{2}}
\ = \ \frac{\kappa_{\Hy}^{2}}{2\sigma^{2}}\cdot\left[\frac{\eta}{2}\cdot t^{2}  +  \widetilde{X}_{0}\cdot t\right].
\eea
\noindent
Let us now turn to the second limit \eqref{fo.prthmentlim1c} for which we 
compute analogously to \eqref{fo.prthmentlim2} 
\bea\label{fo.prthmentlim10}
&& \lim_{\lambda\nearrow1} \frac{1-
D_{\lambda,t}^{U}
}{\lambda\cdot(1-\lambda)} 
\ = \  \lim_{\lambda\nearrow1}\frac{-D_{\lambda,t}^{U}}{1-2\lambda} \ \frac{\partial}{\partial\lambda}\Bigg[-\frac{\Lambda_{\!\lambda}-\kappa_{\lambda}}{\sigma^{2}}\cdot\left[\widetilde{X}_{0}-\frac{\eta}{\frac{1}{2}(\Lambda_{\!\lambda}+\kappa_{\lambda})}\right]\cdot\left(1-e^{-\frac{1}{2}(\Lambda_{\!\lambda}+\kappa_{\lambda})\cdot t}\right)
\notag \\
&& \hspace{5.7cm} -\frac{\eta}{\sigma^{2}}\cdot\left(\Lambda_{\!\lambda}-\kappa_{\lambda}\right)\cdot t
- U_{\lambda}^{(1)}(t)\cdot\widetilde{X}_{0} - \frac{\eta}{\sigma^{2}}\cdot U_{\lambda}^{(2)}(t)\Bigg]\notag\\
& &= \lim_{\lambda\nearrow1}\Bigg\{-\frac{\Lambda_{\!\lambda}-\kappa_{\lambda}}{\sigma^{2}}\Bigg[\left(\widetilde{X}_{0}-\frac{\eta}{\frac{1}{2}(\Lambda_{\!\lambda}+\kappa_{\lambda})}\right)\cdot \frac{t}{2} \cdot e^{-\frac{1}{2}(\Lambda_{\!\lambda}+\kappa_{\lambda})\cdot t} \ \frac{\partial \left(\Lambda_{\!\lambda}+\kappa_{\lambda}\right)}{\partial\lambda}\notag\\
&& \hspace{2.8cm}+ \ \left(1-e^{-\frac{1}{2}(\Lambda_{\!\lambda}+\kappa_{\lambda})\cdot t}\right)\cdot\frac{2\cdot\eta}{(\Lambda_{\!\lambda}+\kappa_{\lambda})^{2}}\cdot\frac{\partial (\Lambda_{\!\lambda}+\kappa_{\lambda})}{\partial\lambda}\Bigg]\notag\\
&& \hspace{1.0cm}- \ \frac{1}{\sigma^{2}}\cdot\frac{\partial \left(\Lambda_{\!\lambda}-\kappa_{\lambda}\right)}{\partial\lambda}\cdot\left(\widetilde{X}_{0}-\frac{\eta}{\frac{1}{2}(\Lambda_{\!\lambda}+\kappa_{\lambda})}\right)\cdot\left(1-e^{-\frac{1}{2}(\Lambda_{\!\lambda}+\kappa_{\lambda})\cdot t}\right)-\frac{\eta\, t}{\sigma^{2}}\cdot\frac{\partial \left(\Lambda_{\!\lambda}-\kappa_{\lambda}\right)}{\partial\lambda}\notag\\
&&\hspace{1.0cm}- \ \frac{\partial \, U_{\lambda}^{(1)}(t)}{\partial\lambda}\cdot\widetilde{X}_{0} \ - \ \frac{\eta}{\sigma^{2}}\frac{\partial \, U_{\lambda}^{(2)}(t)}{\partial\lambda}\Bigg\} \ . 
\eea
\noindent
For the case $\kappa_{\Al}>0$, one can combine  
this with \eqref{fo.prthmentlim1d}, \eqref{fo.prthmentlim4} and 
\eqref{part.all} to end up with
\be\label{fo.prthmentlim12}
\lim_{\lambda\nearrow1} \frac{1-
D_{\lambda,t}^{U}
}{\lambda\cdot(1-\lambda)} \ = \  \frac{\left(\kappa_{\Al}-\kappa_{\Hy}\right)^{2}}{2\sigma^{2}\cdot\kappa_{\Al}}\cdot\left[\left(\widetilde{X}_{0}-\frac{\eta}{\kappa_{\Al}}\right)\cdot\left(1-e^{-\kappa_{\Al}\cdot t}\right)+\eta\cdot t\right].
\ee
For the case $\kappa_{\Al}=0$, we continue the calculation of \eqref{fo.prthmentlim10}
by rearranging terms and by employing the formulae \eqref{part.up}, \eqref{limhilf1}, \eqref{limhilf2}
as well as the obvious
relation 
$\lim_{\lambda\nearrow1}\frac{1}{\Lambda_{\!\lambda}}-\frac{\Lambda_{\!\lambda}-\kappa_{\lambda}}{\Lambda_{\!\lambda}(\Lambda_{\!\lambda}+\kappa_{\lambda})} \ = \ \frac{2}{\kappa_{\Hy}}$	
to obtain 
\bea\label{fo.prthmentlim13}
&& \lim_{\lambda\nearrow1} \frac{1-
D_{\lambda,t}^{U}
}{\lambda\cdot(1-\lambda)} 
\ = \ 
\lim_{\lambda\nearrow1}\Bigg\{\frac{t\cdot\widetilde{X}_{0}}{4\sigma^{2}}\cdot\frac{\Lambda_{\!\lambda}-\kappa_{\lambda}}{\Lambda_{\!\lambda}}\cdot e^{-\frac{1}{2}\left(\Lambda_{\!\lambda}+\kappa_{\lambda}\right)\cdot t}\left(\kappa_{\Hy}^{2}+2\Lambda_{\!\lambda}\kappa_{\Hy}\right) 
\notag \\
&&  
 + \ \frac{\widetilde{X}_{0}}{2\sigma^{2}}\cdot\frac{1-e^{-\frac{1}{2}(\Lambda_{\!\lambda}+\kappa_{\lambda})\cdot t}}{\Lambda_{\!\lambda}}\left(\kappa_{\Hy}^{2}-2\Lambda_{\!\lambda}\kappa_{\Hy}\right)
- \ \frac{\eta\cdot t}{\sigma^{2}}\bigg[
\kappa_{\Hy}\left(1+e^{-\frac{1}{2}\left(\Lambda_{\!\lambda}+\kappa_{\lambda}\right)\cdot t}\ \frac{\Lambda_{\!\lambda}-\kappa_{\lambda}}{\Lambda_{\!\lambda}+\kappa_{\lambda}}\right)
\notag\\
&& 
- \ \frac{\kappa_{\Hy}^{2}}{2}\cdot\left( \frac{1}{\Lambda_{\!\lambda}}-\frac{\Lambda_{\!\lambda}-\kappa_{\lambda}}{\Lambda_{\!\lambda}(\Lambda_{\!\lambda}+\kappa_{\lambda})} 
+  \frac{\Lambda_{\!\lambda}-\kappa_{\lambda}}{\Lambda_{\!\lambda}+\kappa_{\lambda}} \cdot 
\frac{1-e^{-\frac{1}{2}\left(\Lambda_{\!\lambda}+\kappa_{\lambda}\right)\cdot t}}{\Lambda_{\!\lambda}}
\right) \bigg] \notag \\
&& 
+ \ \frac{2\eta}{\sigma^{2}}\cdot\frac{1-e^{-\frac{1}{2}\left(\Lambda_{\!\lambda}+\kappa_{\lambda}\right)\cdot t}}{\Lambda_{\!\lambda}+\kappa_{\lambda}}\left[\kappa_{\Hy}\left(1+\frac{\Lambda_{\!\lambda}-\kappa_{\lambda}}{\Lambda_{\!\lambda}+\kappa_{\lambda}}\right)-\frac{\kappa_{\Hy}^{2}}{2}\left(\frac{1}{\Lambda_{\!\lambda}}-\frac{\Lambda_{\!\lambda}-\kappa_{\lambda}}{\Lambda_{\!\lambda}(\Lambda_{\!\lambda}+\kappa_{\lambda})}\right)\right]\notag\\
&& 
 - \ \frac{\partial \, U_{\lambda}^{(1)}(t)}{\partial\lambda}\cdot\widetilde{X}_{0} \ - \ \frac{\eta}{\sigma^{2}}\frac{\partial \, U_{\lambda}^{(2)}(t)}{\partial\lambda}\Bigg\}\notag\\
&& \hspace{0.2cm}=  \frac{\kappa_{\Hy}^{2}\, t\, \widetilde{X}_{0}}{4\sigma^{2}}  +  \frac{\kappa_{\Hy}^{2}\, t\, \widetilde{X}_{0}}{4\sigma^{2}}  -  \frac{\eta\, t}{\sigma^{2}}\left[
2\kappa_{\Hy} - \kappa_{\Hy} - \frac{\kappa_{\Hy}^2 \, t}{4}
\right]  +  \frac{\eta\, t}{\sigma^{2}}\left[2\kappa_{\Hy}-\kappa_{\Hy}\right]
\ = \ \frac{\kappa_{\Hy}^{2}}{2\sigma^{2}}\left[\frac{\eta}{2}\cdot t^{2}+\widetilde{X}_{0}\cdot t\right].
\eea
Since 
\eqref{fo.prthmentlim11} coincides with \eqref{fo.prthmentlim12} and
\eqref{fo.prthmentlim9} coincides with \eqref{fo.prthmentlim13}, we have finished the proof.
\qed\\

\end{appendix}

\noindent
\textbf{Acknowledgements.} \ 
We are very grateful to the ``Studienstiftung
des deutschen Volkes'' for the generous grant of the first author. 
Furthermore, we 
would also like to thank Andreas Greven 
for some useful remarks.

\end{document}